\RequirePackage{ifpdf}
\ifpdf 
\documentclass[pdftex]{sigma}
\else
\documentclass{sigma}
\fi

\usepackage{amssymb,bbm}    
\usepackage{array}
\usepackage[arrow,curve,frame,matrix,tips]{xy} 
\def\mymatrix#1#2#3{\begin{equation*}\label{eq:#1}
    \xy 0*!C\xybox{\xymatrix#2{#3}} \endxy\end{equation*}}
\usepackage{leftidx} 

\def\C{{\mathcal C}}
\def\D{{\mathcal D}}
\def\E{{\mathcal E}}
\def\H{{\mathcal H}}
\def\I{{\mathcal I}}
\def\Nc#1{\mathcal{N}_{#1}}
\def\Ncp#1#2{\mathcal{N}_{#1}^{(#2)}}
\def\Rc#1{\mathcal{R}_{#1}}
\def\Rcp#1#2{\mathcal{R}_{#1}^{(#2)}}
\def\U{{\mathcal U}}
\def\V{{\mathcal V}}
\def\X{{\mathcal X}}
\def\aq#1#2{\alpha_{#1}^{(#2)}}
\def\bq#1#2{\beta_{#1}^{(#2)}}
\def\fprod#1#2{\underset{#2}{#1}}
\def\je#1{J_{#1}\pi}
\def\mbl#1#2#3#4#5{\leftidx{_{#2}^{#3}}{{#1}}{_{#4}^{#5}}}
\def\sv{\mathbf{s}}
\def\uv{\mathbf{u}}
\def\xv{\mathbf{x}}

\DeclareMathOperator{\cl}{cls}
\DeclareMathOperator{\im}{im}
\DeclareMathOperator{\rank}{rank}

\begin{document}

\allowdisplaybreaks

\renewcommand{\thefootnote}{$\star$}

\renewcommand{\PaperNumber}{092}

\FirstPageHeading

\ShortArticleName{Existence and Construction of Vessiot Connections}

\ArticleName{Existence and Construction of Vessiot Connections\footnote{This paper is a
contribution to the Special Issue ``\'Elie Cartan and Dif\/ferential Geometry''. The
full collection is available at
\href{http://www.emis.de/journals/SIGMA/Cartan.html}{http://www.emis.de/journals/SIGMA/Cartan.html}}}

\Author{Dirk FESSER~$^\dag$ and Werner M. SEILER~$^\ddag$}

\AuthorNameForHeading{D. Fesser and W.M. Seiler}

\Address{$^\dag$~IWR, Universit\"at Heidelberg, INF 368, 69120 Heidelberg,
  Germany}

\EmailD{\href{mailto:dirk.fesser@iwr.uni-heidelberg.de}{dirk.fesser@iwr.uni-heidelberg.de}}

\Address{$^\ddag$~AG ``Computational Mathematics'', Universit\"at Kassel,
  34132 Kassel, Germany}
\EmailD{\href{mailto:seiler@mathematik.uni-kassel.de}{seiler@mathematik.uni-kassel.de}}
\URLaddressD{\url{http://www.mathematik.uni-kassel.de/~seiler/}}

\ArticleDates{Received May 05, 2009, in f\/inal form September 14, 2009;  Published online September 25, 2009}

\Abstract{A rigorous formulation of Vessiot's vector f\/ield approach to the analysis of
  ge\-ne\-ral systems of partial dif\/ferential equations is provided.  It is shown
  that this approach is equivalent to the formal theory of dif\/ferential
  equations and that it can be carried through if, and only if, the given
  system is involutive.  As a by-product, we provide a novel characterisation
  of transversal integral elements via the contact map.}


\Classification{35A07; 35A30; 35N99; 58A20}

\tableofcontents

\renewcommand{\thefootnote}{\arabic{footnote}}
\setcounter{footnote}{0}

\section{Introduction}

Constructing solutions for (systems of) partial dif\/ferential equations is
obviously dif\/f\/icult~-- in particular for non-linear systems.  \'Elie Cartan
\cite{Cartan} proposed to construct f\/irst inf\/initesimal solutions or integral
elements.  These are possible tangent spaces to (prolonged) solutions.  Thus
they always lead to a linearisation of the problem and their explicit
construction requires essentially only straightforward linear algebra.  In the
Cartan--K\"ahler theory \cite{BCGGG:ExteriorDS, IveyLandsberg:Cartan},
dif\/ferential equations are represented by exterior dif\/ferential systems and
integral elements consist of tangent vectors pointwise annihilated by
dif\/ferential forms.

Vessiot \cite{Vessiot:Integration} proposed in the 1920s a dual approach which
does not require the use of exterior dif\/ferential systems.  Instead of
individual integral elements, it always considers distributions of them
generated by vector f\/ields and their Lie brackets replace the exterior
derivatives of dif\/ferential forms.  This approach takes an intermediate
position between the formal theory of dif\/ferential equations
\cite{KrasilshchikLychaginVinogradov, Pommaret:SystemsOfPDEs, Seiler:Buch} and
the Cartan--K\"ahler theory of exterior dif\/ferential systems.  Thus it allows
for the transfer of many techniques from the latter to the former one,
although this point will not be studied here.

Vessiot's approach may be considered a generalisation of the Frobenius
theorem.  Indeed, if one applies his theory to a dif\/ferential equation of
f\/inite type, then one obtains an involutive distribution such that its
integral manifolds are in a one-to-one correspondence with the smooth
solutions of the equation.  For more general equations, Vessiot proposed to
``cover'' the equation with inf\/initely many involutive distributions such that
any smooth solution corresponds to an integral manifold of at least one of
them.

Vessiot's theory has not attracted much attention: presentations in a more
modern language are contained in
\cite{Fackerell:VessiotsVectorFieldFormulation,
  Stormark:LiesStructuralApproach}; applications have mainly appeared in the
context of the Darboux method for solving hyperbolic equations, see for
example \cite{Vassiliou:VessiotStructure}.  While a number of textbooks
provide a~very rigorous analysis of the Cartan--K\"ahler theory, the above
mentioned references (including Vessiot's original work
\cite{Vessiot:Integration}) are somewhat lacking in this respect.  In
particular, the question of under what assumptions Vessiot's construction
succeeds has been ignored.

The purpose of the present article is to close this gap and at the same time
to relate the Vessiot theory with the key concepts of the formal theory like
formal integrability and involution (we will not develop it as a dual form of
the Cartan--K\"ahler theory, but from scratch within the formal theory).  We
will show that Vessiot's construction succeeds if, and only if, it is applied
to an involutive system of dif\/ferential equations.  This result is of course
not surprising, given the well-known fact that the formal theory and the
Cartan--K\"ahler theory are equivalent.  However, to our knowledge an explicit
proof has never been given.  As a by-product, we will provide a~new
characterisation of integral elements based on the contact map, making also
the relations between the formal theory and the Cartan--K\"ahler theory more
transparent.  Furthermore, we simplify the construction of the integral
distributions.  Up to now, quadratic equations had to be considered (under
some assumptions, their solution can be obtained via a sequence of linear
systems.  We will show how the natural geometry of the jet bundle hierarchy
can be exploited for always obtaining a linear system of equations.

This article contains the main results of the f\/irst author's doctoral thesis
\cite{Fesser:2008} (a short summary of the results has already appeared in
\cite{FesserSeiler:Daresbury}).  It is organised as follows.  The next three
sections recall the needed elements from the formal theory of dif\/ferential
equations and provide our new characterisation of transversal integral
elements.  The following two sections introduce the key concepts of Vessiot's
approach: the Vessiot distribution, integral distributions and f\/lat Vessiot
connections.  Two further sections discuss existence theorems for the latter
two based on a step-by-step approach already proposed by Vessiot.  As the
proofs of some results are fairly technical, the main text contains only an
outline of the underlying ideas and full details are given in three
appendices.  The Einstein convention is sometimes used to indicate summation.

\section{The contact structure}

Before we outline the formal theory of partial dif\/ferential equations, we
brief\/ly review its underlying geometry: the jet bundle and its contact
structure. Many dif\/ferent ways exist to introduce these geometric
constructions, see for example \cite{GolubitskyGuillemin:Stable Mappings,
  Olver:LieGroups, Saunders:JetBundles}.  Furthermore, they are discussed in
any book on the formal theory (see the references in the next section).

Let $\pi:\E\rightarrow\X$ be a smooth f\/ibred manifold.  We call coordinates
$\xv=(x^i \colon 1 \le i \le n)$ of~$\X$ \emph{independent variables} and
f\/ibre coordinates $\uv=(u^\alpha \colon 1 \le \alpha \le m)$ in $\E$
\emph{dependent variables}.  Sections $\sigma:\X\rightarrow\E$ correspond
locally to functions $\uv=\sv(\xv)$.  We will use throughout a ``global''
notation in order to avoid the introduction of many local neighbourhoods even
though we mostly consider local sections.

Derivatives are written in the form $u^\alpha_\mu = \partial^{|\mu|}u^\alpha
/\partial x_1^{\mu_1}\cdots\partial x_n^{\mu_n}$ where
$\mu=(\mu_1,\dots,\mu_n)$ is a multi-index.  The set of derivatives
$u^\alpha_\mu$ up to order~$q$ is denoted by $\uv^{(q)}$; it def\/ines a local
coordinate system for the $q$-th order jet bundle~$\je{q}$, which may be
regarded as the space of truncated Taylor expansions of functions $\sv$.

The hierarchy of jet bundles $\je{q}$ with $q=0,1,2,\dots$ possesses many
natural f\/ibrations which correspond to ``forgetting'' higher-order
derivatives.  For us particularly important are
$\pi_{q-1}^q:\je{q}\to\je{q-1}$ and $\pi^q:\je{q}\to\X$.  To each section
$\sigma:\X\rightarrow\E$, locally def\/ined by
$\sigma(\xv)=\bigl(\xv,\sv(\xv)\bigr)$, we may associate its
\emph{prolongation} $j_q\sigma:\X\rightarrow\je{q}$, a section of the
f\/ibration $\pi^q$ locally given by $j_q\sigma(\xv)=\bigl(\xv,\sv(\xv),
\partial_{\xv}\sv(\xv),\partial_{\xv\xv}\sv(\xv),\dots\bigr)$.

The geometry of the jet bundle $\je{q}$ is to a large extent determined by its
\emph{contact structure}.  It can be introduced in various ways.  For our
purposes, three dif\/ferent approaches are convenient.  First, we adopt the
\emph{contact codistribution} $\C_q^{0}\subseteq T^*(\je{q})$; it consists of
all one-forms such that their pull-back by a prolonged section vanishes.
Locally, it is spanned by the \emph{contact forms}
\begin{equation}
	\omega^\alpha_\mu
	=
	d u^\alpha_\mu
	-
	\sum_{i=1}^nu^\alpha_{\mu+1_i}dx^i ,
		\qquad
		0\le|\mu|<q ,
		\quad
		1 \le \alpha \le m .
\nonumber
\end{equation}
Dually, we may consider the \emph{contact distribution} $\C_q\subseteq
T(\je{q})$ consisting of all vector f\/ields annihilated by $\C_q^{0}$.  A
straightforward calculation shows that it is generated by the \emph{contact
  fields}
\begin{gather}
		C_i^{(q)}
		 =
		\partial_i
		+
		\sum_{\alpha=1}^m
			\sum_{0 \le |\mu| < q}
				u^\alpha_{\mu+1_i}\partial_{u^\alpha_\mu} ,
		\qquad
			 1\leq i\leq n ,\nonumber
		\\
		C^\mu_\alpha
		 =
		\partial_{u^\alpha_\mu} ,
			\qquad |\mu|=q ,
			\quad
			1 \le \alpha \le m .\label{gl:ContactFields}
        \end{gather}
Note that the latter f\/ields, $C^\mu_\alpha$, span the vertical bundle
$V\pi^q_{q-1}$ of the f\/ibration $\pi^q_{q-1}$.  Thus the contact distribution
can be split into $\C_q=V\pi^q_{q-1}\oplus\H$.  Here the complement $\H$ is an
$n$-dimensional transversal subbundle of $T(\je{q})$ and obviously not
uniquely determined (though any local coordinate chart induces via the span of
the vectors $C_i^{(q)}$ one possible choice).  Any such complement $\H$ may be
considered the horizontal bundle of a connection on the f\/ibred manifold
$\pi^q:\je{q}\to\X$ (but not for the f\/ibration $\pi^q_{q-1}$).  Following
Fackerell \cite{Fackerell:VessiotsVectorFieldFormulation}, we call any
connection on~$\pi^q$ the horizontal bundle of which consists of contact
f\/ields a \emph{Vessiot connection} (in the literature the terminology
\emph{Cartan connection} is also common, see for example
\cite{Lychagin:1995}).

For later use, we note the structure equations of the contact distribution.
The only non-vanishing Lie brackets of the vector f\/ields
(\ref{gl:ContactFields}) are
\begin{equation}\label{gl:LieOrdnung1kleiner}
	\bigl[C^{\nu + 1_i}_\alpha , C^{(q)}_i\bigr]
	=
	\partial_{u^\alpha_\nu}
	 ,
	\qquad |\nu| = q-1  .
\end{equation}
Note that this observation implies that the vertical bundle $V\pi^q_{q-1}$ is
involutive.

As a third approach to the contact structure we consider, following Modugno
\cite{Modugno:1999}, the \emph{contact map} (see also
\cite{Fesser:2008,Seiler:Buch}).  It is the unique map
$\Gamma_q:\je{q}\fprod{\times}{\X}T\X\rightarrow T(\je{q-1})$ such that the
diagram
\mymatrix{contact}{}{
  *!<0pt,3pt>{\je{q}\fprod{\times}{\X}T\X} \ar[rr]^{\Gamma_q} &&
  {T(\je{q-1})}\\
  & {T\X} \ar[ul]^{((j_q\sigma)\circ\tau_\X)\times \mathrm{id}_{T\X}}
  \ar[ur]_{T(j_{q-1}\sigma)} }
commutes for any section $\sigma$.  Because of its linearity over
$\pi^q_{q-1}$, we may also consider it a map $\Gamma_q:\je{q}\rightarrow
T^*\X\fprod{\otimes}{\je{q-1}}T(\je{q-1})$ with the local coordinate form
\begin{equation}\label{gl:gammaq}
        \Gamma_q: \ (\xv,\uv^{(q)})\mapsto\bigl(\xv,\uv^{(q-1)};
                d x^i\otimes
                (\partial_{x^i}+u^\alpha_{\mu+1_i}\partial_{u^\alpha_\mu})\bigr) .
\end{equation}
Obviously, $\Gamma_q(\rho,\partial_{x^i})=T_\rho\pi^q_{q-1}(C_i^{(q)})$ and
hence $(\C_q)_\rho=(T_\rho\pi^q_{q-1})^{-1}\bigl(\im{\Gamma_q(\rho)}\bigr)$
for any point $\rho\in\je{q}$.  Note that all vectors in the image of
$\Gamma_q(\rho)$ are transversal to the f\/ibration $\pi^{q-1}_{q-2}$.

One of the main applications of the contact structure is given by the
following proposition (for a proof, see \cite[Proposition 2.1.6]{Fesser:2008}
or \cite[Proposition 2.2.7]{Seiler:Buch}).  It characterises those sections of
the jet bundle $\pi^{q}\colon\je{q}\to\X$ which are prolongations of sections
of the underlying f\/ibred manifold $\pi \colon \E \to \X$.

\begin{proposition}\label{prop:contactmapq}
A section\/ $\gamma:\X\to\je{q}$ is of the form\/ $\gamma=j_q\sigma$
for a section\/ $\sigma:\X\to\E$ if, and only if,\/
$\im{\Gamma_q\bigl(\gamma(x)\bigr)}=
T_{\gamma(x)}\pi^q_{q-1}\bigl(T_{\gamma(x)}\im{\gamma}\bigr)$ for
all points\/ $x\in\X$ where\/ $\gamma$ is defined.
\end{proposition}
Thus for any section $\sigma:\X\to\E$ the equality
$\im{\Gamma_{q+1}\bigl(j_{q+1}\sigma(x)\bigr)}=\im{T_{x}(j_q\sigma)}$ holds
and we may say that knowing the $(q+1)$-jet $j_{q+1}\sigma(x)$ of a section
$\sigma$ at some $x\in\X$ is equivalent to knowing its $q$-jet
$\rho=j_q\sigma(x)$ at $x$ as well as the tangent space
$T_\rho(\im{j_q\sigma})$ at this point.  This observation will later be the
key for the Vessiot theory.

\section{The formal theory of dif\/ferential equations}

We are now going to outline the formal theory of partial dif\/ferential
equations to introduce the basic notation.  Our presentation follows
\cite{Seiler:Buch}; other general references are
\cite{KrasilshchikLychaginVinogradov, KruglikovLychagin:2008,
  Pommaret:SystemsOfPDEs}.

\begin{definition}
  A \emph{differential equation} of order $q$ is a f\/ibred submanifold
  $\Rc{q}\subseteq\je{q}$ locally described as the zero set of some smooth
  functions on $\je{q}$:
  \begin{equation}\label{gl:system}
        \Rc{q}\,\colon\,
                \left\{
                  \begin{array}{l}
                    \mathnormal\Phi^\tau\big(\xv,\uv^{(q)}\big)=0 ,\\
                    (\tau=1,\dots,t) .
                  \end{array}
                \right.
  \end{equation}
  Note that we do not distinguish between scalar equations and systems.
\end{definition}

We denote by $\iota:\Rc{q}\hookrightarrow\je{q}$ the canonical inclusion map.
Dif\/ferentiating every equation in the local representation (\ref{gl:system})
leads to the \emph{prolonged equation} $\Rc{q+1}\subseteq\je{q+1}$ def\/ined by
the equations $\mathnormal\Phi^\tau=0$ and $D_i\mathnormal\Phi^\tau=0$ where
the formal derivative $D_i$ is given by
\begin{equation}\label{gl:formderiv}
	D_i\mathnormal\Phi^\tau\big(\xv,\uv^{(q+1)}\big)
	=
    \frac{\partial \mathnormal\Phi^\tau}{\partial x^i}\big(\xv,\uv^{(q)}\big)
    +
    \sum_{0 \le |\mu| \le q}
		\sum_{\alpha=1}^m
			\frac{\partial \mathnormal\Phi^\tau}{\partial u^\alpha_\mu}
				\big(\xv,\uv^{(q)}\big)u^\alpha_{\mu+1_i}.
\end{equation}
Iteration of this process gives the higher prolongations
$\Rc{q+r}\subseteq\je{q+r}$.  A subsequent projection leads to
$\Rcp{q}{1}=\pi^{q+1}_q(\Rc{q+1})\subseteq\Rc{q}$, which is a proper
submanifold whenever integrability conditions appear.

\begin{definition}
  A dif\/ferential equation $\Rc{q}$ is \emph{formally integrable} if at any
  prolongation order $r>0$ the equality $\Rcp{q+r}{1}=\Rc{q+r}$ holds.
\end{definition}

In local coordinates, the following def\/inition coincides with the usual notion
of a solution.

\begin{definition}\label{def:Solution}
  A \emph{solution} is a section $\sigma:\X\to\E$ such that its prolongation
  satisf\/ies $\im{j_q\sigma}\subseteq\Rc{q}$.
\end{definition}

For formally integrable equations it is straightforward to construct order by
order formal power series solutions.  Otherwise it is hard to f\/ind solutions.
A constitutive insight of Cartan was to introduce \emph{infinitesimal solutions}
or \emph{integral elements} at a point $\rho\in\Rc{q}$ as
subspaces $\U_\rho\subseteq T_\rho\Rc{q}$ which are potentially part of the
tangent space of a prolonged solution.
\begin{definition}\label{def:integralelement}
  Let\/ $\Rc{q}\subseteq\je{q}$ be a dif\/ferential equation and\/
  $\iota:\Rc{q}\rightarrow\je{q}$ the canonical inclusion map. Let\/
  $\I[\Rc{q}]=\langle{\iota^*\C_q^0}\rangle_{\mathrm{dif\/f}}$ be the
  dif\/ferential ideal generated by the pull-back of the contact codistribution
  on\/ $\Rc{q}$ (algebraically, $\I[\Rc{q}]$ is then spanned by a basis of
  $\iota^*\C_q^0$ and the exterior derivatives of the forms in this basis).  A
  linear subspace\/ $\U_\rho\subseteq T_\rho\Rc{q}$ is an\/ \emph{integral
    element} at the point\/ $\rho\in\Rc{q}$, if all forms in\/
  $(\I[\Rc{q}])_\rho$ vanish on it.
\end{definition}
The following result provides an alternative characterisation of
\emph{transversal} integral elements via the contact map.  It requires that
the projection $\pi^{q+1}_q:\Rc{q+1}\rightarrow\Rc{q}$ is surjective.
\begin{proposition}\label{prop:integralelement}
  Let\/ $\Rc{q}$ be a differential equation such that\/ $\Rcp{q}{1}=\Rc{q}$.
  A linear subspace\/ $\U_\rho\subseteq T_\rho\Rc{q}$ such that\/ $T_\rho\iota
  (\U_\rho)$ lies transversal to the fibration\/ $\pi^q_{q-1}$ is an\/
  integral element at the point\/ $\rho\in\Rc{q}$ if, and only if,\/ a point\/
  $\hat\rho\in\Rc{q+1}$ exists on the prolonged equation\/ $\Rc{q+1}$ such
  that\/ $\pi^{q+1}_q(\hat\rho)=\rho$ and\/
  $T_\rho\iota(\U_\rho)\subseteq\im{\Gamma_{q+1}(\hat\rho)}$.
\end{proposition}
\begin{proof}
  Assume f\/irst that $\U_\rho$ satisf\/ies the given conditions.  It follows
  immediately from the coordinate form of the contact map that then f\/irstly
  $T_\rho\iota(\U_\rho)$ is transversal to $\pi^q_{q-1}$ and secondly that
  every one-form $\omega\in\iota^*\C_q^0$ vanishes on $\U_\rho$, as
  $\im{\Gamma_{q+1}(\hat\rho)}\subset(\C_q)_\rho$.  Thus there only remains to
  show that the same is true for the two-forms $d\omega\in\iota^*(d\C_q^0)$.

  Choose a section $\gamma:\Rc{q}\rightarrow\Rc{q+1}$ such that
  $\gamma(\rho)=\hat\rho$ and def\/ine a distribution $\D$ of rank $n$ on
  $\Rc{q}$ by setting
  $T\iota(\D_{\tilde\rho})=\im{\Gamma_{q+1}\bigl(\gamma(\tilde\rho)\bigr)}$
  for any point $\tilde\rho\in\Rc{q}$.  Obviously, by construction
  $\U_\rho\subseteq\D_\rho$.  It follows now from the coordinate form
  (\ref{gl:gammaq}) of the contact map that locally the distribution $\D$ is
  spanned by $n$ vector f\/ields $X_i$ such that
  $\iota_*X_i=C_i^{(q)}+\gamma^\alpha_{\mu+1_i}C_\alpha^\mu$ where the
  coef\/f\/icients $\gamma^\alpha_\nu$ are the highest-order components of the
  section $\gamma$.  Thus the commutator of two such vector f\/ields satisf\/ies
  \begin{equation*}
        \iota_*\bigl([X_i,X_j]\bigr)
        =
        \bigl(C_i^{(q)}(\gamma^\alpha_{\mu+1_j})
        -
        C_j^{(q)}(\gamma^\alpha_{\mu+1_i})\bigr)C_\alpha^\mu
        +
        \gamma^\alpha_{\mu+1_j}[C_i^{(q)},C_\alpha^\mu]
        -
        \gamma^\alpha_{\mu+1_i}[C_j^{(q)},C_\alpha^\mu] .
  \end{equation*}
  The commutators on the right hand side vanish whenever $\mu_i=0$ or $\mu_j=0$,
  respectively.  Otherwise we obtain $-\partial_{u^\alpha_{\mu-1_i}}$ and
  $-\partial_{u^\alpha_{\mu-1_j}}$, respectively.  But this fact implies that
  the two sums on the right hand side cancel each other and we f\/ind that
  $\iota_*\bigl([X_i,X_j]\bigr)\in\C_q$.  Thus we f\/ind for any contact form
  $\omega\in\C_q^0$ that
  \begin{gather*}
        \iota^*(d\omega)(X_i,X_j)
         =
        d\omega(\iota_*X_i,\iota_*X_j)
        =\iota_*X_i\bigl(\omega(\iota_*X_j)\bigr)
        -
        \iota_*X_j\bigl(\omega(\iota_*X_i)\bigr)
        +
        \omega\bigl(\iota_*([X_i,X_j])\bigr) .
  \end{gather*}
  Each summand in the last expression vanishes, as all appearing f\/ields are
  contact f\/ields.  Hence any form $\omega\in\iota^*(d\C_q^0)$ vanishes on $\D$
  and in particular on $\U_\rho\subseteq\D_\rho$.

  For the converse, note that any transversal integral element
  $\U_\rho\subseteq T_\rho\Rc{q}$ is spanned by linear combinations of vectors
  $v_i$ such that $T_\rho\iota(v_i)=C_i^{(q)}|_{\rho}+
  \gamma^\alpha_{\mu,i}C_\alpha^\mu|_{\rho}$ where $\gamma^\alpha_{\mu,i}$ are
  real coef\/f\/icients.  Now consider a contact form $\omega^\alpha_\nu$ with
  $|\nu|=q-1$.  Then $d\omega^\alpha_\nu=d x^i\wedge d u^\alpha_{\nu+1_i}$.
  Evaluating the condition $\iota^*(d\omega^\alpha_\nu)|_{\rho}(v_i,v_j)=
  d\omega\bigl(T_\rho\iota(v_i),T_\rho\iota(v_j)\bigr)=0$ yields the equation
  $\gamma^\alpha_{\nu+1_i,j}=\gamma^\alpha_{\nu+1_j,i}$.  Hence the
  coef\/f\/icients are of the form $\gamma^\alpha_{\mu,i}=\gamma^\alpha_{\mu+1_i}$
  and a section $\sigma$ exists such that $\rho=j_q\sigma(x)$ and
  $T_\rho(\im{j_q\sigma})$ is spanned by the vectors
  $T_\rho\iota(v_1),\dots,T_\rho\iota(v_n)$.  This observation implies that
  $\U_\rho$ satisf\/ies the given conditions.
\end{proof}

For many purposes the purely geometric notion of formal integrability is not
suf\/f\/icient, and one needs the stronger algebraic concept of involution.  This
concerns in particular the derivation of uniqueness results but also the
numerical integration of overdetermined systems \cite{Seiler:CompletionSemi}.
An intrinsic def\/inition of involution is possible using the Spencer cohomology
(see for example \cite{Seiler:Spencer} and references therein for a
discussion).  We apply here a simpler approach requiring that one works in
``good'', more precisely:  $\delta$-regular, coordinates $\xv$.  This
assumption represents a mild restriction, as generic coordinates are
$\delta$-regular and it is possible to construct systematically ``good''
coordinates -- see \cite{HausdorfSeiler:geocompl}.  Furthermore, it will turn
out that the use of $\delta$-regular coordinates is essential for Vessiot's
approach.

\begin{definition}
  The (\emph{geometric$)$ symbol} of a dif\/ferential equation $\Rc{q}$ is
  $\Nc{q}=V\pi^q_{q-1}|_{\Rc{q}}\cap T\Rc{q}$.
\end{definition}

Thus, the symbol is the vertical part of the tangent space to $\Rc{q}$.
Locally, $\Nc{q}$ consists of all vertical vector f\/ields
$\sum_{\alpha=1}^{m}\sum_{|\mu|=q}v^{\alpha}_{\mu}\partial_{u^{\alpha}_{\mu}}$
where the coef\/f\/icients $v^{\alpha}_{\mu}$ satisfy the following linear system
of algebraic equations:
\begin{equation}\label{gl:Symbol}
        \Nc{q} \colon
        \left\{\sum_{\alpha=1}^m\sum_{|\mu|=q}
                \left(\frac{\partial\mathnormal\Phi^\tau}{\partial u^\alpha_\mu}
                \right)
                        v^\alpha_\mu=0 .
        \right.
\end{equation}
The matrix of this system is called the \emph{symbol matrix} $M_q$.  The
\emph{prolonged symbols} $\Nc{q+r}$ are the symbols of the prolonged equations
$\Rc{q+r}$ with corresponding symbol matrices $M_{q+r}$.

The \emph{class} of a multi-index $\mu=(\mu_1,\dots,\mu_n)$, denoted $\cl\mu$,
is the smallest~$k$ for which $\mu_k$ is dif\/ferent from zero.  The columns of
the symbol matrix (\ref{gl:Symbol}) are labelled by the $v^\alpha_\mu$.  We
order them as follows.  Let $\alpha$ and $\beta$ denote indices for the
dependent coordinates, and let $\mu$ and $\nu$ denote multi-indices for
marking derivatives.  Derivatives of higher order are greater than derivatives
of lower order: if $|\mu|<|\nu|$, then $u^\alpha_\mu \prec u^\beta_\nu$.  If
derivatives have the same order $|\mu|=|\nu|$, then we distinguish two cases: if
the leftmost non-vanishing entry in $\mu-\nu$ is positive, then $u^\alpha_\mu
\prec u^\beta_\nu$; and if $\mu=\nu$ and $\alpha < \beta$, then $u^\alpha_\mu
\prec u^\beta_\nu$.  This is a class-respecting order: if $|\mu|=|\nu|$ and
$\cl \mu < \cl \nu$, then $u^\alpha_\mu \prec u^\beta_\nu$.  Any set of
objects indexed with pairs $(\alpha,\mu)$ can be ordered in an analogous way.
This order of the multi-indices $\mu$ and $\nu$ is called the \emph{degree
  reverse lexicographic ranking}, and we generalise it in such a way that it
places more weight on the multi-indices $\mu$ and $\nu$ than on the numbers
$\alpha$ and $\beta$ of the dependent variables.  This is called the
\emph{term-over-position lift of the degree reverse lexicographic ranking}.

Now the columns within the symbol matrix are ordered descendingly according to
the degree reverse lexicographic ranking for the multi-indices $\mu$ of the
variables $v^\alpha_\mu$ in equation (\ref{gl:Symbol}) and labelled by the
pairs $(\alpha,\mu)$.  (It follows that, if $v^\alpha_\mu$ and $v^\beta_\nu$
are such that $\cl\mu > \cl\nu$, then the column corresponding to
$v^\alpha_\mu$ is left of the column corresponding to $v^\beta_\nu$.)  The
rows are ordered in the same way with regard to the pairs $(\alpha,\mu)$ of
the variables $u^\alpha_\mu$ which def\/ine the classes of the equations
$\mathnormal\Phi^\tau(\mathbf{x},\mathbf{u}^{(q)})=0$.  If two rows are
labelled by the same pair $(\alpha,\mu)$, it does not matter which one comes
f\/irst.

We compute now a row echelon form of the symbol matrix.  We denote the number
of rows where the pivot is of class $k$ by $\bq{q}{k}$, the \emph{indices} of
the symbol $\Nc{q}$, and associate with each such row its \emph{multiplicative
variables} $x^1,\dots,x^k$.  Prolonging each equation only with respect to
its multiplicative variables yields independent equations of order $q+1$, as
each has a dif\/ferent leading term.

\begin{definition}
  If prolongation with respect to the non-multiplicative variables does not
  lead to additional independent equations of order $q+1$, in other words if
  \begin{equation}\label{gl:invsymb}
        \rank{M_{q+1}}
        =
        \sum_{k=1}^n\,k\bq{q}{k} ,
  \end{equation}
  then the symbol $\Nc{q}$ is \emph{involutive}.  The dif\/ferential equation
  $\Rc{q}$ is called \emph{involutive}, if it is formally integrable and its
  symbol is involutive.
\end{definition}

The criterion (\ref{gl:invsymb}) is also known as \emph{Cartan's test}, as it
is analogous to a similar test in the Cartan--K\"ahler theory of exterior
dif\/ferential systems.  We stress again that it is valid only in
$\delta$-regular coordinates (in fact, in other coordinate systems it will
always fail).

\section{The Cartan normal form}

For notational simplicity, we will consider in our subsequent analysis almost
exclusively f\/irst-order equations $\Rc{1}\subseteq\je{1}$.  At least from a
theoretical standpoint, this is not a restriction, as any higher-order
dif\/ferential equation $\Rc{q}$ can be transformed into an equivalent
f\/irst-order one (see for example \cite[Appendix A.3]{Seiler:Buch}).  For these
we now introduce a convenient local representation.

\begin{definition}\label{def:CNF}
  For a f\/irst-order dif\/ferential equation $\Rc{1}$ the following local
  representation, a~special kind of solved form,
  \begin{subequations}\label{gl:CNF}
  \begin{alignat}{3}
 &   u^\alpha_n
        = \phi^\alpha_n\big(\mathbf{x},u^\beta,u^\gamma_j,u^\delta_n\big)
      \qquad & &
      \left\{
            \begin{array}{l@{}}
            1 \le \alpha \le \beta_1^{(n)}, \\
            1 \le j < n, \\
            \bq{1}{n} < \delta \le m,
            \end{array}
        \right.&
    \label{gl:CNF-a}
    \\
  &  u^\alpha_{n-1}
         = \phi^\alpha_{n-1}\big(\xv,u^\beta,u^\gamma_j,u^\delta_{n-1}\big)
     \qquad   && \left\{
            \begin{array}{l@{}}
            1 \le \alpha \le \beta_1^{(n-1)}, \\
            1 \le j < n-1, \\
            \bq{1}{n-1} < \delta \le m,
            \end{array}
        \right.&
      \label{gl:CNF-b}
    \\
&    \cdots \cdots \cdots\cdots \cdots \cdots\cdots \cdots \cdots &&&
    \nonumber
    \\
 &   u^\alpha_1
        = \phi^\alpha_{1}\big(\xv,u^\beta,u^\delta_1\big)
        &        & \left\{
            \begin{array}{l@{}}
            1 \le \alpha \le \beta_1^{(1)}, \\
            \bq{1}{1} < \delta \le m,
            \end{array}
        \right.&
    \label{gl:CNF-c}
    \\
 &   u^\alpha
        = \phi^\alpha\big(\xv,u^\beta\big)
        &
        & \left\{
            \begin{array}{l@{}}
            1 \le \alpha \le \beta_0, \\
            \beta_0 < \beta \le m,
            \end{array}
        \right.&
    \label{gl:CNF-d}
  \end{alignat}
  \end{subequations}
  is called its \emph{Cartan normal form}.  The equations of zeroth order,
  $u^\alpha = \phi^\alpha(\mathbf{x},u^\beta)$, are called \emph{algebraic}.
  The functions $\phi^\alpha_k$ are called the \emph{right sides} of $\Rc{1}$.
  (If, for some $1 \le k \le n$, the number of equations is $\bq{1}{k}=m$,
  then the condition $\bq{1}{k} < \delta \le m$ is empty and no terms
  $u^\delta_k$ appear on the right sides of those equations.)
\end{definition}
Here, each equation is solved for a \emph{principal} derivative of maximal
class $k$ in such a way that the corresponding right side of the equation may
depend on an arbitrary subset of the independent variables, an arbitrary
subset of the dependent variables $u^\beta$ with $1 \le \beta \le \beta_0$,
those derivatives $u^\gamma_j$ for all $1 \le \gamma \le m$ which are of a
class $j<k$ and those derivatives which are of the same class $k$ but are not
principal derivatives.  Note that a principle derivative $u^\alpha_k$ may
depend on another principle derivative $u^\gamma_l$ as long as $l<k$.  The
equations are grouped according to their class in descending order.

\begin{theorem}[Cartan--K{\"a}hler]\label{thm:Cartan-Kaehler}
  Let the involutive differential equation $\Rc{1}$ be locally represented in
  $\delta$-regular coordinates by the system \eqref{gl:CNF-a}, \eqref{gl:CNF-b},
  \eqref{gl:CNF-c}.  Assume that the following initial conditions are given:
  \begin{subequations}\label{gl:Cartan-Kaehler-AWP}
  \begin{alignat}{3}
     & u^\alpha(x^1, \ldots,x^n) = f^\alpha(x^1,\ldots,x^n)  ,
         \qquad && \bq{1}{n} < \alpha \le m ; &
    \\
&    u^\alpha(x^1,\ldots,x^{n-1},0)  = f^\alpha(x^1,\ldots,x^{n-1})  ,
         \qquad && \bq{1}{n-1} < \alpha \le \bq{1}{n}  ; &
    \\
  & \cdots\cdots\cdots \cdots \cdots \cdots & &  \cdots\cdots\cdots\cdots \cdots \cdots &
    \notag
    \\
 &   u^\alpha(x^1,0,\ldots,0) = f^\alpha(x^1)  ,
         \qquad && \bq{1}{1} < \alpha \le \bq{1}{2} ;&
    \\
&    u^\alpha(0,\ldots,0) = f^\alpha  ,
        \qquad  && 1 \le \alpha \le \bq{1}{1} . &
  \end{alignat}
  \end{subequations}
  If the functions $\phi^\alpha_k$ and $f^\alpha$ are $($real-$)$analytic at the
  origin, then this system has one and only one solution that is analytic at
  the origin and satisfies the initial conditions
\eqref{gl:Cartan-Kaehler-AWP}.
\end{theorem}

\begin{proof}
  For the proof, see \cite{Pommaret:SystemsOfPDEs} or \cite[Section 9.4]{Seiler:Buch} and references therein.  The strategy is to split the
  system into subsystems according to the classes of the equations in it (see
  below).  The solution is constructed step by step; each step renders a
  normal system in fewer independent variables to which the
  Cauchy--Kovalevskaya theorem is applied.  Finally, the condition that
  $\Rc{1}$ is involutive leads to further normal systems ensuring that the
  constructed functions are indeed solutions of the full system with respect
  to all independent variables.
\end{proof}

Under some mild regularity assumptions the algebraic equations can always be
solved locally.  From now on, we will assume that any present algebraic
equation has been explicitly solved, reducing thus the number of dependent
variables.  We simplify the Cartan normal form of a~dif\/ferential equation as
given in Def\/inition \ref{def:CNF} into the \emph{reduced Cartan normal form}.
It arises by solving each equation for a derivative $u^\alpha_j$, the
principal derivative, and eliminating this derivative from all other
equations.  Again, the principal derivatives are chosen in such a manner that
their classes are as great as possible.  Now no principal derivative appears
on a right side of an equation (whereas this was possible with the non-reduced
Cartan normal form of Def\/inition \ref{def:CNF}).  All the remaining,
non-principal, derivatives are called \emph{parametric}.  Ordering the
obtained equations by their class, we again can decompose them into
subsystems:
\begin{gather}\label{gl:reducedCNF}
    u^\alpha_k
        = \phi^\alpha_k\big(\xv,\uv,u^\gamma_j\big)
        \qquad
               \left\{
            \begin{array}{l@{}}
            1 \le j \le k \le n, \\
            1 \le \alpha \le \beta_1^{(k)}, \\
            \beta_1^{(j)} < \gamma \le m.
            \end{array}
        \right.
\end{gather}
Note that the values $\beta_1^{(k)}$ are exactly those appearing in the Cartan
test (\ref{gl:invsymb}), as the symbol matrix of a dif\/ferential equation in
Cartan normal form is automatically triangular with the principal derivatives
as pivots.

\begin{definition}
  The \emph{Cartan characters} of $\Rc{1}$ are def\/ined as
  $\alpha_1^{(k)}=m-\beta_1^{(k)}$ and thus equal the number of parametric
  derivatives of class $k$ and order $1$.
\end{definition}

Provided that $\delta$-regular coordinates are chosen, it is possible to
perform a closed form involution analysis for a dif\/ferential equation $\Rc{1}$
in reduced Cartan normal form.  We remark that an ef\/fective test of involution
proceeds as follows (see for example \cite[Remark 7.2.10]{Seiler:Buch}).  Each
equation in (\ref{gl:reducedCNF}) is prolonged with respect to each of its
non-multiplicative variables.  The arising second-order equations are
simplif\/ied modulo the original system and the prolongations with respect to
the multiplicative variables.  The symbol $\Nc{1}$ is involutive if, and only
if, after the simplif\/ication none of the prolonged equations is of
second-order any more.  The equation $\Rc{1}$ is involutive if, and only if,
all new equations even simplify to zero, as any remaining f\/irst-order equation
would be an integrability condition.

In order to apply this test, we now prove two helpful lemmata.  We introduce
the set $\mathcal{B} := \bigl\{(\alpha,i) \in \mathbbm{N} \times \mathbbm{N}
\colon\bigr. u^\alpha_i \ \text{is a}$ $\bigl.\text{principal} \
\text{derivative}\bigr\}$, and for each $(\alpha,i) \in \mathcal{B}$ we def\/ine
$\mathnormal\Phi^\alpha_i \, := \, u^\alpha_i - \phi^\alpha_i$.  Using the
contact f\/ields (\ref{gl:ContactFields}), any prolongation of some
$\mathnormal\Phi^\alpha_i$ can be expressed in the following form.

\begin{lemma}
  Let the differential equation $\Rc{1}$ be represented in the reduced Cartan
  normal form given by equation \eqref{gl:reducedCNF}.  Then for any
  $(\alpha,i) \in \mathcal{B}$ and $1 \le j \le n$, we have
  \begin{gather}\label{gl:PhiFortsetzung}
    D_j \mathnormal\Phi^\alpha_i
    =
    u_{ij}^\alpha
    -
    C^{(1)}_j(\phi_i^\alpha)
    -
    \sum_{h=1}^i
        \sum_{\gamma=\bq{1}{h}+1}^m
            u_{hj}^\gamma
                C_\gamma^h(\phi_i^\alpha).
  \end{gather}
\end{lemma}

\begin{proof}
  By straightforward calculation; see \cite[Lemma 2.5.5]{Fesser:2008}.
\end{proof}

For $j>i$, the prolongation $D_j \mathnormal\Phi_i^\alpha$ is
non-multiplicative, otherwise it is multiplicative.  Now let $j>i$, so that
equation (\ref{gl:PhiFortsetzung}) shows a non-multiplicative prolongation,
and assume that we are using $\delta$-regular coordinates.  According to our
test, the symbol $\Nc{1}$ is involutive if, and only if, it is possible to
eliminate on the right hand side of (\ref{gl:PhiFortsetzung}) all second-order
derivatives by adding multiplicative prolongations.

If the dif\/ferential equation is not involutive, then the dif\/ference
\[
        D_j\mathnormal\Phi_i^\alpha
        -
        D_i\mathnormal\Phi_j^\alpha
        +
        \sum_{h=1}^i
                \sum_{\gamma=\bq{1}{h}+1}^{\bq{1}{j}}
                        C_\gamma^h(\phi^\alpha_h)
                        D_h\mathnormal\Phi^\gamma_j
\]
does not necessarily vanish but may yield an obstruction to involution for any
$(\alpha,i) \in \mathcal{B}$ and any $i < j \le n$.  The next lemma gives all
these obstructions to involution for a f\/irst-order system in reduced Cartan
normal form.
\begin{lemma}\label{lem:Monster}
  Assume that $\delta$-regular coordinates are used.  For an equation in
  Cartan normal form and indices $i<j$ and $\alpha$ such that $(\alpha,i) \in
  \mathcal{B}$, we have the equality
\begin{gather}
     D_j\mathnormal\Phi_i^\alpha
    -
    D_i\mathnormal\Phi_j^\alpha
    +
    \sum_{h=1}^i
        \sum_{\gamma = \beta_1^{(h)}+1}^{\bq{1}{j}}
            C_\gamma^h(\phi^\alpha_i)
                D_h\mathnormal\Phi^\gamma_j
    %
    \nonumber
    \\
    %
\qquad{}    =
    C^{(1)}_i(\phi_j^\alpha) -
        C^{(1)}_j(\phi_i^\alpha)
        - \sum_{h=1}^i
            \sum_{\gamma = \beta_1^{(h)}+1}^{\bq{1}{j}}
                C_\gamma^h(\phi_i^\alpha)
                C^{(1)}_h(\phi_j^\gamma)
    \label{line:Monster-IntegBed}
\\
\qquad\quad{}
    -
    \sum_{h=1}^{i-1}
        \sum_{\delta = \beta_1^{(h)}+1}^m
            u^\delta_{hh}
                \left[
                \sum_{\gamma = \bq{1}{h}+1}^{\bq{1}{j}}
                    C^h_\gamma(\phi^\alpha_i)
                    C^h_\delta(\phi^\gamma_j)
                \right]
    \label{line:Monster-hh}
\end{gather}
  \begin{subequations}\label{line:Monster-hk1}
\begin{gather}
\qquad\quad{}
    -
    \sum_{1 \le h < k < i}
        \left\{
            \sum_{\delta = \bq{1}{h}+1}^{\bq{1}{k}}
        \right.
                u^\delta_{hk}
                                \left.
                    \left[
                        \sum_{\gamma = \bq{1}{k}+1}^{\bq{1}{j}}
                            C^k_\gamma(\phi^\alpha_i)
                            C^h_\delta(\phi^\gamma_j)
                    \right]
                \right.
    \label{line:Monster-hk-a}
\\
 \qquad\quad{}   +
    \sum_{\delta = \bq{1}{k}+1}^m
        u^\delta_{hk}
                \left.
            \left[
                \sum_{\gamma = \bq{1}{h}+1}^{\bq{1}{j}}
                    C^h_\gamma(\phi^\alpha_i)
                    C^k_\delta(\phi^\gamma_j) +
                \sum_{\gamma = \bq{1}{k}+1}^{\bq{1}{j}}
                    C^k_\gamma(\phi^\alpha_i)
                    C^h_\delta(\phi^\gamma_j)
            \right]
        \right\}
    \label{line:Monster-hk-b}
\end{gather}
 \end{subequations}
  \begin{subequations}\label{line:Monster-hi}
  \begin{gather}
\qquad\quad{}
    -
    \sum_{h=1}^{i-1}
        \left\{
            \sum_{\delta = \bq{1}{h}+1}^{\bq{1}{i}}\!\!\!\!
                u^\delta_{hi}
        \right.
                            \left.
                    \left[
                        -C^h_\delta(\phi^\alpha_j) +\!\!\!\!\!
                            \sum_{\gamma = \bq{1}{i}+1}^{\bq{1}{j}}
                                C^i_\gamma(\phi^\alpha_i)
                                C^h_\delta(\phi^\gamma_j)
                    \right]
                \right.
    \label{line:Monster-hi-a}
\\
\qquad\quad{}    +
    \sum_{\delta = \bq{1}{i}+1}^m
        u^\delta_{hi}
        \left.
            \left[
                -C^h_\delta(\phi^\alpha_j) +\!\!\!\!\!
                    \sum_{\gamma = \bq{1}{i}+1}^{\bq{1}{j}}
                        C^i_\gamma(\phi^\alpha_i)
                        C^h_\delta(\phi^\gamma_j) +\!\!\!\!\!
                    \sum_{\gamma = \bq{1}{h}+1}^{\bq{1}{j}}
                        C^h_\gamma(\phi^\alpha_i)
                        C^i_\delta(\phi^\gamma_j)
            \right]
        \right\}
    \label{line:Monster-hi-b}
  \end{gather}
  \end{subequations}
  \begin{gather}
\qquad\quad{}
    -
    \sum_{\substack{h=1 \\ i+1 \le k < j}}^{i-1}
        \sum_{\delta = \bq{1}{k}+1}^m
            u^\delta_{hk}
                \left[
                    \sum_{\gamma = \bq{1}{h}+1}^{\bq{1}{j}}
                        C^h_\gamma(\phi^\alpha_i)
                        C^k_\delta(\phi^\gamma_j)
                \right]
    \label{line:Monster-hk2}
\\
\qquad\quad{}-
    \sum_{k=i}^{j-1}
        \sum_{\delta = \bq{1}{k}+1}^m
            u^\delta_{ik}
                \left[
                    -C_\delta^k(\phi^\alpha_j)
                        +\sum_{\gamma = \bq{1}{i}+1}^{\bq{1}{j}}
                            C_\gamma^i(\phi^\alpha_i)
                            C_\delta^k(\phi^\gamma_j)
                \right]
    \label{line:Monster-ik}
\\
\qquad\quad{} -
    \sum_{h = 1}^{i-1}
        \sum_{\delta = \bq{1}{j}+1}^m
            u^\delta_{hj}
                \left[ C^h_\delta(\phi^\alpha_i) +
                        \sum_{\gamma = \bq{1}{h}+1}^{\bq{1}{j}}
                            C^h_\gamma(\phi^\alpha_i)
                            C^j_\delta(\phi^\gamma_j)
                \right]
    \label{line:Monster-hj}
\\
\qquad\quad{}-
    \sum_{\delta = \bq{1}{j}+1}^m
        u^\delta_{ij}
            \left[C^i_\delta(\phi^\alpha_i)
                -C^j_\delta(\phi^\alpha_j)
                +\sum_{\gamma = \bq{1}{i}+1}^{\bq{1}{j}}
                    C^i_\gamma(\phi^\alpha_i)
                    C^j_\delta(\phi^\gamma_j)
            \right].
    \label{line:Monster-ij}
  \end{gather}
\end{lemma}
\begin{proof}
  By a tedious but fairly straightforward computation, see \cite[pages 48--55]{Fesser:2008}.
\end{proof}

In line (\ref{line:Monster-IntegBed}) we have collected all terms which are of
lower than second order.  Furthermore, none of the appearing second-order
derivatives is of a form that it could be eliminated by adding some
multiplicative prolongation.  Hence, under our assumption of $\delta$-regular
coordinates, the symbol $\Nc{1}$ is involutive if, and only if, all the
expressions in square brackets vanish.  The dif\/ferential equation $\Rc{1}$ is
involutive if, and only if, in addition line (\ref{line:Monster-IntegBed})
vanishes, as it represents an integrability condition.  Thus Lemma
\ref{lem:Monster} gives us all obstructions to involution for $\Rc{1}$ in
explicit form.  They will reappear in the proof of the existence theorem for
integral distributions in Section \ref{section:ExThmIntegralDist}.

\section{The Vessiot distribution}

By Proposition \ref{prop:contactmapq}, the tangent spaces
$T_\rho(\im{j_q\sigma})$ of prolonged sections at points $\rho\in\je{q}$ are~always subspaces of the contact distribution $\C_q|_{\rho}$.  If the section
$\sigma$ is a solution of the dif\/ferential equation $\Rc{q}$, then by
def\/inition it furthermore satisf\/ies $\im{j_q\sigma}\subseteq\Rc{q}$, and
therefore $T(\im{j_q\sigma})\subseteq T\Rc{q}$.  Hence, the following
construction suggests itself.
\begin{definition}\label{def:vessiot}
  The\/ \emph{Vessiot distribution} of a dif\/ferential equation
  $\Rc{q}\subseteq\je{q}$ is the distribution\/ $\V[\Rc{q}]\subseteq T\Rc{q}$
  def\/ined by
\begin{equation*}\label{gl:vessiot}
	T\iota\bigl(\V[\Rc{q}]\bigr)
	=
	T\iota\bigl(T\Rc{q}\bigr)\cap\C_q|_{\Rc{q}} .
\end{equation*}
\end{definition}

Note that the Vessiot distribution is not necessarily of constant rank along
$\Rc{q}$ (just like the symbol $\Nc{q}$); for simplicity, we will almost
always assume here that this is the case.  This def\/inition of the Vessiot
distribution is not the one usually found in the literature.  But the
equivalence to the standard approach is an elementary exercise in computing
with pull-backs.
\begin{proposition}\label{prop:vessclass}
The Vessiot distribution satisfies
	$\V[\Rc{q}]
	=
	(\iota^*\mathcal{C}_q^0)^0$.
\end{proposition}

For a dif\/ferential equation given in explicitly solved form, the inclusion map
$\iota \colon \Rc{q} \to J_q\pi$ is available in closed form and can be used
to calculate the pull-back of the contact forms.  This has the advantage of
keeping the calculations within a space of smaller dimension, namely the
submanifold $\Rc{q}$.  Thereby regarding the dif\/ferential equation as a
manifold in its own right, we bar its coordinates to distinguish them from
those of the jet bundle.

\begin{example}\label{bsp:Welle1}
  Consider the f\/irst-order system given by the representation $\Rc{1} \colon
  u_t-v = v_t-w_x = u_x-w = 0.$ Then from the prolongations $u_{xt}-v_x=0$ and
  $u_{xt}-w_t=0$ follows the integrability condition $w_t=v_x$ by elimination
  of the second-order derivative $u_{xt}$.  Thus, the f\/irst projection of
  $\Rc{1}$ is (in Cartan normal form) represented by
  \[
  \Rcp{1}{1} \colon \ \left\{
        \begin{array}{l}
        u_t  = v, \qquad  v_t  = w_x  , \qquad  w_t  = v_x , \\
        u_x  = w .
        \end{array}
    \right.
  \]
  It is not dif\/f\/icult to verify that the projected equation $\Rcp{1}{1}$ is
  involutive.  For coordinates on~$J_1\pi$ choose
  $x$, $t$; $u$, $v$, $w$; $u_x$, $v_x$, $w_x$, $u_t$, $v_t$, $w_t$.  Since $\Rcp{1}{1}$ is represented by a
  system in solved form, it is natural to choose appropriate local coordinates
  for $\Rcp{1}{1}$, which we bar to distinguish them:
  $\overline{x}$, $\overline{t}$; $\overline{u}$, $\overline{v}$, $\overline{w}$;
  $\overline{v_x}$, $\overline{w_x}$.  The contact codistribution for $J_1\pi$ is
  generated by
  \[
    \omega^1 = du - u_xdx - u_tdt , \qquad
    \omega^2 = dv - v_xdx - v_tdt , \qquad
    \omega^3 = dw - w_xdx  - w_tdt .
  \]
  The tangent space $T\Rcp{1}{1}$ is spanned by
  $\partial_{\overline{x}}$, $\partial_{\overline{t}}$, $\partial_{\overline{u}}$,
  $\partial_{\overline{v}}$, $\partial_{\overline{w}}$, $\partial_{\overline{v_x}}$,
  $\partial_{\overline{w_x}}$ and $T\iota(T\Rcp{1}{1})$ therefore by the f\/ields
  $\partial_x$, $\partial_t$, $\partial_u$, $\partial_v+\partial_{u_t}$,
  $\partial_w+\partial_{u_x}$, $\partial_{v_x}+\partial_{w_t}$ and
  $\partial_{w_x}+\partial_{v_t}$.  This space is annihilated by
  \begin{displaymath}
    \omega^4 = du_x - dw,\qquad
    \omega^5 = dv_x - dw_t , \qquad
    \omega^6 = dw_x - dv_t  ,\qquad
    \omega^7 = du_t - dv  .
  \end{displaymath}
  These seven one-forms $\omega^1,\dots,\omega^7$ annihilate the Vessiot
  distribution $\V[\Rcp{1}{1}]$, which is spanned by the four vector f\/ields
    \begin{alignat*}{3}
    & X_1  = \partial_x + u_x \partial_u + v_x \partial_v + w_x \partial_w +
    v_x\partial_{u_t} + w_x\partial_{u_x} , \qquad &&
    X_3  = \partial_{v_x} + \partial_{w_t} , & \\
 &   X_2  = \partial_t + u_t \partial_u + v_t \partial_v + w_t \partial_w +
    v_t\partial_{u_t} + w_t\partial_{u_x} , \qquad & &
    X_4  = \partial_{v_t} + \partial_{w_x} .
  \end{alignat*}
  In local coordinates on $\Rcp{1}{1}$, these four vector f\/ields become
  \begin{alignat*}{3}
 &   \bar{X}_1 =
    \partial_{\overline{x}}
    + \overline{w} \partial_{\overline{u}}
    + \overline{v_x} \partial_{\overline{v}}
    + \overline{w_x} \partial_{\overline{w}} ,\qquad & &
    \bar{X}_3  =
    \partial_{\overline{v_x}}  ,&  \\
  &   \bar{X}_2
   =
    \partial_{\overline{t}}
    + \overline{v} \partial_{\overline{u}}
    + \overline{w_x} \partial_{\overline{v}}
    + \overline{v_x} \partial_{\overline{w}} , \qquad & &
    \bar{X}_4  =
    \partial_{\overline{w_x}} . &
  \end{alignat*}
  They satisfy $\iota_*\bar{X}_i = X_i$, as a simple calculation using the
  Jacobian matrix for $T\iota$ shows.  The vector f\/ields $\bar{X}_i$ are
  annihilated by the pull-backs of the contact forms, $\iota^*\omega^1 =
  d\overline{u} - \overline{u_x}d\overline{x} - \overline{u_t}d\overline{t}$,
  $\iota^*\omega^2 = d\overline{v} - \overline{x_x}d\overline{x} -
  \overline{x_t}d\overline{t}$ and $\iota^*\omega^3 = d\overline{w} -
  \overline{w_x}d\overline{x} - \overline{w_t}d\overline{t}$ (the pullbacks of
  the remaining four one-forms $\omega^4,\dots,\omega^7$ trivially vanish on
  $\Rcp{1}{1}$).
\end{example}

For a fully nonlinear dif\/ferential equation $\Rc{q}$, in particular an
implicit one, this approach to compute its Vessiot distribution $\V[\Rc{q}]$
via a pull-back is in general not ef\/fectively feasible.  However, applying
directly our def\/inition of $\V[\Rc{q}]$, it is easily possible even for such
equations to determine ef\/fectively $T\iota(\V[\Rc{q}])$, in other words: to
realize it as a subbundle of $T(\je{q})|_{\Rc{q}}$.  The contact f\/ields
(\ref{gl:ContactFields}) form a basis for $\C_q$.  It follows that for any
vector f\/ield $\bar{X} \in \V[\Rc{q}]$, coef\/f\/icients $a^i, b^\alpha_\mu \in
\mathcal{F}(\Rc{q})$, where $1 \le i \le n$, $1 \le \alpha \le m$ and
$|\mu|=q$, exist such that
\begin{equation}\label{gl:VessiotAnsatz}
  \iota_*\bar{X} =
      a^i C_i^{(q)} + b^\alpha_\mu C_\alpha^\mu .
\end{equation}
If the dif\/ferential equation $\Rc{q}$ is locally represented by
$\mathnormal\Phi^\tau = 0$, where $1 \le \tau \le t$, it follows from the
tangency of the vector f\/ields in $\V[\Rc{q}]$ that
$d\mathnormal\Phi^\tau(\iota_*\bar{X}) = \iota_*\bar{X}(\mathnormal\Phi^\tau)
= 0$ and thus the coef\/f\/icient functions must satisfy the following system
of linear equations:
\begin{equation}\label{gl:VessiotfeldAllgemein}
        C_i^{(q)}(\mathnormal\Phi^\tau) a^i
        +
        C_\alpha^\mu (\mathnormal\Phi^\tau) b^\alpha_\mu
        = 0,
\end{equation}
where $1 \le \tau \le t$.  Note that this approach to determine the Vessiot
distribution is closely related to prolonging the dif\/ferential equation
$\Rc{q}$ and requires essentially the same computations.  Indeed, the formal
derivative (\ref{gl:formderiv}) can be written in the form
\begin{equation}\label{gl:prolcont}
 D_i\Phi^\tau=C_i^{(q)}(\Phi^\tau)+
              C^\mu_\alpha(\Phi^\tau)u^\alpha_{\mu+1_i}=0
\end{equation}
and in the context of the order-by-order determination of formal power series
solutions (see for example \cite[Section 2.3]{Seiler:Buch}) these equations
are considered an inhomogeneous system for the Taylor coef\/f\/icients of order
$q+1$ depending on the lower order coef\/f\/icients.  Taking this point of view,
we may call (\ref{gl:VessiotfeldAllgemein}) the ``projective'' version of
(\ref{gl:prolcont}).  In fact for $n=1$, that is, for ordinary dif\/ferential
equations, this is even true in a rigorous sense.

\begin{example}\label{bsp:kugel}
  We consider the fully nonlinear f\/irst-order ordinary dif\/ferential equation
  $\Rc{1}$ locally def\/ined by $(u')^2+u^2+x^2=1$.  The contact distribution
  $\C_1$ is spanned by the two vector f\/ields $X_1=\partial_x+u'\partial_u$ and
  $X_2=\partial_{u'}$ and the Vessiot distribution $T\iota(\V[\Rc{q}])$
  consists of all linear combinations of these two f\/ields which are tangent to
  $\Rc{1}$.  Setting $\omega=xdx+udu+u'du'$, we thus have to solve the linear
  equation $\omega(aX_1+bX_2)=0$ in order to determine $T\iota(\V[\Rc{q}])$.
  Its solution requires a case distinction (which is typical for fully
  nonlinear dif\/ferential equations). If $u'\neq0$, then we f\/ind
  \begin{displaymath}
    T\iota(\V[\Rc{q}])=
    \langle u'\partial_x+(u')^2\partial_u-(x+u'u)\partial_{u'}\rangle .
  \end{displaymath}
  For $u'=0$ and $x\neq0$, the Vessiot distribution is spanned by the vertical
  contact f\/ield $X_2$.  Finally, for $x=u'=0$ the rank of the Vessiot
  distribution jumps, as at these points the whole contact plane is contained
  in it.
\end{example}

The def\/inition of the symbol $\Nc{q}$ and of the Vessiot distribution
$\V[\Rc{q}]$, respectively, of a~dif\/fe\-ren\-tial equation $\Rc{q}\subseteq\je{q}$
immediately imply the following generalisation of the above discussed
splitting of the contact distribution $\C_q=V\pi^q_{q-1}\oplus\H$ (such a
splitting of the Vessiot distribution is also discussed by Lychagin and
Kruglikov \cite{KruglikovLychagin:2008, Lychagin:1995} where the Vessiot
distribution is called ``Cartan distribution'').

\begin{proposition}\label{prop:SymbolVessiotDistribution}
  For any differential equation\/ $\Rc{q}$, its symbol is contained in the
  Vessiot distribution:\/ $\Nc{q} \subseteq \V[\Rc{q}]$.  The Vessiot
  distribution can therefore be decomposed into a direct sum
  \begin{gather}\label{gl:VessiotDistributionSymbolH}
        \V[\Rc{q}] = \Nc{q} \oplus \H
  \end{gather}
  for some complement $\H$ $($which is not unique$)$.
\end{proposition}

Such a complement $\H$ is necessarily transversal to the f\/ibration
$\Rc{q}\rightarrow\X$ and thus leads naturally to connections:  provided
$\dim{\H}=n$, it may be considered the horizontal bundle of a~connection of
the f\/ibred manifold $\Rc{q}\rightarrow\X$.
\begin{definition}
  Any such connection is called a \emph{Vessiot connection} for $\Rc{q}$.
\end{definition}
In general, the Vessiot distribution $\V[\Rc{q}]$ is not involutive (that is,
closed under the Lie bracket; an exception are formally integrable equations
of f\/inite type \cite[Remark 9.5.8]{Seiler:Buch}), but it may contain
involutive subdistributions.  If these are furthermore transversal (to the
f\/ibration $\mathcal{R}_q\rightarrow\X$) and of dimension $n$, then they def\/ine
a \emph{flat} Vessiot connection.
\begin{lemma}
  If the section $\sigma:\X\rightarrow\E$ is a solution of the equation
  $\Rc{q}$, then the tangent bundle $T(\im{j_q\sigma})$ is an $n$-dimensional
  transversal involutive subdistribution of $\V[\Rc{q}]|_{\im{j_q\sigma}}$.
  Conversely, if $\U\subseteq\V[\Rc{q}]$ is an $n$-dimensional transversal
  involutive subdistribution, then any integral manifold of $\U$ has locally
  the form $\im{j_q\sigma}$ for a solution $\sigma$ of $\Rc{q}$.
\end{lemma}
\begin{proof}
  Let $\sigma$ be a local solution of the equation $\Rc{q}$.  Then it
  satisf\/ies by Def\/inition \ref{def:Solution} $\im j_q\sigma \subseteq \Rc{q}$
  and thus $T(\im j_q\sigma) \subseteq T\Rc{q}$.  Besides, by the def\/inition
  of the contact distribution, for $\xv \in \X$ with $j_q\sigma(\xv)=\rho \in
  J_q\pi$, the tangent space $T_\rho(\im j_q\sigma)$ is a subspace of
  $\C_q|_\rho$.  By def\/inition of the Vessiot distribution, it follows
  $T_\rho(\im j_q\sigma) \subseteq T\iota(T_\rho\Rc{q}) \cap \C_q|_\rho$,
  which proves the f\/irst claim.

  Now let $\U\subseteq\V[\Rc{q}]$ be an $n$-dimensional transversal involutive
  subdistribution.  Then according to the Frobenius theorem, $\U$ has
  $n$-dimensional integral manifolds.  By def\/inition, $T\iota(\V[\Rc{q}])
  \subseteq \C_q|_{\Rc{q}}$; this characterises prolonged sections.  Hence,
  for any integral manifold of $\U$ there is a local section $\sigma$ such
  that the integral manifold is of the form $\im j_q\sigma$.  Furthermore, the
  integral manifold is a subset of $\Rc{q}$.  Thus it corresponds to a local
  solution of $\Rc{q}$.
\end{proof}

This simple observation forms the basis of Vessiot's approach to the analysis
of $\Rc{q}$: he proposed to construct all f\/lat Vessiot connections.  Before we
do this, we f\/irst show how integral elements appear in this program.
\begin{proposition}\label{prop:intdist}
  Let $\U\subseteq\V[\Rc{q}]$ be a transversal subdistribution of the Vessiot
  distribution of constant rank $k$.  Then the spaces $\U_\rho$ are
  $k$-dimensional integral elements for all points $\rho\in\Rc{q}$ if, and
  only if, $[\U,\U] \subseteq \V[\Rc{q}]$.
\end{proposition}
\begin{proof}
  Let $\{\omega_1,\dots,\omega_r\}$ be a basis of the codistribution
  $\iota^*\C_q^0$.  Then an algebraic basis of the ideal $\I[\Rc{q}]$ is
  $\{\omega_1,\dots,\omega_r, d\omega_1,\dots,d\omega_r\}$.  Any vector f\/ield
  $X\in\U$ trivially satisf\/ies $\omega_i(X)=0$ by Proposition
  \ref{prop:vessclass}.  For arbitrary f\/ields $X_1,X_2\in\U$, we have
  \[
        d\omega_i(X_1,X_2)
        =
        X_1\bigl(\omega_i(X_2)\bigr)
        -
        X_2\bigl(\omega_i(X_1)\bigr)
        +
        \omega_i\bigl([X_1,X_2]\bigr)  .
  \]
  The f\/irst two summands on the right hand side vanish trivially and the
  remaining equation implies our claim.
\end{proof}

We call a subdistribution $\U\subseteq\V[\Rc{q}]$ satisfying the conditions of
Proposition \ref{prop:intdist} an \emph{integral distribution} for the
dif\/ferential equation $\Rc{q}$.  In the literature
\cite{Stormark:LiesStructuralApproach}, the terminology ``involution'' is
common for such distributions which, however, is confusing.  Note that
generally an integral distribution is \emph{not} integrable; the name only
ref\/lects the fact that it consists of integral elements.

A general dif\/ferential equations $\Rc{q}$ does not necessarily possess any
Vessiot connection (not even a non-f\/lat one).  Their existence is linked to
the absence of integrability conditions.  More precisely, we obtain the
following characterisation.

\begin{proposition}\label{prop:Existenz-H-notwendighinreichend}
  Let $\Rc{q}$ be a differential equation.  Then its Vessiot distribution
  possesses locally a direct decomposition with an $n$-dimensional complement
  $\mathcal{H}$ such that $\V[\Rc{q}] = \Nc{q} \oplus \mathcal{H}$ if, and
  only if, there are no integrability conditions which arise as the
  prolongation of equations of lower order in the system.
\end{proposition}

\begin{proof}
  Let $\Rc{q}$ be locally represented by
\[
    \Rc{q} \colon \
        \left\{
                \begin{array}{l}
                        \mathnormal\Phi^\tau\big(\xv,\uv^{(q)}\big),  \\
                        \mathnormal\Psi^\sigma\big(\xv,\uv^{(q-1)}\big),
                \end{array}
        \right.
\]
  such that the equations $\mathnormal\Phi^\tau(\xv,\uv^{(q)}) = 0$ do not
  imply lower-order equations which are independent of the equations
  $\mathnormal\Psi^\sigma(\xv,\uv^{(q-1)}) = 0$.  Let $\uv_{(q)}$ denote the
  subset of all derivatives of order~$q$ only; then the Jacobi matrix
  $\bigl(\partial \mathnormal\Phi^\tau(\xv,\uv^{(q)})/\partial
  \uv_{(q)}\bigr)$ has maximal rank.  If we proceed as in the last proof, then
  the ansatz (\ref{gl:VessiotAnsatz}) for the determination of the Vessiot
  distribution yields for the above representation the linear system
  \begin{gather}
        C_i^{(q)}(\mathnormal\Phi^\tau) a^i
        +
        C_\alpha^\mu (\mathnormal\Phi^\tau) b^\alpha_\mu
         = 0 ,\label{gl:Spezialsystem}
        \\
        C^{(q)}_i(\mathnormal\Psi^\sigma)a^i
                 = 0 .
        \nonumber
  \end{gather}
  Here, the matrix $C_\alpha^\mu (\mathnormal\Phi^\tau)$ has maximal rank,
  too; thus the equations $C_i^{(q)}(\mathnormal\Phi^\tau) a^i + C_\alpha^\mu
  (\mathnormal\Phi^\tau) b^\alpha_\mu = 0$ can be solved for a subset of the
  unknowns $b^\alpha_\mu$.  And since no terms of order $q$ are present in
  $\mathnormal\Psi^\sigma(\xv,\uv^{(q-1)}) = 0$, we have $C^{(q)}_i
  (\mathnormal\Psi^\sigma) = D_i\mathnormal\Psi^\sigma$.  Recall that we
  consider the Vessiot distribution, and thus the linear system
  (\ref{gl:Spezialsystem}), only on $\Rc{q}$.  It follows that the subsystem
  $C^{(q)}_i(\mathnormal\Psi^\sigma)a^i=0$ becomes trivial if, and only if, no
  integrability conditions arise from the prolongation of lower order
  equations.  And if, and only if, this is the case, then
  (\ref{gl:Spezialsystem}) has for each $1 \le j \le n$ a solution where $a^j
  = 1$ while all other $a^i$ are zero.  The existence of such a solution is
  equivalent to the existence of an $n$-dimensional transversal complement
  $\H$.
\end{proof}

From the proof of this proposition follows that for the determination of the
Vessiot distribution $\V[\Rc{q}]$, equations of order less than $q$ in the
local representation of $\Rc{q}$ can be ignored if there are no integrability
conditions which arise from equations of lower order.  It is the integrability
conditions which arise as prolongations of lower order equations that hinder
the construction of $n$-dimensional complements, while those which follow from
the relations between cross derivatives do not inf\/luence this approach.

\begin{remark}
  If one does not care about the distinction between dif\/ferent kinds of
  integrability conditions and simply requires that $\Rc{q}=\Rcp{q}{1}$
  (meaning that no integrability conditions at all appear in the f\/irst
  prolongation of $\Rc{q}$), then one can provide a more geometric proof for
  the existence of an $n$-dimensional complement (of course, in contrast to
  Proposition \ref{prop:Existenz-H-notwendighinreichend}, the converse is not
  true then).

  The assumption $\Rc{q}=\Rcp{q}{1}$ implies that to every point
  $\rho\in\Rc{q}$ at least one point $\hat\rho\in\Rc{q+1}$ with
  $\pi^{q+1}_q(\hat\rho)=\rho$ exists.  We choose such a point and consider
  $\im{\Gamma_{q+1}(\hat\rho)}\subset T_\rho(\je{q})$.  By def\/inition of the
  contact map $\Gamma_{q+1}$, this is an $n$-dimensional transversal subset of
  $\C_q|_{\rho}$. Thus there only remains to show that it is also tangential
  to $\Rc{q}$, as then we can def\/ine a complement by
  $T_\rho\iota(H_\rho)=\im{\Gamma_{q+1}(\hat\rho)}$.  But this tangency is a
  trivial consequence of $\hat\rho\in\Rc{q+1}$; using for example the local
  coordinates expression (\ref{gl:gammaq}) for the contact map and a local
  representation $\Phi^\tau=0$ of $\Rc{q}$, one immediately sees that the
  vector $v_i=\Gamma_{q+1}(\hat\rho,\partial_{x^i})\in T_\rho(\je{q})$
  satisf\/ies $d\Phi^\tau|_{\rho}(v_i)=D_i\Phi^\tau(\hat\rho)=0$ and thus is
  tangential to $\Rc{q}$.

  Hence it is possible to construct for each point $\rho\in\Rc{q}$ a
  transversal complement $\H_\rho$ such that
  $\V_\rho[\Rc{q}]=(\Nc{q})_\rho\oplus\H_\rho$.  There remains to show that
  these complements can be chosen so that they form a smooth distribution.
  Our assumption $\Rc{q}=\Rcp{q}{1}$ implies that the restricted projection
  $\hat\pi^{q+1}_q:\Rc{q+1}\rightarrow\Rc{q}$ is a surjective submersion,
  that is, it def\/ines a f\/ibred manifold.  Thus if we choose a section
  $\gamma:\Rc{q}\rightarrow\Rc{q+1}$ and then always take
  $\hat\rho=\gamma(\rho)$, it follows immediately that the corresponding
  complements $\H_\rho$ def\/ine a smooth distribution as required.
\end{remark}

\begin{example}
  Consider again the dif\/ferential equation $\Rc{1}$ in Example
  \ref{bsp:Welle1}.  It is locally represented by the same equations as
  $\Rcp{1}{1}$, except that the integrability condition $w_t=v_x$ is missing.
  The matrix of $T\iota$ for the system $\Rc{1}$ has eleven rows and eight
  columns---one column more than the symbol matrix for the system
  $\Rcp{1}{1}$.  The symbol $T\iota(\Nc{1})$ of $\Rc{1}$ is $3$-dimensional,
  spanned by $\partial_{v_x}$, $\partial_{w_t}$ and
  $\partial_{w_x}+\partial_{v_t}$, while the symbol $T\iota(\Ncp{1}{1})$ of
  $\Rcp{1}{1}$ has dimension $2$ and is spanned by
  $\partial_{v_x}+\partial_{w_t}$ and $\partial_{w_x}+\partial_{v_t}$.  But
  the one-forms $\omega^1$, $\omega^2$ and $\omega^3$ (and their pull-backs
  $\iota^*\omega^1=
  d\overline{u}-\overline{u_x}d\overline{x}-\overline{u_t}d\overline{t}$,
  $\iota^*\omega^2=
  d\overline{v}-\overline{v_x}d\overline{x}-\overline{v_t}d\overline{t}$ and
  $\iota^*\omega^3=
  d\overline{w}-\overline{w_x}d\overline{x}-\overline{w_t}d\overline{t}$) are
  the same, and therefore the coordinate expressions for the vector f\/ields
  $X_1$ and $X_2$ (and their representations
  $\bar{X}_1=\partial_{\overline{x}}+\overline{u_x}\partial_{\overline{u}}+
  \overline{v_x}\partial_{\overline{v}}+\overline{w_x}\partial_{\overline{w}}$
  and
  $\bar{X}_2=\partial_{\overline{t}}+\overline{u_t}\partial_{\overline{u}}+
  \overline{v_t}\partial_{\overline{v}}+\overline{w_t}\partial_{\overline{w}}$
  in $T\Rc{1}$ and $T\Rcp{1}{1}$) look alike (see Example \ref{bsp:Welle1} for
  their representations).  The integrability condition $w_t=v_x$ does not
  inf\/luence the results as it stems from the equality of the cross
  derivatives, $u_{tx}=v_x$ and $u_{xt}=w_t$, and not from the prolongation of
  a lower order equation.

  Now consider for comparison the dif\/ferential equation which is locally
  represented by
  \[
  \tilde{\mathcal{R}}_{1} \colon \
    \left\{
        \begin{array}{l}
        u_t = v  , \qquad v_t = w_x , \qquad w_t = v_x  , \\
        u_x  = w  , \\
        u  = x  .
        \end{array}
    \right.
  \]
  It arises from the system $\Rcp{1}{1}$ in Example~\ref{bsp:Welle1} by adding
  the algebraic equation $u=x$.  Proceeding as in Example \ref{bsp:Welle1}, we
  f\/ind that the Vessiot distribution $\V[\tilde{\mathcal{R}}_{1}]$ is spanned
  by the three vector f\/ields
  \[
        \partial_{\overline{v_x}},
        \qquad
        \partial_{\overline{w_x}}
        \qquad \text{and} \qquad
        \overline{v}\partial_{\overline{x}}
        +
        (1-\overline{w})\partial_{\overline{t}}
        +
        \bigl(\overline{w_x}(1-\overline{w})+\overline{v}\overline{v_x}\bigr)
            \partial_{\overline{v}}
        +
        \bigl(\overline{v_x}(1-\overline{w})+\overline{v}\overline{w_x}\bigr)
            \partial_{\overline{w}}
   .
  \]
  The vector f\/ields $\partial_{\overline{v_x}}$ and
  $\partial_{\overline{w_x}}$ generate the symbol $\tilde{\mathcal{N}}_{1}$;
  any of its complements in $\V[\tilde{\mathcal{R}}_{1}]$ is one-dimensional
  and, as the dimension of the base space is two, none of them has the right
  dimension to be the horizontal space of a connection.

  The reason for this is that $\tilde{\mathcal{R}}_{1}$ is not formally
  integrable, as the prolongation of the algebraic equation $u=x$ leads to the
  integrability conditions $u_x=1$ and $u_t=0$.  Projecting the prolonged
  equation gives
  \[
  \tilde{\mathcal{R}}_{1}^{(1)} \colon \
    \left\{
        \begin{array}{l}
        u_t  = v = v_t = w_x = w_t = v_x = 0 , \\
        u_x = w = 1 , \\
        u   = x  .
        \end{array}
    \right.
  \]
  Now the symbol vanishes, and so do the pull-backs of the contact forms:
  $\iota^*\omega^1=d\overline{x}-d\overline{x}=0$,
  $\iota^*\omega^2=d\overline{v}=0$,
  $\iota^*\omega^3=d\overline{w}=0$.
  Therefore we f\/ind
  $\V[\tilde{\mathcal{R}}_{1}^{(1)}]=\tilde{\mathcal{N}}_{1}^{(1)} \oplus \H =
  \{0\} \oplus \langle \partial_{\overline{x}},
  \partial_{\overline{t}}\rangle$.  As the Lie brackets of
  $\partial_{\overline{x}}$ and $\partial_{\overline{t}}$ trivially vanish,
  $T\tilde{\mathcal{R}}_{1}^{(1)} = \V[\tilde{\mathcal{R}}_{1}^{(1)}] = \H =
  \langle\partial_{\overline{x}}, \partial_{\overline{t}}\rangle$ is a
  two-dimensional involutive distribution.
\end{example}

Any $n$-dimensional complement $\H$ is obviously a transversal subdistribution
of $\V[\Rc{q}]$, but not necessarily involutive.  Conversely, any
$n$-dimensional subdistribution $\H$ of $\V[\Rc{q}]$ is a possible choice as a
complement.  Choosing a ``reference'' complement $\H_0$ with a basis $(X_i
\colon 1 \le i \le n)$, a basis for any other complement $\H$ arises by adding
some vertical f\/ields to the vectors $X_i$.  We will follow this approach in
the next section.  For the remainder of this section we turn our attention to
the choice of a convenient basis of $\V[\Rc{q}]$ that will facilitate our
computations.

Let $r:=\dim{\Nc{q}}$.  As an intersection of two involutive distributions,
the symbol $\Nc{q}$ is an involutive distribution, too.  Hence, there exists a
basis $(Y_k \colon 1 \le k \le r)$ for it whose Lie brackets vanish:
$[Y_k,Y_\ell]=0$ for all $1 \le k,\ell \le r$.  Since the vertical bundle
$V\pi^q_{q-1}$ is also involutive, we can decompose it into
\begin{displaymath}
	V\pi^q_{q-1}
	=
	\Nc{q} \oplus \mathcal{W},
\end{displaymath}
where $\mathcal{W}$ is again an involutive distribution.  It can be spanned by
vector f\/ields $W_1,W_2,\ldots,W_s$ (where $s = \sum_{k=1}^n \beta_q^{(k)}$
equals the number of principal derivatives) which are chosen such that we have
$[W_a,W_b]=0$ for all $1 \le a,b \le s$.  In local coordinates, a particularly
convenient choice for the f\/ields $Y_k$ and $W_a$ exists.  We f\/irst choose for
any $1 \le k \le r$ a \emph{parametric} derivative $u^\alpha_\mu$, that is
$(\alpha,\mu) \notin \mathcal{B}$, and set $Y_k := Y^\alpha_\mu :=
\iota_*(\partial_{\overline{u^\alpha_\mu}})$; then we choose for any $1 \le a
\le s$ a \emph{principal} derivative $u^\alpha_\mu$, that is $(\alpha,\mu) \in
\mathcal{B}$, and set $W_a := W^\alpha_\mu := \partial_{u^\alpha_\mu}$.

The reference complement $\H_0$ is chosen as follows.  Any basis of it must
consist of $n$ transversal contact f\/ields.  Since the f\/ields $C^\mu_\alpha$
are vertical, we can always use a basis $(\tilde{X}_1 \colon 1 \le i \le n)$
of the form
\begin{displaymath}
	\tilde{X}_i
	=
	C^{(q)}_i
	+
	\xi^\alpha_{i\mu} C^\mu_\alpha
\end{displaymath}
with some coef\/f\/icient functions $\xi^\alpha_{i\mu}$ chosen such that
$\tilde{X}_i$ is tangential to $\Rc{q}$.  The f\/ields $C^\mu_\alpha$ also span
the vertical bundle $V\pi^q_{q-1}$ and hence we may exploit the above
decomposition for a further simplif\/ication of the basis.  By subtracting from
each $\tilde{X}_i$ a suitable linear combination of the f\/ields $Y_k$ spanning
the symbol $\Nc{q}$, we arrive at a basis $(X_i \colon 1 \le i \le n)$ where
\begin{equation}\label{gl:Xi}
	X_i
	=
	C_i^{(q)}
	+
	\xi_i^a W_a .
\end{equation}

As already mentioned above, the Vessiot distribution $\V[\Rc{q}]$ is not
necessarily involutive.  So it is not surprising that its structure equations
are going to be important later.  We may extend the above chosen basis
$(X_i,Y_k)$ of $\V[\Rc{q}]$ to a basis of the derived Vessiot distribution,
$\V'[\Rc{q}]$, by adding vector f\/ields $Z_c$, $1 \le c \le
C:=\dim\V'[\Rc{q}]-\dim\V[\Rc{q}]$, where, using
(\ref{gl:LieOrdnung1kleiner}), for each $c$ we have some coef\/f\/icients
$\kappa^\alpha_{c \nu}$ such that $Z_c = \kappa^\alpha_{c \nu}
\partial_{u^\alpha_\nu}$ with $|\nu| = q-1$.  By construction, the
non-vanishing structure equations of $\V[\Rc{q}]$ take now the form
\begin{gather}\label{gl:Strukturgleichungen}
    [X_i,X_j] = \mathnormal\Theta_{ij}^cZ_c \qquad \text{and} \qquad
    [X_i,Y_k] = \mathnormal\Xi_{ik}^cZ_c
\end{gather}
for $1 \le i,j \le n$ and $1 \le k \le r$, with smooth functions
$\mathnormal\Theta_{ij}^c$ and $\mathnormal\Xi_{ik}^c$.  (For the complete set
of structure equations, we have to add $[Y_k,Y_l] = 0$ for $1 \le k,l \le r$.)

\begin{remark}\label{bem:blowup}
  Since the vector f\/ields $Z_c$, which appear on the right sides of the
  structure equations~(\ref{gl:Strukturgleichungen}), may, for $q=1$, span a
  proper subspace in $\langle \partial_{u^\alpha} \colon 1 \le \alpha \le m
  \rangle$, about the exact size of which we know nothing, we write them as
  linear combinations $Z_c =: \kappa_c^\alpha \partial_{u^\alpha}$.  The
  structure equations then become
  \setcounter{equation}{25}
\begin{subequations}\label{*****}
  \begin{gather}
        [X_i,X_j]
                = \mathnormal\Theta^c_{ij} \kappa_c^\alpha \partial_{u^\alpha}
        =: \mathnormal\Theta^\alpha_{ij} \partial_{u^\alpha}
        , \qquad 1 \le i,j \le n ,
        \tag{\ref{gl:Strukturgleichungen}a$^{\prime}$}\label{*****a}
        \\
    [X_i,Y_k]
        = \mathnormal\Xi^c_{ik}    \kappa_c^\alpha \partial_{u^\alpha}
        =: \mathnormal\Xi^\alpha_{ik} \partial_{u^\alpha}
        , \qquad 1 \le i \le n , \ \  1 \le k \le r  .
        \tag{\ref{gl:Strukturgleichungen}b$^{\prime}$}\label{*****b}
  \end{gather}
 \end{subequations}

  Knowing the larger sets of coef\/f\/icients $\mathnormal\Theta_{ij}^\alpha$,
  $\mathnormal\Xi_{ik}^\alpha$, we can reconstruct the true structure
  coef\/f\/icients $\mathnormal\Theta_{ij}^c$, $\mathnormal\Xi_{ik}^c$ by solving
  the overdetermined systems of linear equations
  \begin{equation*}
        \mathnormal\Theta^c_{ij} \kappa_c^\alpha
        =
        \mathnormal\Theta^\alpha_{ij}
                \qquad \text{and} \qquad
        \mathnormal\Xi^c_{ik}    \kappa_c^\alpha
        =
        \mathnormal\Xi^\alpha_{ik}.
  \end{equation*}
  This is always possible since the f\/ields $Z_c$ are assumed to be part of a
  basis for the derived Vessiot distribution $\V^\prime[\Rc{1}]$ and therefore
  linearly independent.  Thus there exist some coef\/f\/icient functions
  $\kappa_\alpha^{c}$ such that
  \begin{equation*}
        \mathnormal\Theta_{ij}^c
        =
        \mathnormal\Theta^\alpha_{ij}\kappa_\alpha^{c}
                \qquad \text{and} \qquad
        \mathnormal\Xi_{ik}^c
        =
        \mathnormal\Xi^\alpha_{ik}\kappa_\alpha^{c}
        .
  \end{equation*}
  For our later proof of an existence theorem for integral distributions, we
  will have to analyse certain matrices with the coef\/f\/icients
  $\mathnormal\Theta_{ij}^c$ and $\mathnormal\Xi_{ik}^c$ as their entries.  It
  turns out that this analysis becomes simpler, if we use the extended sets of
  coef\/f\/icients $\mathnormal\Theta_{ij}^\alpha$ and
  $\mathnormal\Xi_{ik}^\alpha$ instead.
\end{remark}

For a f\/irst-order equation $\Rc{1}$ with Cartan normal form (\ref{gl:CNF})
satisfying the assumptions of Proposition
\ref{prop:Existenz-H-notwendighinreichend} it is possible to perform this
process explicitly.  We choose as a reference complement $\H_0$ the linear
span of the vector f\/ields
\begin{equation}\label{gl:KomplementBasis}
        \iota_*\bar{X}_i
        =
        C^{(q)}_i
        +
        \sum_{(\alpha,\mu) \in \mathcal{B}}
                C^{(q)}_i(\phi^\alpha_\mu)
                C^\mu_\alpha
\end{equation}
for $1 \le i \le n$.  One verif\/ies in a straightforward computation that
(\ref{gl:KomplementBasis}) represents a valid choice (see \cite[Proposition
3.1.19]{Fesser:2008}).  Using this reference complement, we can explicitly
evaluate the Lie brackets (\ref{gl:Strukturgleichungen}) on $\Rc{1}$.  As we
are not able to determine a simple expression for the derived Vessiot
distribution, we follow the approach of Remark \ref{bem:blowup} and consider
the equations
(\ref{*****}$^{\prime}$) instead.

\begin{lemma}
  Let $i<j$, without loss of generality.  Then we obtain for the extended set
  of structure coefficients $\mathnormal\Theta^\alpha_{ij}$ in local
  coordinates on $\Rc{1}$ the following results:
  \begin{equation}\label{gl:VektorTheta}
  \mathnormal\Theta^\alpha_{ij} =
    \left\{
    \begin{array}{ll}
        C^{(1)}_i(\phi^\alpha_j) - C^{(1)}_j(\phi^\alpha_i), \qquad
            &
            (\alpha,i) \in \mathcal{B}
            \ \text{and} \
            (\alpha,j) \in \mathcal{B},
            \\
        C^{(1)}_i(\phi^\alpha_j),
            &
            (\alpha,i) \not\in \mathcal{B}
            \ \text{and} \
            (\alpha,j) \in \mathcal{B},
            \\
        0,
            &
            (\alpha,i) \not\in \mathcal{B}
            \ \text{and} \
            (\alpha,j) \not\in \mathcal{B}.
    \end{array}
    \right.
  \end{equation}
\end{lemma}

\begin{proof}
  By a rather straightforward computation; see \cite[pages 75 and
  76]{Fesser:2008}.
\end{proof}

We collect these coef\/f\/icients into vectors $\mathnormal\Theta_{ij}$ which have
$m$ rows each where the entries are ordered according to increasing $\alpha$.
It is useful to denote the symbol f\/ields
$Y_k=\iota_{*}(\partial_{\overline{u^{\beta}_{h}}})$ by using double
indices: $Y_k=Y^\beta_h$ for any $(\beta,h) \not\in \mathcal{B}$.

\begin{lemma}
  Set $\bar{Y}_k=:\bar{Y}^\beta_h$, then the extended set of structure
  coefficients $\mathnormal\Xi_{ik}^\alpha$ in local coordinates on $\Rc{1}$
  are
  \begin{equation}\label{gl:MatrixXi}
        \mathnormal\Xi^\alpha_{ik}
        =
        \left\{
        \begin{array}{ll}
                -C^h_\beta(\phi^\alpha_i),\qquad
                        &
                        (\alpha,i) \in \mathcal{B},
                        \\
                -1,
                        &
                        (\alpha,i) \notin \mathcal{B}
                        \ \text{and} \
                        (\alpha,i) = (\beta,h),
                        \\
                0,
                        &
                        (\alpha,i) \notin \mathcal{B}
                        \ \text{and} \
                        (\alpha,i) \not= (\beta,h).
        \end{array}
        \right.
  \end{equation}
\end{lemma}

\begin{proof}
  Again by a rather straightforward computation; see \cite[page 76]{Fesser:2008}.
\end{proof}

\begin{example}\label{bsp:Strukturgleichungen}
  We calculate the structure equations for Example \ref{bsp:Welle1}.  Here,
  $\bq{1}{1}=1$ and $\bq{1}{2}=3$.  We set
  $\bar{X}_3=:\bar{Y}_1=\bar{Y}^2_1$ and $\bar{X}_4=:\bar{Y}_2=\bar{Y}^3_1$.
  Then besides the trivial structure equations we get:
  \begin{gather*}
    [\bar{X}_1, \bar{X}_2] =0, \qquad
    [\bar{X}_1, \bar{Y}_1] =
    [\bar{X}_2, \bar{Y}_2] = -\partial_{\overline{v}} , \qquad
    [\bar{X}_2, \bar{Y}_1] =
    [\bar{X}_1, \bar{Y}_2] = -\partial_{\overline{w}} .
  \end{gather*}
  Noting that here the set $\mathcal{B}$ contains the pairs $(1,1) \equiv u_x$,
  $(1,2) \equiv u_t$, $(2,2) \equiv v_t$ and $(3,2) \equiv w_t$ while $(2,1) \equiv v_x$ and $(3,1) \equiv w_x$ are not in $\mathcal{B}$, one easily
  verif\/ies that all coef\/f\/icients are as given in equations~(\ref{gl:VektorTheta}) and~(\ref{gl:MatrixXi}).
\end{example}

We end this section with two technical remarks on how these coef\/f\/icients are
collected into matrices $\mathnormal\Xi_i$.  The examination of the ranks of
these matrices $\mathnormal\Xi_i$ is basic for the proof of the existence Theorem~\ref{thm:ExIntegralDist}
for integral distributions.

\begin{figure}[t]
\centering
\includegraphics{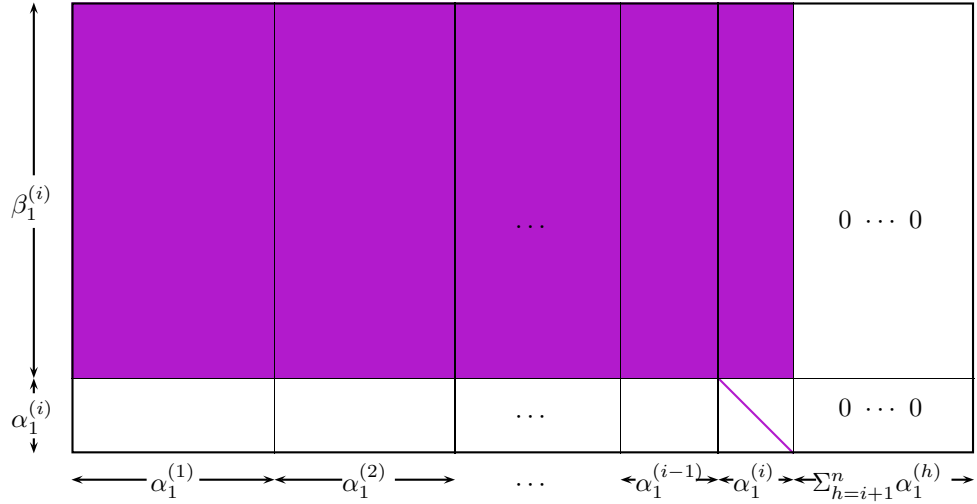}

    \caption{A sketch for the matrix $\mathnormal\Xi_i$ in equation
      (\ref{gl:Xi-Struktur}).}
    \label{bild:Xi-Struktur}
\end{figure}

\begin{remark}\label{bem:Subblocks}
  Some of the terms $-C^h_\beta(\phi^\alpha_i)$ where $(\alpha,i)\in\mathcal{B}$
  vanish, too: all of the parametric derivatives on the right side of an
  equation $\mathnormal\Phi^\alpha_i$ in the reduced Cartan normal form
  (\ref{gl:reducedCNF}) are of a class lower than that of the equation's left
  side as otherwise we would solve this equation for the derivative of highest
  class.  This means $-C^h_\beta(\phi^\alpha_i) = 0$ whenever
  $i=\cl(u^\alpha_i) < \cl(u^\beta_h)=h$.

  We collect the coef\/f\/icients $\mathnormal\Xi^\alpha_{ik}$ into matrices
  $\mathnormal\Xi_i$ using $i$ as the number of the matrix to which the entry
  $\mathnormal\Xi^\alpha_{ik}$ belongs, $\alpha$ as the row index of the entry
  and $k$ as its column index.  These matrices have $m$ rows each, ordered
  according to increasing $\alpha$; and they have $r=\dim{\Nc{1}}$ columns
  each of which can be labelled by pairs $(\beta,h)\not\in\mathcal{B}$ or the
  symbol f\/ields $\bar{Y}_k=\bar{Y}^\beta_h$.  More precisely, for $1 \le h \le
  n$, we set
  \begin{gather}\label{def:Xi-oben}
    \begin{pmatrix}
        -C^h_{\bq{1}{h}+1}(\phi^1_i)           & -C^h_{\bq{1}{h}+2}(\phi^1_i)
        & \cdots & -C^h_{m}(\phi^1_i) \\
        -C^h_{\bq{1}{h}+1}(\phi^2_i)           & -C^h_{\bq{1}{h}+2}(\phi^2_i)
        & \cdots & -C^h_{m}(\phi^2_i) \\
        \vdots                                 & \vdots
        & \ddots & \vdots \\
        -C^h_{\bq{1}{h}+1}(\phi^{\bq{1}{i}}_i) &
        -C^h_{\bq{1}{h}+2}(\phi^{\bq{1}{i}}_i) & \cdots &
        -C^h_{m}(\phi^{\bq{1}{i}}_i) \\
    \end{pmatrix}
    =: [\mathnormal\Xi_i]^h .
  \end{gather}
  Such a matrix with an upper index $h$ collects all those
  $\mathnormal\Xi^\alpha_{ik}$ into a block where $(\alpha,i) \in
  \mathcal{B}$.  For any $1 \le i \le n$, we have $m-\bq{1}{i} = \aq{1}{i}$,
  so such a matrix has $\bq{1}{i}$ rows and $\aq{1}{i}$ columns.  Since, for
  any $h$ with $i<h$ and for all $\bq{1}{h}+1 \le \beta \le m$, we have
  $-C^h_\beta(\phi^\alpha_i)=0$, such matrices $[\mathnormal\Xi_i]^h$ where
  $i<h$ are zero.  Furthermore, for $1 \le h \le n$, we assemble the remaining
  terms $\mathnormal\Xi^\alpha_{ik}$ (which are those where $(\alpha,i) \notin
  \mathcal{B}$) in a matrix.  As above, let $\bar{Y}_k=\bar{Y}^\beta_h$, and
  denote any entry $\mathnormal\Xi^\alpha_{ik}$ accordingly, for the sake of
  introducing the following matrix, by
  $\leftidx{_\beta^h}{\mathnormal\Xi}{_i^\alpha}$.  Now set
  \begin{gather*}
    \begin{pmatrix}
        \leftidx{_{\bq{1}{h}+1}^h}{\mathnormal\Xi}{_i^{\bq{1}{i}+1}} &
            \leftidx{_{\bq{1}{h}+2}^h}{\mathnormal\Xi}{_i^{\bq{1}{i}+1}} &
                \cdots &
                    \leftidx{_m^h}{\mathnormal\Xi}{_i^{\bq{1}{i}+1}} \\
        \leftidx{_{\bq{1}{h}+1}^h}{\mathnormal\Xi}{_i^{\bq{1}{i}+2}} &
            \leftidx{_{\bq{1}{h}+2}^h}{\mathnormal\Xi}{_i^{\bq{1}{i}+2}} &
                \cdots &
                    \leftidx{_m^h}{\mathnormal\Xi}{_i^{\bq{1}{i}+2}} \\
        \vdots                                                       &
            \vdots                                                       &
                \ddots &
                    \vdots \\
        \leftidx{_{\bq{1}{h}+1}^h}{\mathnormal\Xi}{_i^m}             &
            \leftidx{_{\bq{1}{h}+2}^h}{\mathnormal\Xi}{_i^m}             &
                \cdots &
                    \leftidx{_m^h}{\mathnormal\Xi}{_i^m} \\
    \end{pmatrix}
    =: [\mathnormal\Xi_i]_h .
  \end{gather*}
  For any $1 \le h \le n$, such a matrix with the index $h$ written below has
  $\aq{1}{i}$ rows and $\aq{1}{h}$ columns.  Let for any natural numbers $a$
  and $b$ denote $0_{a \times b}$ the $a \times b$ zero matrix.  According to
  equation~(\ref{gl:MatrixXi}), we have
  \begin{gather*}
    [\mathnormal\Xi_i]_h
     =
    \left\{
    \begin{array}{ll}
        -\mathbbm{1}_{\aq{1}{i}},     &     h=i, \\
        0_{\aq{1}{i}\times\aq{1}{h}}, & h\not=i .
    \end{array}
    \right.
  \end{gather*}
  Using the matrices $[\mathnormal\Xi_i]^h$ and $[\mathnormal\Xi_i]_h$ as
  blocks, we now build the matrix
  \begin{gather*}
    \mathnormal\Xi_i =
    \begin{pmatrix}
        [\mathnormal\Xi_i]^1 & [\mathnormal\Xi_i]^2 \cdots & [\mathnormal\Xi_i]^n \\
        [\mathnormal\Xi_i]_1 & [\mathnormal\Xi_i]_2 \cdots & [\mathnormal\Xi_i]_n
    \end{pmatrix} .
  \end{gather*}
  Taking into account what we have noted on its entries, this means
  \begin{gather}\label{gl:Xi-Struktur}
    \mathnormal\Xi_i
    =
        \begin{pmatrix}
                [\mathnormal\Xi_i]^1 & \cdots & [\mathnormal\Xi_i]^{i-1} &
        [\mathnormal\Xi_i]^i           & 0 \cdots 0 \\
        0                    & \cdots & 0                        & -{\mathbbm
        1}_{\alpha_1^{(i)}} & 0 \cdots 0
    \end{pmatrix}  .
  \end{gather}
  A sketch of such a matrix $\mathnormal\Xi_i$ where the entries which may be
  dif\/ferent from zero are marked as shaded areas and $-{\mathbbm
    1}_{\aq{1}{i}}$ as a diagonal line is given in Fig.~\ref{bild:Xi-Struktur}.
\end{remark}

For all $h$ where $1 \le h \le n$, we call the block $[\mathnormal\Xi_i]^h$ in
$\mathnormal\Xi_i$ stacked upon the block $[\mathnormal\Xi_i]_h$ in
$\mathnormal\Xi_i$ \emph{the $h$th block of columns in $\mathnormal\Xi_i$}.
For those $h$ with $\bq{1}{h}=m$ the $h$th block of columns is empty.  Now the
symbol f\/ields $\bar{Y}^\beta_h$, or equivalently the pairs $(\beta,h) \not\in
\mathcal{B}$, are used to order the $\dim \Nc{1}= r$ columns of
$\mathnormal\Xi_i$, according to increasing $h$ into $n$ blocks (empty for
those $h$ with $\aq{1}{h}=0$) and within each block according to increasing
$\beta$ (with $\bq{1}{j}+1 \le \beta \le m$).

\begin{remark}
  This means, the columns in $\mathnormal\Xi_i$ are ordered increasingly with
  respect to the term-over-position lift of the degree reverse lexicographic
  ranking.  Therefore, the entry $-C^h_{\bq{1}{h}+\gamma}(\phi^\alpha_i)$
  stands in the matrix $\mathnormal\Xi_i$ in line $\alpha$, in the $h$th block
  of columns of which it is the $\gamma$th one from the left.  Entries
  dif\/ferent from zero and from $-1$ may appear in $\mathnormal\Xi_i$ only in a
  $[\mathnormal\Xi_i]^h$ for $h \le i$.  To be exact, for any class $i$, the
  matrix $\mathnormal\Xi_i$ has $\aq{1}{i}$ rows where all entries are zero
  with only one exception: for each $1 \le \ell_i < \aq{1}{i}$ we have
  \begin{displaymath}
        \mathnormal\Xi_{i\,k}^{\beta_1^{(i)}+\ell_i}
        =
        -\delta_{\ell \, k},
  \end{displaymath}
  where $\ell := \sum_{h=1}^{i-1} \aq{1}{h} + \ell_i$.  The entries in the
  remaining upper $\beta_1^{(i)}$ rows are $-C^h_\beta(\phi^\alpha_i)$.  The
  potentially non-trivial ones of them are marked as shaded areas in Fig.~\ref{bild:Xi-Struktur}.
\end{remark}

Note that for a dif\/ferential equation with constant coef\/f\/icients all vectors
$\mathnormal\Theta_{ij}$ vanish and for a maximally over-determined equation
there are no matrices $\mathnormal\Xi_i$.

The unit block of $\aq{1}{i}$ rows, ${\mathbbm 1}_{\alpha_1^{(i)}}$, leads
immediately to the estimate
\begin{equation*}
        \alpha_1^{(i)}
        \le
        \rank \mathnormal\Xi_i
        \le
        \min{\left\{m,\sum_{h=1}^i \alpha_1^{(h)}\right\}}.
\end{equation*}

\begin{example}
  For Example \ref{bsp:Welle1}, we have $n=2$ and therefore two matrices
  $\mathnormal\Xi_1$ and $\mathnormal\Xi_2$.  Their entries follow immediately
  from our results in Example \ref{bsp:Strukturgleichungen}.  For $i=1$, we
  have
  \begin{gather*}
	\mathnormal\Xi_1	
	=
	\begin{pmatrix}
		 0 &  0
        \\
        -1 &  0
        \\
         0 & -1
    \end{pmatrix}
.
  \end{gather*}
  The f\/irst line is $[\mathnormal\Xi_1]^1 = (0,0)$.  The unit block below it
  is $[\mathnormal\Xi_1]_1$.  (We have $\bq{1}{1}=1$ and therefore
  $\aq{1}{1}=2$.  There is only $h=1$ to consider, as both parametric
  derivatives $v_x$ and $w_x$ are with respect to $x$ only.)  Because of the
  very simple nature of our system, we f\/ind here by accident that
  $\mathnormal\Xi_1 = \mathnormal\Xi_2$.  However, $[\mathnormal\Xi_2]^1 =
  \mathnormal\Xi_2$ is the whole matrix while there is no
  $[\mathnormal\Xi_2]_1$ because $\bq{1}{2}=m=3$ and therefore $\aq{1}{2}=0$.
  Finally, we obtain $\mathnormal\Theta_{12} = (v_x-w_t, 0, 0)^t$.  (The $t$
  top right marks the transpose.)  Note that its f\/irst entry is in fact the
  integrability condition.
\end{example}

\section{Flat Vessiot connections}

In this section, we develop an approach for constructing f\/lat Vessiot
connections which improves recent approaches
\cite{Fackerell:VessiotsVectorFieldFormulation,
  Stormark:LiesStructuralApproach, Vassiliou:VessiotStructure}: we exploit the
splitting of $\V[\Rc{q}]$ suggested by (\ref{gl:VessiotDistributionSymbolH})
to introduce convenient bases for integral distributions which yield structure
equations that are particularly simple.  Finally, we give necessary and
suf\/f\/icient conditions for Vessiot's approach to succeed.

Let $\Rc{q}$ locally be represented by the system
$\mathnormal\Phi^\tau(\xv,\uv^{(q)}) = 0$ where $1 \le \tau \le t$.  Our goal
is the construction of all $n$-dimensional transversal involutive
subdistributions $\U$ within the Vessiot distribution $\V[\Rc{q}]$.  Taking
for the Vessiot distribution the basis $(X_i,Y_k)$ where the vector f\/ields~$Y_k$ are the above mentioned basis of the symbol $T\iota(\Nc{q})$ with
vanishing Lie brackets and the vector f\/ields $X_i$ are given in (\ref{gl:Xi}),
we make for the basis $(U_i \colon 1 \le i \le n)$ of such a~subdistribution
$\U$ the ansatz
\begin{equation*}
    U_i = X_i + \zeta_i^kY_k
\end{equation*}
with yet undetermined coef\/f\/icient functions $\zeta_i^k \in
\mathcal{F}(\Rc{q})$.  This ansatz follows naturally from our considerations
so far, as the f\/ields $X_i$ are transversal to the f\/ibration over $\X$ and
span a reference complement to the symbol and in $\Nc{q}$ all f\/ields $Y_k$ are
vertical.

Since the f\/ields $U_i$ are in triangular form, the distribution $\U$ is
involutive if, and only if, their Lie brackets vanish, and using the structure
equations (\ref{gl:Strukturgleichungen}) this means:
\begin{gather}
		{[U_i,U_j]}
		=
			[X_i,X_j]
			+
			\zeta_i^k[Y_k,X_j]
			+
			\zeta_j^k[X_i,Y_k]
			+
			\bigl(U_i(\zeta_j^k) - U_j(\zeta_i^k)\bigr)Y_k
		\nonumber\\
\phantom{{[U_i,U_j]}}{} =
			\bigl(\mathnormal\Theta_{ij}^c
			-
			\mathnormal\Xi_{jk}^c\zeta_i^k
			+
			\mathnormal\Xi_{ik}^c\zeta_j^k\bigr)Z_c
			+
			\bigl(U_i(\zeta_j^k)
			-
			U_j(\zeta_i^k)\bigr)Y_k
		= 0.\label{gl:UU}
	\end{gather}
It follows from the def\/inition of the f\/ields $Y_k$ and $Z_c$ that they are
linearly independent and so their coef\/f\/icients must vanish for $\U$ to be
involutive.  Thus the Lie bracket (\ref{gl:UU}) yields two sets of conditions
for the coef\/f\/icient functions $\zeta_i^k$: a system of algebraic equations
\begin{gather}\label{gl:AlgebraischeBed}
	G_{ij}^c
	:=
	\mathnormal\Theta_{ij}^c
	-
	\mathnormal\Xi_{jk}^c\zeta_i^k
	+
	\mathnormal\Xi_{ik}^c\zeta_j^k
		 = 0 ,
			\qquad
						\left\{
				\begin{array}{l}
					1 \le c \le C, \\
					1 \le i < j \le n,
				\end{array}
			\right.
\end{gather}
and a system of dif\/ferential equations
\begin{gather}\label{gl:DifferentielleBed}
	H_{ij}^p
	:=
	U_i(\zeta_j^p)
	-
	U_j(\zeta_i^p)
		 = 0,
			  \qquad
						\left\{
				\begin{array}{l}
					1 \le p \le r, \\
					1 \le i < j \le n.
				\end{array}
			\right.
\end{gather}
In the algebraic system (\ref{gl:AlgebraischeBed}) the true structure
coef\/f\/icients $\mathnormal\Theta_{ij}^c$, $\mathnormal\Xi_{jk}^c$ appear.  For
our subsequent analysis we follow Remark \ref{bem:blowup} and replace them by
the extended set of coef\/f\/icients $\mathnormal\Theta_{ij}^\alpha$,
$\mathnormal\Xi_{jk}^\alpha$.  This corresponds to replacing
(\ref{gl:AlgebraischeBed}) by an equivalent but larger linear system of
equations which is, however, simpler to analyse.

Obviously, (\ref{gl:AlgebraischeBed}) is an inhomogeneous \emph{linear} system
to which any solution method for linear systems can be applied.  We shall see
in the next section that its structure (induced by the structure equations of
the Vessiot distribution) allows us to decompose it into simpler subsystems
which are considered step by step.  Many papers on the Vessiot theory (see for
example \cite{Fackerell:VessiotsVectorFieldFormulation,
  Stormark:LiesStructuralApproach}) study at this stage a \emph{quadratic}
system which \emph{only} by such a step-by-step approach can be reduced to a
series of linear problems.  The linearity of (\ref{gl:AlgebraischeBed}) is a
simple consequence of our choice of a~basis for $\V[\Rc{q}]$ which in turn
exploits the natural splitting $\V[\Rc{q}]=\Nc{q}\oplus\H$.

\begin{remark}\label{bem:Contraction}
  The vector f\/ields $Y_k$ lie in the Vessiot distribution $\V[\Rc{1}]$.  Thus,
  according to Proposition \ref{prop:intdist}, the subdistribution $\U$ is an
  integral distribution if, and only if, the coef\/f\/icients~$\zeta_k^i$ satisfy
  the algebraic conditions (\ref{gl:AlgebraischeBed}).  This observation
  permits us immediately to reduce the number of unknowns in our ansatz.
  Assume that we have values $1\leq i,j\leq n$ and $1\leq\alpha\leq m$ such
  that both $(\alpha,i)$ and $(\alpha,j)$ are not contained in $\mathcal{B}$,
  that is, $u^\alpha_i$ and $u^\alpha_j$ are both parametric derivatives (and
  thus obviously the second-order derivative $u^\alpha_{ij}$, too).  Then
  there exist two symbol f\/ields $Y_k=\iota_*(\partial_{u^\alpha_i})$ and
  $Y_l=\iota_*(\partial_{u^\alpha_j})$.  Now it follows from the coordinate
  form (\ref{gl:gammaq}) of the contact map that the subdistribution $\U$,
  spanned by the vector f\/ields $U_h = X_h + \zeta^k_h Y_k$ for $1 \le h \le
  n$, can be an integral distribution if, and only if,
  \begin{equation}\label{gl:zeta-Vektor-Beziehungen-1}
    \zeta^k_j=\zeta^l_i
    \qquad\text{or, equivalently,}\qquad
    \zeta_j^{(\alpha,i)}=\zeta_i^{(\alpha,j)}.
  \end{equation}
  for all $1 \le i < j \le n$ and $1 \le k,l \le r$.

  As the unknowns $\zeta_j^k$ may be understood as labels for the columns of
  the matrices $\mathnormal\Xi_h$, this identif\/ication leads to a contraction
  of these matrices.  The contracted matrices, denoted by
  $\hat{\mathnormal\Xi}_h$, arise as follows: whenever $\zeta^k_j=\zeta^l_i$
  then the corresponding columns of $\mathnormal\Xi_h$ are added; their sum
  replaces the f\/irst of these columns, while the second column is dropped.
  Similarly, we introduce reduced vectors $\hat\zeta_h$ where the redundant
  components are left out.  From now on, we always understand that in the
  equations above this reduction has been performed.  Otherwise some results
  would not be correct (see Example \ref{bsp:contraction} below).

  For a more thorough outline of these technical details, see \cite[Lemma
  3.3.5 and Remark 3.3.6]{Fesser:2008}.
\end{remark}

\begin{example}
  We consider a f\/irst-order equation in one dependent variable
  \begin{equation*}
    \Rc{1} \colon\ \left\{
      \begin{array}{l}
        u_n =\phi_n(\xv,u,u_1,\dots,u_{n-r-1}) ,\\
        \cdots\cdots\cdots\cdots\cdots\cdots \cdots \cdots\cdots \cdots \cdots\\
        u_{n-r} =\phi_{n-r}(\xv,u,u_1,\dots,u_{n-r-1}).
      \end{array}\right.
  \end{equation*}
  It is easy to show that the symbol of such an equation is always involutive
  and that the used coordinates are $\delta$-regular.  In order to simplify
  the notation, we use the following convention for the indices: $1\leq
  k,\ell<n-r$ and $n-r\leq a,b\leq n$.  Thus we may express our system in the
  concise form $u_a=\phi_a(\xv,u,u_k)$ and local coordinates on the
  submanifold $\Rc{1}$ are $(\overline{\xv},\overline{u},\overline{u}_k)$.

  The pull-back of the contact form $\omega=d u-u_id x^i$ generating the
  contact codistribution $\C_1^0$ is $\iota^*\omega=
  d\overline{u}-\overline{u}_kd\overline{x}^k-\phi_ad\overline{x}^a$ and
  $\V[\Rc{1}]=\langle
  \bar{X}_1,\dots,\bar{X}_n,\bar{Y}_1,\dots,\bar{Y}_{n-r-1}\rangle$ where
  \[
    \bar{X}_k=\partial_{\overline{x}^k}+
              \overline{u}_k\partial_{\overline{u}},\qquad
    \bar{X}_a=\partial_{\overline{x}^a}+\phi_a\partial_{\overline{u}},\qquad
    \bar{Y}_k=\partial_{\overline{u}_k}.
  \]
  The f\/ields $\bar{Y}_k$ span the symbol $\Nc{1}$ and the f\/ields $\bar{X}_i$
  our choice of a reference complement $\H_0$.  Setting
  $\bar{Z}=\partial_{\overline{u}}$, the structure equations of $\V[\Rc{1}]$
  are
	\begin{alignat*}{3}
    	&	[\bar{X}_k,\bar{X}_\ell]
      		=0,\qquad &&
		[\bar{Y}_k,\bar{Y}_\ell]
			=0,&
		\\
		& [\bar{X}_k,\bar{X}_a]
      		=\bar{X}_k(\phi_a)\bar{Z},\qquad
		&&
		[\bar{X}_a,\bar{X}_b]
      		=\bigl(\bar{X}_a(\phi_b)-\bar{X}_b(\phi_a)\bigr)\bar{Z},&
		\\
		&[\bar{X}_k,\bar{Y}_\ell]
      		=-\delta_{k\ell}\bar{Z},\qquad
		&&
		[\bar{X}_a,\bar{Y}_k]
      		=-\bar{Y}_k(\phi_a)\bar{Z}.&
	\end{alignat*}

  Now we make the above discussed ansatz $U_i=X_i+\zeta^k_iY_k$ for the
  generators of a transversal complement $\H$.  Modulo the Vessiot
  distribution $\V[\Rc{1}]$ we obtain for their Lie brackets
  \begin{subequations}
    \begin{gather}
      [U_k,U_\ell] \equiv(\zeta^\ell_k-\zeta^k_\ell)Z\mod\V[\Rc{1}],
           \label{eq:ukl}\\
      [U_a,U_k] \equiv
          \bigl(\zeta^k_a-\zeta^\ell_kY_\ell(\phi_a)-X_k(\phi_a)\bigr)Z
      \mod\V[\Rc{1}],\label{eq:ual}\\
      [U_a,U_b] \equiv
         \bigl(\zeta^k_aY_k(\phi_b)-\zeta^\ell_bY_\ell(\phi_a)+X_a(\phi_b)-X_b(\phi_a)\bigr)Z\mod\V[\Rc{1}].
\label{eq:uab}
    \end{gather}
  \end{subequations}

  The algebraic system (\ref{gl:AlgebraischeBed}) is now obtained by requiring
  that all the expressions in parentheses on the right hand sides vanish.  Its
  solution is straightforward.  The f\/irst subsystem (\ref{eq:ukl}) implies the
  equalities $\zeta^\ell_k=\zeta^k_\ell$.  This result was to be expected by
  the discussion in Remark~\ref{bem:Contraction}: both~$u_k$ and~$u_\ell$ are
  parametric derivatives for $\Rc{1}$ and thus we could have made this
  identif\/ication already in our ansatz for the complement.  The second
  subsystem (\ref{eq:ual}) yields that
  $\zeta^k_a=\zeta^\ell_kY_\ell(\phi_a)+X_k(\phi_a)$.  If we enter these
  results into the third subsystem (\ref{eq:uab}), then all unknowns $\zeta$
  drop out and the solvability condition
  \begin{displaymath}
    X_a(\phi_b)-X_b(\phi_a)+X_k(\phi_a)Y_k(\phi_b)-X_k(\phi_b)Y_k(\phi_a)=0
  \end{displaymath}
  arises.  Thus in this example the algebraic system
  (\ref{gl:AlgebraischeBed}) has a solution if, and only if, this condition is
  satisf\/ied.

  Comparing with the classical theory of such systems, one easily verif\/ies
  that this solvability condition is equivalent to the vanishing of the Mayer
  or Jacobi bracket $[u_a-\phi_a,u_b-\phi_b]$ on the submanifold $\Rc{1}$,
  which in turn is a necessary and suf\/f\/icient condition for the formal
  integrability of the dif\/ferential equation $\Rc{1}$.  Thus we may conclude
  that $\Rc{1}$ possesses $n$-dimensional integral distributions if, and only
  if, it is formally integrable (which in our case is also equivalent to
  $\Rc{1}$ being involutive).

  Thus here involution can be decided solely on the basis of the algebraic
  system (\ref{gl:AlgebraischeBed}).  We will show in the next section that
  this observation does not represent a special property of a very particular
  class of dif\/ferential equations but a general feature of the Vessiot theory.
\end{example}

\section{The existence theorem for integral distributions}
\label{section:ExThmIntegralDist}

Now the question arises, when the combined system (\ref{gl:AlgebraischeBed}),
(\ref{gl:DifferentielleBed}) has solutions?  We begin by analysing the
algebraic part (\ref{gl:AlgebraischeBed}).  We use for this analysis a
step-by-step approach originally proposed by Vessiot
\cite{Vessiot:Integration}.  Our analysis will automatically reveal necessary
and suf\/f\/icient assumptions for it to succeed.  As outlined in Remark~\ref{bem:Contraction}, we replace $\zeta_2$ by $\hat\zeta_2$ since for the
entries~$\zeta_2^{(\beta,1)}$ where $\bq{1}{2}+1 \le \beta \le m$ we know
already from equation (\ref{gl:zeta-Vektor-Beziehungen-1}) that
$\zeta_2^{(\beta,1)}=\zeta_1^{(\beta,2)}$.  Thus we begin the construction of
the integral distribution $\U$ by f\/irst choosing an arbitrary vector f\/ield
$U_1$ and then aiming for another vector f\/ield $U_2$ such that $[U_1,U_2] \in
T\iota(\V[\Rc{q}])$.  During the construction of the f\/ield $U_2$ we regard the
components of the vector $\zeta_1=\hat\zeta_1$ as given parameters and the
components of $\hat\zeta_2$ as the only unknowns of the system
\begin{gather}\label{gl:Xi-System-1.Schritt}
	\hat{\mathnormal\Xi}_1 \hat\zeta_2
	=
	\hat{\mathnormal\Xi}_2 \hat\zeta_1 - \mathnormal\Theta_{12} .
\end{gather}
Since the components of $\hat\zeta_1$ are not considered unknowns, the system
(\ref{gl:Xi-System-1.Schritt}) must not lead to any restrictions for the
coef\/f\/icients $\hat\zeta_1^k$.  Obviously, this is the case if, and only if,
\begin{gather}\label{gl:Rang-1.Schritt}
	\rank{\hat{\mathnormal\Xi}_1}
	=
	\rank{(\hat{\mathnormal\Xi}_1\ \ \hat{\mathnormal\Xi}_2)}.
\end{gather}
Assuming that (\ref{gl:Rang-1.Schritt}) holds, the system
(\ref{gl:Xi-System-1.Schritt}) is solvable if, and only if, it satisf\/ies
the augmented rank condition
\begin{equation}\label{gl:Rang-1.Schritt-Aug}
	\rank{\hat{\mathnormal\Xi}_1}
	=
	\rank{(\hat{\mathnormal\Xi}_1 \ \ \hat{\mathnormal\Xi}_2 \ \
         {-\mathnormal\Theta_{12}})}.
\end{equation}
When we have succeeded in constructing $U_2$, the next step is to
seek a further vector f\/ield $U_3$ such that $[U_1,U_3] \in
T\iota(\V[\Rc{q}])$ and $[U_2,U_3] \in T\iota(\V[\Rc{q}])$.  Now the
components of both vectors $\hat\zeta_1$ and $\hat\zeta_2$ are
regarded as given, and the components of $\hat\zeta_3$ are regarded
as the unknowns of the system
\begin{gather*}
	\hat{\mathnormal\Xi}_1
	\hat\zeta_3
	=
	\hat{\mathnormal\Xi}_3
	\hat\zeta_1
	-
	\mathnormal\Theta_{13}  ,
	\qquad
	\hat{\mathnormal\Xi}_2
	\hat\zeta_3
	=
	\hat{\mathnormal\Xi}_3
	\hat\zeta_2
	-
	\mathnormal\Theta_{23}  .
\end{gather*}
Now this system is not to restrict the components of both $\hat\zeta_1$ and
$\hat\zeta_2$ any further.  The interrelations between the $\hat\zeta_i$
following from the condition (\ref{gl:zeta-Vektor-Beziehungen-1}) for the
existence of integral distributions, given in Remark~\ref{bem:Contraction},
are taken care of by contracting $\mathnormal\Xi_3$ into
$\hat{\mathnormal\Xi}_3$.  This implies that now the rank condition
\begin{gather*}
	\mathrm{rank}
		\begin{pmatrix}
			\hat{\mathnormal\Xi}_1 \\
			\hat{\mathnormal\Xi}_2
		\end{pmatrix}
	=
	\mathrm{rank}
		\begin{pmatrix}
			\hat{\mathnormal\Xi}_1 & \hat{\mathnormal\Xi}_3 & 0 \\
			\hat{\mathnormal\Xi}_2 & 0 & \hat{\mathnormal\Xi}_3
		\end{pmatrix}
\end{gather*}
has to be satisf\/ied.  If it is, then for $1 \le c \le C =
\dim\V'[\Rc{q}]-\dim\V[\Rc{q}]$ the system
\begin{displaymath}
	\mathnormal\Theta_{13}^c
	-
	\hat{\mathnormal\Xi}_{3k}^c\zeta_1^k
	+
	\hat{\mathnormal\Xi}_{1k}^c\zeta_3^k
	= 0
 , \qquad
	\mathnormal\Theta_{23}^c
	-
	\hat{\mathnormal\Xi}_{3k}^c\zeta_2^k
	+
	\hat{\mathnormal\Xi}_{2k}^c\zeta_3^k
	= 0
\end{displaymath}
is solvable if, and only if, the augmented rank condition
\begin{displaymath}
	\mathrm{rank}
		\begin{pmatrix}
			\hat{\mathnormal\Xi}_1 \\
			\hat{\mathnormal\Xi}_2
		\end{pmatrix}
	=
	\mathrm{rank}
		\begin{pmatrix}
			\hat{\mathnormal\Xi}_1 & \hat{\mathnormal\Xi}_3 & 0 &
				-\mathnormal\Theta_{13} \\
			\hat{\mathnormal\Xi}_2 & 0 & \hat{\mathnormal\Xi}_3 &
				-\mathnormal\Theta_{23}
		\end{pmatrix}
\end{displaymath}
holds.  Now we proceed by iteration.  Given $j-1$ vector f\/ields $U_1,
U_2,\ldots,U_{j-1}$ of the required form spanning an involutive
subdistribution of $T\iota(\V[\Rc{1}])$, we construct the next vector f\/ield
$U_j$ by solving the system
\begin{gather}
		\hat{\mathnormal\Xi}_1 \hat\zeta_j
			=
			\hat{\mathnormal\Xi}_j \hat\zeta_1
			-
			\mathnormal\Theta_{1j}, \nonumber\\
	\cdots \cdots\cdots\cdots\cdots\cdots\cdots\label{gl:Schritti}\\
			\hat{\mathnormal\Xi}_{j-1} \hat\zeta_j
			=
			\hat{\mathnormal\Xi}_j \hat\zeta_{j-1}
			-
			\mathnormal\Theta_{j-1,j}.\nonumber
        \end{gather}
Again we consider only the components of the vector $\hat\zeta_j$ unknowns,
and the system (\ref{gl:Schritti}) must not imply any further restrictions on
the components of the vectors $\hat\zeta_i$ for $1 \le i<j$.  The
corresponding rank condition is
\begin{gather}
	\rank{
		\begin{pmatrix}
			\hat{\mathnormal\Xi}_1 \\
			\hat{\mathnormal\Xi}_2 \\
			\vdots \\
			\hat{\mathnormal\Xi}_{j-1}
		\end{pmatrix}}
	=
	\rank{
		\begin{pmatrix}
			\hat{\mathnormal\Xi}_1 & \hat{\mathnormal\Xi}_j &
				&  &    \\
			\hat{\mathnormal\Xi}_2 &     & \hat{\mathnormal\Xi}_j
				&  & 0  \\
				\vdots  & 0   &     & \ddots &    \\
			\hat{\mathnormal\Xi}_{j-1} & & & & \hat{\mathnormal\Xi}_j
		\end{pmatrix}}.
	\label{gl:Rangbed}
\end{gather}
Assuming that it holds, the equations (\ref{gl:Schritti}) are solvable and
yield solutions for the components of $\hat\zeta_j$ if, and only if, the system
satisf\/ies the augmented rank condition
\begin{gather}
	\rank{
		\begin{pmatrix}
			\hat{\mathnormal\Xi}_1 \\
			\hat{\mathnormal\Xi}_2 \\
			\vdots \\
			\hat{\mathnormal\Xi}_{j-1}
		\end{pmatrix}}
	=
	\rank{
		\begin{pmatrix}
			\hat{\mathnormal\Xi}_1 & \hat{\mathnormal\Xi}_j
 				& & & & -\mathnormal\Theta_{1j} \\
			\hat{\mathnormal\Xi}_2 & & \hat{\mathnormal\Xi}_j
				& & 0 & -\mathnormal\Theta_{2j} \\
				\vdots  & 0   &     & \ddots &     & \vdots  \\
			\hat{\mathnormal\Xi}_{j-1} &  & &
                                & \hat{\mathnormal\Xi}_j
				& -\mathnormal\Theta_{j-1,j}
		\end{pmatrix}}.
	\label{gl:RangbedAugmented}
\end{gather}

\begin{remark}
  We can simplify our calculations for a dif\/ferential equation $\Rc{q}$,
  represented by $u^\alpha_\mu =
  \phi^\alpha_\mu(\mathbf{x},\mathbf{u},\hat{\mathbf{u}}^{(q)})$ where each
  equation is solved for a principal derivative (all pairwise dif\/ferent)
  $u^\alpha_\mu$ with $|\mu|=q$, and where $\hat{\mathbf{u}}^{(q)}$ denotes
  the set of the remaining, thus parametric, derivatives of order $q$ and
  less.  (Example \ref{bsp:Welle1} discusses a system of that kind.)  Now we
  can use the local coordinates on $\Rc{q}$.  For the generators of the symbol
  $\Nc{q}$ we may choose for $(\alpha,\mu)\not\in\mathcal{B}$ the vector
  f\/ields $\bar{Y}^\alpha_\mu = \partial_{\overline{u^\alpha_\mu}}$; as a basis
  for the complement $\H \subseteq \V[\Rc{q}]$, we can take the vector f\/ields
  \[
        \bar{X}_i
        =
        \partial_{\overline{x^i}}
        +
        \sum_{\alpha = 1}^m
                \sum_{0 \le |\mu| < q}
                        \overline{u^\alpha_{\mu+1_i}}
                                \partial_{\overline{u^\alpha_\mu}}  ,
  \]
  which satisfy equation (\ref{gl:KomplementBasis}).  For the vector f\/ields
  $W_a$ which appear in equation (\ref{gl:Xi}) we may choose the contact
  vector f\/ields $C_\alpha^\mu$ where $(\alpha,\mu) \in \mathcal{B}$.  The
  further procedure, solving f\/irst the structure equations, which now take the
  simple form (\ref{gl:Strukturgleichungen}), to f\/ind the generators of
  $\V^\prime(\Rc{q})$ and then solving equation (\ref{gl:AlgebraischeBed}) for
  the coef\/f\/icient functions $\zeta_j$, is the same as in the case where the
  representation is not in the above solved form.
\end{remark}

\begin{example}\label{bsp:noninvol}
  Before we proceed with our theoretical analysis, let us demonstrate with a
  concrete dif\/ferential equation that in general we cannot expect that the
  above outlined step-by-step construction of integral distributions works.
  In order to keep the size of the example reasonably small, we use the
  second-order equation\footnote{It could always be rewritten as a f\/irst-order
    one satisfying the assumptions made above, and the phenomenon we want to
    discuss is independent of this transformation.} $\Rc{2}$ def\/ined by the
  system $u_{xx}=\alpha u$ and $u_{yy}=\beta u$ with two real constants
  $\alpha$, $\beta$.  Its symbol $\Nc{2}$ is \emph{not} involutive.  However,
  one easily proves that $\Rc{2}$ is formally integrable for arbitrary choices
  of the constants $\alpha$, $\beta$.  One readily computes that the Vessiot
  distribution $\V[\Rc{2}]$ is generated by the following three vector f\/ields
  on $\Rc{2}$:
  \begin{gather*}
      X_1=\partial_{\overline{x}}+
           \overline{u_x}\partial_{\overline{u}}+
           \alpha\overline{u}\partial_{\overline{u_x}}+
           \overline{u_{xy}}\partial_{\overline{u_y}},\qquad
      X_2=\partial_{\overline{y}}+
           \overline{u_y}\partial_{\overline{u}}+
           \overline{u_{xy}}\partial_{\overline{u_x}}+
           \beta\overline{u}\partial_{\overline{u_y}},\qquad
      Y_1=\partial_{\overline{u_{xy}}}.\!
  \end{gather*}
  They yield as structure equations for $\V[\Rc{2}]$:
  \begin{gather*}
    [X_1,X_2]=\beta\overline{u_x}\partial_{\overline{u_y}}-
              \alpha\overline{u_y}\partial_{\overline{u_x}},\qquad
    [X_1,Y_1]=-\partial_{\overline{u_y}},\qquad
    [X_2,Y_1]=-\partial_{\overline{u_x}}.
  \end{gather*}

  For the construction of a two-dimensional integral distribution
  $\U\subset\V[\Rc{2}]$ we make as above the ansatz $U_i=X_i+\zeta_i^1Y_1$
  with two coef\/f\/icients~$\zeta_i^1$.  As we want to perform a step-by-step
  construction, we assume that we have chosen some f\/ixed value for~$\zeta_1^1$
  and try now to determine~$\zeta_2^1$ such that
  $[U_1,U_2]\equiv0\bmod\V[\Rc{2}]$.  Evaluation of the Lie bracket yields the
  equation
  \begin{equation}\label{eq:U1U2}
    [U_1,U_2]\equiv
        (\beta\overline{u_x}-\zeta_2^1)\partial_{\overline{u_y}}-
        (\alpha\overline{u_y}-\zeta_1^1)\partial_{\overline{u_x}}
        \mod\V[\Rc{2}].
  \end{equation}
  A necessary condition for the vanishing of the right side is that
  $\zeta_1^1=\alpha\overline{u_y}$.  Hence one cannot choose this coef\/f\/icient
  arbitrarily, as we assumed, but~(\ref{eq:U1U2}) determines \emph{both}
  functions $\zeta_i^1$ uniquely.  Note that the conditions on the
  coef\/f\/icients $\zeta_i^1$ imposed by~(\ref{eq:U1U2}) are trivially solvable
  and thus~$\Rc{2}$ possesses an integral distribution (which is even
  involutive).  The problem is that it could not be constructed systematically
  with the above outlined \emph{step-by-step} process.
\end{example}

The following theorem links the satisfaction of the rank conditions
(\ref{gl:Rangbed}) and (\ref{gl:RangbedAugmented}), and thus the solvability
of the algebraic system (\ref{gl:AlgebraischeBed}) by the above described
step-by-step process, with intrinsic properties of the dif\/ferential equation
$\Rc{q}$ and its symbol $\Nc{q}$.  It represents an existence theorem for
integral distributions.
\begin{theorem}\label{thm:ExIntegralDist}
  Assume that $\delta$-regular coordinates have been chosen for the local
  representation of the differential equation $\Rc{q}$.  Then the rank
  condition \eqref{gl:Rangbed} is satisfied for all $1\leq j\leq n$ if, and
  only if, the symbol $\Nc{q}$ is involutive.  The augmented rank condition
\eqref{gl:RangbedAugmented} holds for all $1\leq j\leq n$ if, and only if,
  the differential equation $\Rc{q}$ is involutive.
\end{theorem}

The proof of Theorem \ref{thm:ExIntegralDist} is given in Appendix
\ref{anhang:ExThmIntegralDist}, as parts of it are rather technical.  First,
we transform the dif\/ferential equation into an equivalent f\/irst-order system
with a representation in the reduced Cartan normal form (\ref{gl:reducedCNF})
and the same Cartan characters.  This is always possible; see, for example,
\cite[Proposition A.3.1]{Seiler:Buch}.  The basic idea of the proof is then
that the rank condition (\ref{gl:Rangbed}) is equivalent to the vanishing of
the obstructions to involution of the symbol in Lemma \ref{lem:Monster}, and
the augmented rank condition (\ref{gl:RangbedAugmented}) is equivalent to the
vanishing of the remaining integrability conditions (recall that Lemma
\ref{lem:Monster} gives all obstructions to involution only if the used
coordinates are $\delta$-regular so that this assumption is necessary in
Theorem \ref{thm:ExIntegralDist}).  The concrete implementation of this idea
requires some technical considerations concerning the transformation of the
matrices (\ref{gl:Rangbed}) and (\ref{gl:RangbedAugmented}) into row echelon
form, working out their contractions and analysing the interrelation between
these operations.

\begin{example}
  We demonstrate the role of $\delta$-regularity of the coordinates with the
  wave equation $u_{xy}=0$ in characteristic coordinates which are not
  $\delta$-regular.  A straightforward computation yields as generators for
  the Vessiot distribution the f\/ields $X_1=C_1^{(2)}$, $X_2=C_2^{(2)}$,
  $Y_1=\partial_{u_{xx}}$ and $Y_2=\partial_{u_{yy}}$.  Following our
  approach, we make the ansatz $U_1=X_1+\zeta_1^1Y_1+\zeta_1^2Y_2$ and
  $U_2=X_2+\zeta_2^1Y_1+\zeta_2^2Y_2$.  Evaluation of the Lie bracket
  $[U_1,U_2]$ yields now the algebraic equations $\zeta_1^2=0$ and
  $\zeta_2^1=0$ and the dif\/ferential equations
  $U_1(\zeta_2^1)-U_2(\zeta_1^1)=0$ and $U_1(\zeta_2^2)-U_2(\zeta_1^2)=0$.  As
  in Example \ref{bsp:noninvol}, we see that the step-by-step process is
  broken, since we obtain a condition on the coef\/f\/icient $\zeta_1^2=0$ which
  we should be able to consider a parameter here.
\end{example}

At this point, we have proven that integral distributions within the Vessiot
distribution exist if, and only if, the algebraic conditions
(\ref{gl:AlgebraischeBed}) are solvable, and that this is equivalent to the
augmented rank condition (\ref{gl:Rangbed}) being satisf\/ied.  This in turn is
the case precisely if the dif\/ferential equation is involutive.  Thus we have
characterised the existence of Vessiot connections for Vessiot's
\cite{Vessiot:Integration} step-by-step approach.

\section{The existence theorem for f\/lat Vessiot connections}
\label{section:ExThmFlatVConn}

There remains to analyse the solvability, if we add the dif\/ferential system
(\ref{gl:DifferentielleBed}).  Its solvability is equivalent to the existence
of \emph{flat} Vessiot connections in that each f\/lat Vessiot connection of~$\Rc{1}$ corresponds to a solution of the combined system
(\ref{gl:AlgebraischeBed}), (\ref{gl:DifferentielleBed}).  We f\/irst note that
the set of dif\/ferential conditions (\ref{gl:DifferentielleBed}) alone is again
an involutive system.

\begin{proposition}\label{prop:DiffBedInvolutiv}
  The differential conditions \eqref{gl:DifferentielleBed} alone represent an
  involutive differential equation of first order.
\end{proposition}

The proof, which is somewhat technical, is given in Appendix
\ref{anhang:DiffBedInvolutiv}. It is based on an evaluation of the Jacobi
identity for the Lie brackets of the f\/ields $U_i$.

If the original equation $\Rc{1}$ is analytic, then the quasi-linear system
(\ref{gl:DifferentielleBed}) is analytic, too.  Thus we may apply the
Cartan--K\"ahler theorem (Theorem \ref{thm:Cartan-Kaehler}) to it which
guarantees the existence of solutions.  The problem is that the combined
system (\ref{gl:AlgebraischeBed}), (\ref{gl:DifferentielleBed}) is not
necessarily involutive, as the prolongation of the algebraic equations
(\ref{gl:AlgebraischeBed}) may lead to additional dif\/ferential equations.
Instead of analysing the ef\/fect of these integrability conditions, we proceed
as follows.  If we assume that~$\Rc{1}$ is involutive, then we know from
Theorem \ref{thm:ExIntegralDist} that the algebraic equations
(\ref{gl:AlgebraischeBed}) are solvable.  Actually, the proof of the Theorem
\ref{thm:ExIntegralDist} produces an explicit row echelon form of the system
matrix given on the right hand side of equation (\ref{gl:Rangbed}).  Now we
use the interrelations between the unknowns $\zeta_\ell^{(\alpha,h)}$, where
$(\alpha,h) \not\in \mathcal{B}$, to eliminate in the dif\/ferential conditions
(\ref{gl:DifferentielleBed}) some of them by expressing them as linear
combinations of the remaining ones.  Let \mbox{$1 \le i < j \le n$}; then for all
$(\alpha,i)$ where $\bq{1}{i}+1 \le \alpha \le \bq{1}{j}$ we f\/ind
\begin{gather}\label{gl:zeta-rel}
	\zeta_j^{(\alpha,i)}
	=
	\sum_{k=1}^j
		\sum_{\gamma=\bq{1}{k}+1}^m
			C^k_\gamma
				(\phi_j^\alpha)
					\zeta_i^{(\gamma,k)}
			+
			C_i^{(1)}(\phi^\alpha_j)
	 ,
\end{gather}
while for all $(\alpha,i)$ where $\bq{1}{j}+1 \le \alpha \le m$ the
correspondence (\ref{gl:zeta-Vektor-Beziehungen-1}), given in Remark
\ref{bem:Contraction}, holds again (see \cite[Corollary 3.3.23]{Fesser:2008}
for the detailed calculation).  In this way, we plug the algebraic conditions
into the dif\/ferential conditions and thus get rid of them.
\begin{example}
For Example \ref{bsp:Welle1}, where $i=1$ and $j=2$, we deduce for
\[
	U_1 = X_1 + \zeta_1^{v_x} \partial_{v_x} + \zeta_1^{w_x} \partial_{w_x}
	\qquad \text{and} \qquad
	U_2 = X_2 + \zeta_2^{v_x} \partial_{v_x} + \zeta_2^{w_x} \partial_{w_x}
\]
from the algebraic condition (\ref{gl:Xi-System-1.Schritt}) that
$-\zeta_2^{v_x} = -\zeta_1^{v_x}$ and $-\zeta_2^{w_x} = -\zeta_1^{w_x}$.  This
is equation~(\ref{gl:zeta-rel}) for $(\alpha,i) \equiv v_x$ and $(\alpha,i)
\equiv w_x$.  The rank conditions (\ref{gl:Rang-1.Schritt}) and
(\ref{gl:Rang-1.Schritt-Aug}) are satisf\/ied because $v_x-w_t=0$, so that the
matrix $(\mathnormal\Xi_1 \ \mathnormal\Xi_2 \ -\mathnormal\Theta_{12})$ has a
vanishing f\/irst row.  (There is no contraction of matrices for $j=2$.)  There
are no interrelations like those in equation
(\ref{gl:zeta-Vektor-Beziehungen-1}) because there is no index value~$\alpha$
for which both $\iota_*(\partial_{u^\alpha_1})$ and
$\iota_*(\partial_{u^\alpha_2})$ would be symbol f\/ields.  The dif\/ferential
conditions~(\ref{gl:DifferentielleBed}) are
\[
	U_1 (\zeta_2^{v_x}) - U_2(\zeta_1^{v_x}) = 0
	 , \qquad
	U_1 (\zeta_2^{w_x}) - U_2(\zeta_1^{w_x}) = 0
	 ,
\]
and by substituting $\zeta_2^{v_x} = \zeta_1^{v_x}$ and $\zeta_2^{w_x} =
\zeta_1^{w_x}$, we can drop the algebraic conditions from the system.  Since
$n=2$, there is no further procedure.
\end{example}

We can now prove the following existence theorem for f\/lat Vessiot connections.
\begin{theorem}\label{thm:CombinedSolvable}
  Assume that\/ $\delta$-regular coordinates have been chosen for the local
  representation of the analytic differential equation~$\Rc{1}$.  Then the
  combined system \eqref{gl:AlgebraischeBed}, \eqref{gl:DifferentielleBed} is
  solvable.
\end{theorem}

The idea of the proof is this: following the strategy we have just outlined,
we eliminate some of the unknowns $\hat\zeta^k_i$.  Because of the simple
structure of (\ref{gl:DifferentielleBed}), it turns out that we must take a~closer look only at those equations where the leading derivative is of one of
the unknowns we eliminate.  A somewhat lengthy but straightforward computation
shows that these equations actually vanish.  The remaining equations still
form an involutive system.  Thus we eventually arrive at an analytic
involutive dif\/ferential equation for the coef\/f\/icient functions $\hat\zeta^k_i$
which is solvable according to the Cartan--K\"ahler theorem.  The details are
worked out in Appendix~\ref{anhang:ExThmFlatVConn}.

\begin{example}\label{bsp:contraction}
  Consider the f\/irst-order equation
  \begin{equation*}
	\Rc{1} \colon \
		\left\{
			\begin{array}{l}
				u_t=v_t=w_t=u_s=0 , \qquad
				v_s=2u_x+4u_y , \\
				w_s=-u_x-3u_y , \qquad
				u_z=v_x+2w_x+3v_y+4w_y
			\end{array}
		\right.
  \end{equation*}
  in the f\/ive independent variables $x$, $y$, $z$, $s$, $t$ and the three dependent
  variables $u$, $v$, $w$.  It is formally integrable, and its symbol is involutive
  with $\dim{\Nc{1}}=8$.  Thus $\Rc{1}$ is an involutive equation.  For the
  matrices $\mathnormal\Xi_i$, all of which are $3 \times 8$-matrices, we f\/ind
  \begin{alignat*}{3}
	& \mathnormal\Xi_1
		=
		\left(
			\begin{smallmatrix}
            -1 & 0 & 0 & 0 & 0 & 0 & 0 & 0 \\
            0 & -1 & 0 & 0 & 0 & 0 & 0 & 0 \\
            0 & 0 & -1 & 0 & 0 & 0 & 0 & 0
			\end{smallmatrix}
		\right), \qquad
		&&
	\mathnormal\Xi_2
		=
		\left(
			\begin{smallmatrix}
            0 & 0 & 0 & -1 & 0 & 0 & 0 & 0 \\
            0 & 0 & 0 & 0 & -1 & 0 & 0 & 0 \\
            0 & 0 & 0 & 0 & 0 & -1 & 0 & 0
			\end{smallmatrix}
		\right),&
	\\
	& \mathnormal\Xi_3
		=
		\left(
			\begin{smallmatrix}
            0 & -1 & -2 & 0 & -3 & -4 &  0 &  0 \\
            0 &  0 &  0 & 0 &  0 &  0 & -1 &  0 \\
            0 &  0 &  0 & 0 &  0 &  0 &  0 & -1
			\end{smallmatrix}
		\right), \qquad
		&&
    \mathnormal\Xi_4
		=
		\left(
			\begin{smallmatrix}
             0 & 0 & 0 &  0 & 0 & 0 & 0 & 0 \\
            -2 & 0 & 0 & -4 & 0 & 0 & 0 & 0 \\
             1 & 0 & 0 &  3 & 0 & 0 & 0 & 0
			\end{smallmatrix}
		\right), \qquad
	\mathnormal\Xi_5
		\
		&= \,
		0_{3 \times 8} .&
  \end{alignat*}

  For the f\/irst two steps in the construction of the f\/ields $U_i$, the rank
  conditions are trivially satisf\/ied even for the non-contracted matrices.
  But not so in the third step where we have in the row echelon form of the
  arising $9 \times 32$-matrix in the 7th row zero entries throughout except
  in the 12th column (where we have $-2$) and in the 17th column (where we
  have $2$).  As a consequence, we obtain the equality $\zeta_1^4=\zeta_2^1$
  and the rank condition for this step does not hold.  However, since both
  $u_x$ and $u_y$ are parametric derivatives and in our ordering
  $Y_1=\iota_*(\partial_{u_x})$ and $Y_4=\iota_*(\partial_{u_y})$, this
  equality is already taken into account in our reduced ansatz and for the
  matrices $\hat{\mathnormal\Xi}_i$ the rank condition is satisf\/ied.

  Note that the rank condition is f\/irst violated when the rank reaches the
  symbol dimension (which is $8$).  From then on, the rank of the left matrix
  in (\ref{gl:Rangbed}) stagnates at $\dim{\Nc{1}}$ while the rank of the
  augmented matrix rises further.  The entries breaking the rank condition
  dif\/fer by their sign, while their corresponding coef\/f\/icients in Lemma
  \ref{lem:Monster} are collected into one sum and thus vanish.  Here it shows
  that contracting the matrices, as explained in Remark \ref{bem:Contraction},
  is necessary.
\end{example}

\section{Conclusions}

Vessiot's \cite{Vessiot:Integration} original motivation for the introduction
of his theory was to provide an alternative proof of the Cartan--K\"ahler
theorem (the same holds for Stormark's presentation
\cite{Stormark:LiesStructuralApproach} of the theory).  If one takes this
point of view, then one may say that our proof of Theorem~\ref{thm:CombinedSolvable} is a~``cheat'', as it uses the Cartan--K\"ahler
theorem instead of proving it.  However, in our opinion, this point of view is
the main reason why the corresponding proofs in
\cite{Stormark:LiesStructuralApproach, Vessiot:Integration} are so dif\/f\/icult
to read (and actually incomplete as they neglect that the involutivity of the
system is a necessary condition).  The direct proof of the Cartan--K\"ahler
theorem (for dif\/ferential equations) given in textbooks on the formal theory
like \cite{Pommaret:SystemsOfPDEs,Seiler:Buch} is much simpler and more
transparent; in particular, it makes clear where involution (as opposed to
mere formal integrability) is needed.  Thus, if the only goal consisted of
proving such an existence and uniqueness theorem, then there would be no need
to bother with Vessiot's theory.

We believe that the real value of Vessiot's theory lies in the fact that it
provides via the distribution $\V[\Rc{q}]$ an additional geometric structure
on an involutive dif\/ferential equation $\Rc{q}$ which is very useful for the
\emph{further} analysis of the equation, that is, after its solvability has
been established.  One possible application is the investigation of certain
forms of singular behaviour of solutions.  In fact, the classical works of
Arnold and collaborators (see \cite{Arnold:GeoMethODE} for an elementary
introduction) on implicit ordinary dif\/ferential equations are based on the
Vessiot distribution (without using this terminology); for some further works
also in the context of partial dif\/ferential equations or numerical analysis
see, for example, \cite{Lychagin:Singularities, Tuomela:1997, Tuomela:1998}.
The basic idea here is that such behaviour mainly stems from the fact that at
some points the considered involutive distributions cease to be transversal to
the f\/ibration $\Rc{q}\rightarrow\X$ so that integral manifolds can no longer
be interpreted as prolonged solutions in the classical sense.  However, it
often makes sense to consider them a generalised form of (potentially
multi-valued) solutions.

Another application concerns symmetry theory.  In classical symmetry theory
one always considers the dif\/ferential equation $\Rc{q}$ as a submanifold of~$\je{q}$ and looks then for dif\/feomorphisms of $\je{q}$ which are $(i)$~compatible with the contact structure of $\je{q}$ and $(ii)$~leave $\Rc{q}$
invariant.  It was only fairly late realized \cite{AKO:Symmetries} that this
approach yields only what is now called external symmetries and that in
general further useful symmetries may exist.  Using Vessiot's theory, we may
represent a dif\/ferential equation as the pair $(\Rc{q},\V[\Rc{q}])$, that is,
as a manifold together with a distribution on it~-- without any recourse to an
ambient space (all relevant properties of~$\je{q}$ are captured in the Vessiot
distribution~$\V[\Rc{q}]$).  This point of view yields automatically all
internal symmetries (and is equivalent to the symmetry theory of exterior
dif\/ferential systems).

We mentioned already in the Introduction that Vessiot's theory takes an
intermediate position between the formal theory and the theory of exterior
dif\/ferential systems.  It allows for the transfer of ideas from the latter to
the former one without the need to rewrite a dif\/ferential equation as an
exterior system.  For example, \cite[Section 9.5]{Seiler:Buch} contains a
formulation of prolongation structures without exterior forms using the
Vessiot distribution instead\footnote{Both symmetry theory and prolongation
  structures were also discussed by Fackerell
  \cite{Fackerell:VessiotsVectorFieldFormulation} via Vessiot theory.}.  It
should also be possible to give an \emph{explicit} proof of the equivalence of
the formal theory and the Cartan--K\"ahler theory on the basis of the results
presented here.

Finally, we comment on the practical side of the Vessiot theory.  Throughout
this article, we made a number of assumptions on the treated dif\/ferential
equation: we assumed that it is f\/irst rewritten as a f\/irst-order system, then
that all present algebraic equations are explicitly solved and f\/inally that
each equation is solved for its principal derivative.  While, from a
theoretical point of view, these operations are always possible under fairly
mild regularity assumptions, usually they cannot be performed ef\/fectively.
However, the made assumptions are only needed in order to be able to prove
results like Lemma~\ref{lem:Monster} on the explicit form of the obstructions
to involution.  These expressions are already bad enough with our simplifying
assumptions; for more general systems they would hardly be manageable.

A concrete dif\/ferential equation need not be transformed to such a special
normal form in order to apply the Vessiot theory.  Even for fully
implicit equations of arbitrary order, the determination of all integral
distributions requires only linear algebra and is easily implemented in a
computer algebra system.  In fact, already Fackerell
\cite{Fackerell:VessiotsVectorFieldFormulation} reported of an implementation
of Vessiot's theory.  Within a recent project thesis
\cite{Globke:Studienarbeit}, Globke implemented large parts of the theory in
\textsf{MuPAD}.  The problems and limitations are here exactly the same as in
implementations of the Cartan--K\"ahler theory (see, for example,
\cite{HartleyTucker:ConstructiveCartanKaehler}).

\appendix

\vspace{-1pt}

\section{Proof of Theorem \ref{thm:ExIntegralDist}}
\label{anhang:ExThmIntegralDist}

In this section of the appendix, we prove Theorem \ref{thm:ExIntegralDist},
the existence theorem for integral distributions.  It is in principle by
straightforward matrix calculation and involves a tedious distinction of
several cases and subcases: we have to compare the entries in the matrix on
the right hand side of equation (\ref{gl:RangbedAugmented}) for step $j$ after
turning it into row echelon form with the integrability conditions and the
obstructions to involution as they are given in Lemma \ref{lem:Monster}.  For
this purpose, we f\/ix $1 \le i < j$ and consider the block
$\hat{\mathnormal\Xi}_i$ and the entries to its right in the complete matrix.

Since for a dif\/ferential equation with an involutive symbol the
obstructions to involution va\-nish
and for an involutive dif\/ferential equation the integrability
conditions vanish, too, it follows that the augmented rank condition,
stated in equation (\ref{gl:RangbedAugmented}), is equivalent to the
dif\/ferential equation being involutive, which in turn is the case if,
and only if, the algebraic conditions (\ref{gl:AlgebraischeBed}) are
satisf\/ied, which is necessary and suf\/f\/icient for the existence of
integral distributions within the Vessiot distribution.

If the order of the dif\/ferential equation is $q>1$, transform it into an
equivalent f\/irst-order equation; for the details of this procedure, see
\cite[Subsection 2.5.1]{Fesser:2008} or \cite{Seiler:Buch}.  If $\Rc{q}$ is
involutive, then so is $\Rc{1}$ (see \cite{Seiler:Buch} for a straightforward
proof of this).  We assume this f\/irst-order equation is represented in reduced
Cartan normal form (\ref{gl:reducedCNF}).  We f\/irst prove the rank condition
(\ref{gl:Rangbed}) for the homogeneous system.  (The proof for the augmented
rank condition (\ref{gl:RangbedAugmented}) follows.) We proceed in two steps:
We consider the \emph{complete matrix} at the $j$th step,
\begin{equation}\label{def:CompleteMatrix}
	\begin{pmatrix}
		\mathnormal\Xi_1     & \mathnormal\Xi_j &                  &        &    \\
		\mathnormal\Xi_2     &                  & \mathnormal\Xi_j &        & 0  \\
		\vdots               & 0                &                  & \ddots &    \\
		\mathnormal\Xi_{j-1} &                  &                  &        & \mathnormal\Xi_j
	\end{pmatrix}.
\end{equation}
(It dif\/fers from the matrix on the right hand side of
equation (\ref{gl:Rangbed}) in that it is not yet contracted.)
The complete matrix at the $j$th step is built from $(j-1)j$ blocks:
the stack of $j-1$ matrices~$\mathnormal\Xi_i$, $1 \le i \le j-1$, on
the left, with each of the $\mathnormal\Xi_i$ having another $j-1$
blocks to its right, the $i$th of which being $\mathnormal\Xi_j$ and
all the others being zero.

For easier reference, let, for $1 \le i,k \le j-1$, be $[i,k]$ the
$k$th block right of a $\mathnormal\Xi_i$.  Then, for all $1 \le i,k \le
j-1$, we have
\[
[i,k] =
    \left\{
        \begin{array}{ll}
            \mathnormal\Xi_j, & i=k, \\
            0_{m \times r},   & i\not=k .
        \end{array}
    \right.
\]
For convenience, we set $[i,0]:=\mathnormal\Xi_i$.  Let, for $1 \le g \le h \le n$,
\begin{equation}\label{def:Xi-Bloecke_g_bis_h}
    \mbl{[\mathnormal\Xi_i]}{}{}{}{g\dots h}
\end{equation}
denote the matrix that results from writing the block matrices
$\mbl{[\mathnormal\Xi_i]}{}{}{}{g}$,
$\mbl{[\mathnormal\Xi_i]}{}{}{}{g+1}$, \ldots, $\mbl{[\mathnormal\Xi_i]}{}{}{}{h}$ (from left to right) next to each
other.  Let $\mbl{[i,k]}{}{}{}{h}$ be the $h$th block of columns in $[i,k]$,
and let, for $1 \le g \le h \le n$, in analogy to the shorthand
(\ref{def:Xi-Bloecke_g_bis_h}),
\begin{equation*}
\mbl{[i,k]}{}{}{}{g\dots h}
\end{equation*}
denote the matrix that results from writing the matrices
$\mbl{[i,k]}{}{}{}{g}$, $\mbl{[i,k]}{}{}{}{g+1}$, \ldots,
$\mbl{[i,k]}{}{}{}{h}$ (from left to right) next to each other.

If $M$ is any $a\times b$-matrix, then let $\mbl{[M]}{c}{d}{}{}$ be
the matrix made from the \emph{rows} with indices (that is, labels)
$c$ to $d$, $\mbl{[M]}{}{}{e}{f}$ the matrix made from the
\emph{columns} with indices $e$ to $f$ and $\mbl{[M]}{c}{d}{e}{f}$
the matrix made from the entries in the rows with indices $c$ to $d$
and in the columns with indices from $e$ to $f$.  If the columns of
$M$ are grouped into blocks and $g$ denotes which blocks are meant,
then we write $\mbl{[M]}{1}{c}{}{g}$ (with $g$ up right) to show that
the f\/irst $c$ upper rows are being selected, and we write
$\mbl{[M]}{d}{b}{g}{}$ (with $g$ below right) to show that the last
$b-d+1$ lower rows are being selected.  For the block matrix made from
$M$ by selecting the rows indexed $c$ to $d$ within the block of
columns labelled $g$, we write $\mbl{[M]g}{c}{d}{}{}$.  This notation
is redundant in that the position of~$g$ does not give new
information, but in the calculations to come it increases
readability.

We f\/irst turn the complete matrix into row echelon form.  Then we
contract columns and consider the ef\/fect.

The matrices $\mathnormal\Xi_i$ are shown in equation (\ref{gl:Xi-Struktur})
and have a block-structure of the form
\begin{gather*}
    \begin{pmatrix}
        a_{\alpha r} & b_{\alpha s} & c_{\alpha t} \\
        0_{S \times R} & -\mathbbm{1}_S & d_{st}
    \end{pmatrix}
     =:  (g_{ij} \colon 1 \le i \le A+S , 1 \le j \le R+S+T),
\end{gather*}
where $A$, $R$, $S$, $T$ are adequate natural numbers or zero, and the blocks are
\begin{alignat*}{3}
    &\begin{pmatrix}
        a_{\alpha r} \colon 1 \le \alpha \le A , 1 \le r \le R
    \end{pmatrix} ,\qquad &
    &\begin{pmatrix}
        b_{\alpha s} \colon 1 \le \alpha \le A , 1 \le s \le S
    \end{pmatrix} , &
    \\
    &\begin{pmatrix}
        c_{\alpha t} \colon 1 \le \alpha \le A , 1 \le t \le T
    \end{pmatrix} ,&
    &\begin{pmatrix}
        d_{st} \colon 1 \le s \le S , 1 \le t \le T
    \end{pmatrix} , &
\end{alignat*}
the $S \times S$ unit matrix $\mathbbm 1_S$ or the zero matrix $0_{S \times R}$.
(Here the index $r$ is just some index and
not supposed to be the dimension of any symbol.) Then the
substitution
\begin{gather}\label{lem:substitution}
        \begin{pmatrix}
                a_{\alpha r} & b_{\alpha s} & c_{\alpha t}
        \end{pmatrix}
         \leftarrow
        \begin{pmatrix}
                a_{\alpha r} & b_{\alpha s} & c_{\alpha t}
        \end{pmatrix}
        +
        (b_{\alpha s}) \cdot
        \begin{pmatrix}
                0_{S \times R} & -\mathbbm{1}_S & d_{st}
        \end{pmatrix}
\end{gather}
transforms the matrix $(g_{ij})$ into
\begin{gather*}
	\begin{pmatrix}
		a_{\alpha r} & 0_{A \times S} & c_{\alpha t}+\sum_{s=1}^S b_{\alpha s}d_{s t} \\
		0_{S \times R} & -\mathbbm{1}_S & d_{st}
	\end{pmatrix}
	 .
\end{gather*}
For the transformation of the complete matrix into row echelon form, we use
this obvious method.  Let $1 < j \le n$ be given.  Choose $1 \le i < j$.  Now
we eliminate all the non-trivial entries in the rows
\[
	\mbl{[\mathnormal\Xi_i \, [i,1] \, [i,2]\dots[i,j-1]]}{1}{\bq{1}{i}}{}{} \ .
\]
We turn to the block
$\mbl{[\mathnormal\Xi_i]}{1}{\bq{1}{i}}{}{}$ f\/irst.  Non-trivial
entries therein may appear according to equation~(\ref{gl:Xi-Struktur})
only for $1 \le h \le i$ in the blocks~$\mbl{[\mathnormal\Xi_i]}{1}{\bq{1}{i}}{}{h}$.  So we f\/ix $1 \le h \le
i$ and consider the block~$\mbl{[\mathnormal\Xi_i]}{1}{\bq{1}{i}}{}{h}$.  Now set
\begin{subequations}
\begin{gather}
    (a_{\alpha r})
         = \mbl{[\mathnormal\Xi_i]}{1}{\bq{1}{i}}{}{1\dots h-1} ,
        \nonumber
        \\
    (b_{\alpha s})
         = \mbl{[\mathnormal\Xi_i]}{1}{\bq{1}{i}}{}{h} ,
        \nonumber
        \\
    (c_{\alpha t})
         = \mbl{[[\mathnormal\Xi_i]^{h+1\dots n}[i,1][i,2]\dots [i,j-1]]}{1}{\bq{1}{i}}{}{} ,
        \label{alg:cBlock}
        \\
    (d_{st})
         = \mbl{[[\mathnormal\Xi_h]_{h+1\dots n}[h,1][h,2]\dots [h,j-1]]}{\bq{1}{h}+1}{m}{}{} ,
        \label{alg:dBlock}
        \\
    -\mathbbm{1}_S
         = -\mathbbm{1}_{\aq{1}{h}}  ,
        \nonumber
        \\
    0_{S \times R}
         = 0_{\aq{1}{h}\times\sum_{l=1}^{h-1}\aq{1}{l}}  .
        \nonumber
\end{gather}
\end{subequations}
Then it follows that the substitution (\ref{lem:substitution}) leaves
$(a_{\alpha r})$ unchanged, turns $(b_{as})$ into zero as required
and makes $(c_{\alpha t})$ into $(c_{\alpha t}+\sum_{s=1}^S b_{\alpha
s}d_{s t})$.  According to (\ref{alg:dBlock}) the entries
$\sum_{s=1}^S b_{\alpha s}d_{st}$
here are of two types: from
$\mbl{[\mathnormal\Xi_h]}{\bq{1}{h}+1}{m}{h+1\dots n}{}$, the left part
of $(d_{st})$, we have the entries in $(b_{\alpha s}) \cdot
\mbl{[\mathnormal\Xi_h]}{\bq{1}{h}+1}{m}{h+1\dots n}{}$, and from
$\mbl{[[h,1][h,2]\dots [h,j-1]]}{\bq{1}{h}+1}{m}{}{}$, the right part of
$(d_{st})$, we have those of $(b_{\alpha s}) \cdot
\mbl{[h,l]}{\bq{1}{h}+1}{m}{}{}$ where $1 \le l \le j-1$.

Consider the $c_{\alpha t}+\sum_{s=1}^S b_{\alpha s}d_{s t}$ where
the factors $d_{s t}$ are from the left part of $(d_{st})$: according
to equation (\ref{gl:Xi-Struktur}),
$\mbl{[\mathnormal\Xi_h]}{\bq{1}{h}+1}{m}{h+1\dots n}{} = 0$, so
\begin{gather*}
    (b_{\alpha s}) \cdot \mbl{[\mathnormal\Xi_h]}{\bq{1}{h}+1}{m}{h+1\dots n}{}
        = 0   .
\end{gather*}
It follows $c_{\alpha t}+\sum_{s=1}^S b_{\alpha s}d_{s t} = c_{\alpha
t}$, so that all entries in
$\mbl{[\mathnormal\Xi_i]}{1}{\bq{1}{i}}{}{h+1\dots n}$, the left part of
$(c_{st})$, remain unchanged through the elimination of $(b_{\alpha
s}) = \mbl{[\mathnormal\Xi_i]}{1}{\bq{1}{i}}{}{h}$.

For the rest of the proof, we consider the remaining entries
$c_{\alpha t}+\sum_{s=1}^S b_{\alpha s}d_{s t}$ ,where the factors
$d_{s t}$ are from
$\mbl{[[h,1][h,2]\dots [h,j-1]]}{\bq{1}{h}+1}{m}{}{}$, the right part of
$(d_{st})$.  Fix $1 \le l \le j-1$ and consider $(b_{\alpha s}) \cdot
\mbl{[h,l]}{\bq{1}{h}+1}{m}{}{}$.

There are two cases: $h=l$ and $h\not=l$.  For $h\not=l$ we have
$\mbl{[h,l]}{}{}{}{} = 0$ according to the def\/inition of the complete matrix
(\ref{def:CompleteMatrix}), thus $\mbl{[h,l]}{\bq{1}{h}+1}{m}{}{} = 0$ and so
\begin{gather*}
    (b_{\alpha s}) \cdot \mbl{[h,l]}{\bq{1}{h}+1}{m}{}{}
          = 0   .
\end{gather*}
Again it follows $c_{\alpha t}+\sum_{s=1}^S b_{\alpha s}d_{s t} =
c_{\alpha t}$, so that for $h\not=l$ all entries in
\begin{displaymath}
	\mbl{[[i,1][i,2]\dots [i,j-1]]}{1}{\bq{1}{i}}{}{}  ,
\end{displaymath}
the right part of $(c_{\alpha t})$, remain unchanged through the
elimination of $(b_{\alpha s}) =
\mbl{[\mathnormal\Xi_i]}{1}{\bq{1}{i}}{}{h}$, too.

For $h=l$, we have $[h,h]=\mathnormal\Xi_j$.  The structure of
$\mathnormal\Xi_j$, which is (\ref{gl:Xi-Struktur}) with $i$ replaced by
$j$, implies that non-vanishing entries are possible in the blocks
$\mbl{[h,h]}{\bq{1}{h}+1}{m}{k}{}$ where $1 \le k \le j$.  In fact
\begin{subequations}\label{gl:reason}
\begin{gather}
    \mbl{[h,h]}{\bq{1}{j}+1}{m}{k}{} =0 \qquad \text{for} \ \  1 \le k \le j-1 \  \ \ \text{and} \\
    \mbl{[h,h]}{\bq{1}{j}+1}{m}{j}{} =-\mathbbm{1}_{\aq{1}{j}}.
\end{gather}
\end{subequations}
The entries $\sum_{s=1}^S b_{\alpha s}d_{st}$ for those $d_{st}$
within $\mbl{[h,h]}{\bq{1}{h}+1}{m}{k}{}$ are
\begin{gather}
    (b_{\alpha s}) \cdot \mbl{[h,h]}{\bq{1}{h}+1}{m}{k}{}
         = \mbl{[\mathnormal\Xi_i]}{1}{\bq{1}{i}}{}{h} \cdot \mbl{[\mathnormal\Xi_j]}{\bq{1}{h}+1}{m}{k}{}
        = \left(\sum_{s=1}^S C_{\bq{1}{h}+s}^h(\phi^\alpha_i) C_{\bq{1}{k}+t}^k(\phi^{\bq{1}{h}+s}_j)\right)  .
        \label{gl:Eintrag1}
\end{gather}
This matrix has $A=\bq{1}{i}$ rows and $T=\aq{1}{k}$ columns.  We
consider its entry in row~$\alpha$ and column $t$.  Setting $\gamma
:= \bq{1}{h}+s$, $\delta := \bq{1}{k}+t$ and using $S=\aq{1}{h}$
in~(\ref{gl:Eintrag1}), this entry is
\begin{gather}
    \sum_{\gamma=\bq{1}{h}+1}^{\bq{1}{h}+\aq{1}{h}} C_\gamma^h(\phi^\alpha_i) C_{\delta}^k(\phi^\gamma_j)
    \label{gl:Eintrag2}
\end{gather}
for all $1 \le k \le j$.  Some of the $\bq{1}{h}+\aq{1}{h}=m$ summands
vanish
because of the special entries~(\ref{gl:reason}).  As a consequence,
from $\bq{1}{j}+1$ on, of all the summands $C^h_\gamma(\phi^\alpha_i)
C^k_\delta(\phi^\gamma_j)$ in (\ref{gl:Eintrag2}) at most one remains:
none for $k\not=j$, exactly one for $k=j$, namely the one for
$\gamma=\bq{1}{j}+t=\delta$.  Using the Kronecker-delta, for all $\bq{1}{j}+1 \le \delta \le m$ we
have
\begin{gather*}
    C^h_\delta(\phi^\alpha_i) C^j_\delta(\phi^\delta_j)  = \delta_{kj} \cdot C^h_\delta(\phi^\alpha_i) .
\end{gather*}
Now (\ref{gl:Eintrag2}) becomes for all $1 \le \alpha \le \bq{1}{i}$
and all $\bq{1}{k}+1 \le \delta \le m$
\begin{gather}\label{gl:Allgemein1}
    \sum_{\gamma=\bq{1}{h}+1}^{\bq{1}{j}}
        C_\gamma^h(\phi^\alpha_i)
            C_\delta^k(\phi^\gamma_j)
        + \delta_{kj}
        \cdot
        C^h_\delta(\phi^\alpha_i)  .
    \tag{\ref{gl:Eintrag2}$^\prime$}
\end{gather}
These are the terms $\sum_{s=1}^S b_{\alpha s}d_{st}$ in the case
$h=l$ when $1 \le k \le j$.  Now for the subcase $j+1 \le k \le n$:
these blocks $\mbl{[h,h]}{\bq{1}{h}+1}{m}{k}{} =
\mbl{[\mathnormal\Xi_j]}{\bq{1}{h}+1}{m}{k}{}$, again on account of
the structure of $\mathnormal\Xi_j$, which is~(\ref{gl:Xi-Struktur})
with $i$ replaced by $j$, are zero (if they exist at all; they do not
for $k$ with $\aq{1}{k}=0$).  So in this case we have
\begin{gather*}
    (b_{\alpha s}) \cdot \mbl{[\mathnormal\Xi_j]}{\bq{1}{h}+1}{m}{j+1\dots n}{}
         = 0   ,
\end{gather*}
which contains the terms $\sum_{s=1}^S b_{\alpha s}d_{st}=0$.  As a
consequence, here, in the case $h=l$, the terms $c_{\alpha
t}+\sum_{s=1}^S b_{\alpha s}d_{s t}$ are of the following form.  The
$c_{\alpha t}$ of interest in line (\ref{alg:cBlock}) are those in
the block $\mbl{[i,h]}{1}{\bq{1}{i}}{}{{}}$.  For $h\not=i$, we have
$[i,h]=0$, thus $\mbl{[i,h]}{1}{\bq{1}{i}}{}{}=0$.  Non-trivial
entries $c_{\alpha t}$ are possible only for $h=i$; in that case
$[i,i]=\mathnormal\Xi_j$, thus
$\mbl{[i,h]}{1}{\bq{1}{i}}{}{}=\mbl{[\mathnormal\Xi_j]}{1}{\bq{1}{i}}{}{}$.
According to the structure of $\mathnormal\Xi_j$, which is
(\ref{gl:Xi-Struktur}) with $i$ replaced by $j$, the only entries
$c_{\alpha t}$ within the block
$\mbl{[\mathnormal\Xi_j]}{1}{\bq{1}{i}}{}{{}}$ which may not vanish
are, for $1 \le k \le j$, those of the form
$-C^k_\delta(\phi_j^\alpha)$.  They make up the block
$\mbl{[\mathnormal\Xi_j]}{1}{\bq{1}{i}}{}{1\dots j}=\mbl{[i,i]}{1}{\bq{1}{i}}{}{1\dots j}$.

So for any row index $\alpha$ and any column index
$\delta=\bq{1}{k}+t$ we have $c_{\alpha t}=-\delta_{hi}
C^k_\delta(\phi_j^\alpha)$ as the most general form of an entry.

Since (\ref{gl:Allgemein1}) is $\sum_{s=1}^Sb_{\alpha s}d_{st}$, it
follows that $c_{\alpha t}+\sum_{s=1}^Sb_{\alpha s}d_{st}$ is
\begin{gather}
    -\delta_{hi} C^k_\delta(\phi_j^\alpha)
        + \sum_{\gamma=\bq{1}{h}+1}^{\bq{1}{j}}
            C_\gamma^h(\phi^\alpha_i)
            C_\delta^k(\phi^\gamma_j)
                + \delta_{kj} C^h_\delta(\phi^\alpha_i) .
    \label{gl:Allgemein2}
\end{gather}
Since $1 \le \alpha \le \bq{1}{i}$ and $\bq{1}{k}+1 \le \delta \le
m$, any term in the row echelon form of the complete matrix without a
$-1$ in a negative unit block somewhere to its left has this form or
is a zero in $\mbl{[\mathnormal\Xi_i]}{1}{\bq{1}{i}}{}{i+1\dots n}$ for
some $i<j$ (in which case it remains zero throughout the elementary
row transformations and so does not inf\/luence the rank).  This means, if
all these expressions vanish (when contracted), the rank condition is
satisf\/ied.  To show that they do vanish (when contracted) for a system
with an involutive symbol, we consider them as the new entries in
$\mbl{[i,h]}{1}{\bq{1}{i}}{}{}$ with $1 \le h \le j$.  Now with regard
to the relation between $h$ and $i$ there are three cases: $h>i$,
$h=i$ and $h<i$.  We consider them in that order.
\begin{enumerate}\itemsep=0pt
\item[1.] Let $h>i$.  Then according to the structure of
$\mathnormal\Xi_i$, (\ref{gl:Xi-Struktur}), all the
$C_\gamma^h(\phi^\alpha_i)$  and $C_\delta^h(\phi^\alpha_i)$ in
(\ref{gl:Allgemein2}) vanish.  Since $h\not=i$, $\delta_{hi}=0$, thus all
of (\ref{gl:Allgemein2}) vanishes.
%
%
\item[2.] Let $h=i$.  We shall consider several subcases for $h=i$,
and for $h<i$ after that, which may be labelled by second-order
derivatives $u^\delta_{hk}$ where $h$ and $k$ are combined as shown in Fig.~\ref{bild:FaelleMatrix}.
For f\/ixed $\delta$, any $u^\delta_{hk}$ belongs, according to its indices $h$
and $k$, to exactly one of the blocks marked by the case denominations 2.(a),
2.(b) and 3.(a) to 3.(d2).

\begin{figure}
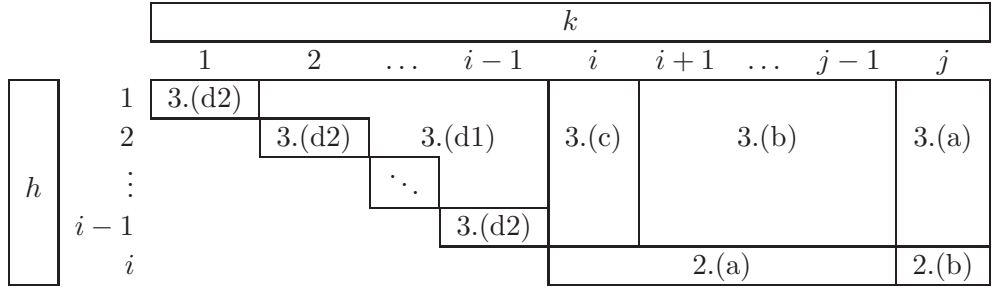

\centering
    \extrarowheight1pt
    \begin{tabular}{crccccccccc}
        \cline{3-11}
           &
           & \multicolumn{9}{|c|}{$k$} \\
        \cline{3-11}
           &
           & $1$ & $2$ &\ldots& $i-1$ & $i$ & $i+1$ &\ldots& $j-1$ & $j$  \\
    \cline{1-1} \cline{3-11}
    \multicolumn{1}{|c|}{} &
    1      & \multicolumn{1}{|c|}{3.(d2)}
           & \multicolumn{3}{c}{}
           & \multicolumn{1}{|c|}{}
           & \multicolumn{3}{c|}{}
           & \multicolumn{1}{c|}{} \\
    \cline{3-4}
    \multicolumn{1}{|c|}{} &
    2      &
           & \multicolumn{1}{|c|}{3.(d2)}
           & \multicolumn{2}{c}{3.(d1)}
           & \multicolumn{1}{|c|}{3.(c)}
           & \multicolumn{3}{c|}{3.(b)}
           & \multicolumn{1}{c|}{3.(a)} \\
    \cline{4-5}
    \multicolumn{1}{|c|}{$h$} &
    \vdots &
           &
           & \multicolumn{1}{|c|}{$\ddots$}
           &
           & \multicolumn{1}{|c|}{}
           & \multicolumn{3}{c|}{}
           & \multicolumn{1}{c|}{} \\
    \cline{5-6}
    \multicolumn{1}{|c|}{} &
    $i-1$  &
           &
           &
           & \multicolumn{1}{|c}{3.(d2)}
           & \multicolumn{1}{|c|}{}
           & \multicolumn{3}{c|}{}
           & \multicolumn{1}{c|}{} \\
    \cline{6-11}
    \multicolumn{1}{|c|}{} &
    $i$    &
           &
           &
           & \multicolumn{1}{c}{}
           & \multicolumn{4}{|c|}{2.(a)}
           & \multicolumn{1}{c|}{2.(b)} \\
    \cline{1-1} \cline{7-11}
    \end{tabular}
    \caption{Possible combinations of $h$ and $k$ in the terms $u^\delta_{hk}$
      used as labels for columns in the contracted matrix.  Each term
      $u^\delta_{hk}$ corresponds to one subcase in the consideration of cases
      2. $h=i$ and 3. $h<i$.  Each of these subcases refers to one of the
      terms in square brackets in Lemma \ref{lem:Monster} as follows:
      2.(a)~-- line~(\ref{line:Monster-ik}), 2.(b)~-- line~(\ref{line:Monster-ij}), 3.(a)~-- line~(\ref{line:Monster-hj}),
      3.(b)~-- line~(\ref{line:Monster-hk2}), 3.(c)~-- line~(\ref{line:Monster-hi}), 3.(d1)~-- line~(\ref{line:Monster-hk1}),
      3.(d2)~-- line~(\ref{line:Monster-hh}).}
    \label{bild:FaelleMatrix}
\end{figure}

It turns out that (\ref{gl:Allgemein2}) is the common form of all the sums
that appear in the squared brackets of Lemma \ref{lem:Monster} as the
coef\/f\/icients of second-order derivatives $u^\delta_{hk}$.  Not only can the
case distinctions of the following argument, 2.(a), 2.(b) and 3.(a) to 3.(d2),
be labelled by these $u^\delta_{hk}$, but in fact a case labelled
$u^\delta_{hk}$ is dealt with by using the fact that according to Lemma
\ref{lem:Monster} the coef\/f\/icient in square brackets of that same
$u^\delta_{hk}$ vanishes for an involutive system.  The correspondence of the
cases and the terms $u^\delta_{hk}$ is given in the caption of Fig.~\ref{bild:FaelleMatrix}.

The cases $h=i$ with subcase $k < i$ and $h<i$ with subcase $k < i$ need not
be considered because the columns of the complete matrix are to be contracted
when proving the rank condition.  This contraction concerns those columns of
the complete matrix which are labelled by the same second-order derivative,
and each second-order derivative is used exactly once in the argument.
\item[(a)] For the f\/irst subcase of $h=i$ let $i \le k < j$.  Then
(\ref{gl:Allgemein2}) becomes
\begin{gather*}
    -C^k_\delta(\phi_j^\alpha)+\sum_{\gamma=\bq{1}{i}+1}^{\bq{1}{j}}
    C_\gamma^i(\phi^\alpha_i)
    C_\delta^k(\phi^\gamma_j)
\end{gather*}
where $\bq{1}{k}+1 \le \delta \le m$.  This vanishes for a system with
an involutive symbol according to Lemma \ref{lem:Monster}, line
(\ref{line:Monster-ik}).

\item[(b)] For the second subcase of $h=i$ choose $k=j$.  Then
(\ref{gl:Allgemein2}) becomes
\begin{gather*}
    -C^j_\delta(\phi_j^\alpha)+\sum_{\gamma=\bq{1}{i}+1}^{\bq{1}{j}}
    C_\gamma^i(\phi^\alpha_i) C_\delta^j(\phi^\gamma_j) + C^i_\delta(\phi^\alpha_i)
\end{gather*}
where $\bq{1}{j}+1 \le \delta \le m$.  According to Lemma \ref{lem:Monster},
line (\ref{line:Monster-ij}), this vanishes for a system with an
involutive symbol.
%
%
\item[3.] Let $h<i$.  We consider several subcases with regard to the
relation between~$h<i$,~$j$ and~$k$.

\item[(a)] First choose $k=j$.  Then (\ref{gl:Allgemein2}) becomes
\begin{gather*}
    \sum_{\gamma=\bq{1}{h}+1}^{\bq{1}{j}} C_\gamma^h(\phi^\alpha_i)
    C_\delta^j(\phi^\gamma_j) +
    C^h_\delta(\phi^\alpha_i)
\end{gather*}
where $\bq{1}{j}+1 \le \delta \le m$.  According to Lemma
\ref{lem:Monster}, line (\ref{line:Monster-hj}), this vanishes for a system with an
involutive symbol.

\item[(b)] For the second subcase choose $i<k<j$.  Then
(\ref{gl:Allgemein2}) becomes
\begin{gather*}
    \sum_{\gamma=\bq{1}{h}+1}^{\bq{1}{j}} C_\gamma^h(\phi^\alpha_i)
    C_\delta^k(\phi^\gamma_j)
\end{gather*}
where $\bq{1}{k}+1 \le \delta \le m$.  According to Lemma \ref{lem:Monster},
line (\ref{line:Monster-hk2}), this vanishes for a system with an involutive symbol.

\item[(c)] For the third subcase choose $k=i$.  Then (\ref{gl:Allgemein2})
  becomes
\begin{gather}\label{gl:Contraction1}
    \sum_{\gamma=\bq{1}{h}+1}^{\bq{1}{j}}
        C_\gamma^h(\phi^\alpha_i)
        C_\delta^i(\phi^\gamma_j)
\end{gather}
where $\bq{1}{k}+1 \le \delta \le m$.  Since we have $h<k=i$, for any
such $\delta$ the cross derivative $u^\delta_{ih}=u^\delta_{hi}$
labels two columns in the complete matrix: the one with label
$\delta$ in $\mbl{[i,h]}{1}{\bq{1}{i}}{}{i}$ and the one with label
$\delta$ in $\mbl{[i,i]}{1}{\bq{1}{i}}{}{h}$ which according to
(\ref{gl:Allgemein2}) has the new entries
\begin{gather}
    -\delta_{ki} C^h_\delta(\phi_j^\alpha)
        + \sum_{\gamma=\bq{1}{k}+1}^{\bq{1}{j}}
            C_\gamma^k(\phi^\alpha_i)
            C_\delta^h(\phi^\gamma_j)
                + \delta_{hj} C^k_\delta(\phi^\alpha_i)
    \label{gl:Allgemein4}
\end{gather}
where $\bq{1}{h}+1 \le \delta \le m$.  In the current subcase,
(\ref{gl:Allgemein4}) becomes
\begin{gather}\label{gl:Spaltung1}
    -C^h_\delta(\phi_j^\alpha)
        + \sum_{\gamma=\bq{1}{i}+1}^{\bq{1}{j}}
            C_\gamma^i(\phi^\alpha_i)
            C_\delta^h(\phi^\gamma_j)
\end{gather}
where $\bq{1}{h}+1 \le \delta \le m$.  Since the columns of both
(\ref{gl:Contraction1}) and (\ref{gl:Spaltung1}) are labelled by the same
second-order derivatives, we have to contract them, which means
adding their new entries.  For all $\bq{1}{h}+1 \le \delta \le
\bq{1}{i}$ this yields
\begin{gather*}
    -C^h_\delta(\phi_j^\alpha)
        + \sum_{\gamma=\bq{1}{i}+1}^{\bq{1}{j}}
            C_\gamma^i(\phi^\alpha_i)
            C_\delta^h(\phi^\gamma_j) ,
\end{gather*}
which, for a system with an involutive symbol, vanishes according to
Lemma \ref{lem:Monster}, line (\ref{line:Monster-hi-a}), and for
all $\bq{1}{i}+1 \le \delta \le m$
\begin{gather*}
    -C^h_\delta(\phi_j^\alpha)
        + \sum_{\gamma=\bq{1}{i}+1}^{\bq{1}{j}}
            C_\gamma^i(\phi^\alpha_i)
            C_\delta^h(\phi^\gamma_j)
        + \sum_{\gamma=\bq{1}{h}+1}^{\bq{1}{j}}
            C_\gamma^h(\phi^\alpha_i)
            C_\delta^i(\phi^\gamma_j) ,
\end{gather*}
which, for a system with an involutive symbol, again vanishes
according to Lemma~\ref{lem:Monster}, line~(\ref{line:Monster-hi-b}).

\item[(d)] For the fourth subcase choose $k<i$.  Under this
assumption, we have to distinguish two further subcases: $k=h$ and
$k<h$.  The subcase $k>h$ need not be considered since
$u^\delta_{kh}=u^\delta_{hk}$ and any second-order derivative is used
only once to label a column in the contracted matrix.

\item[(d1)] First consider $k<i$ and $k = h$.  Then still $h<i$,
furthermore $k<j$, and (\ref{gl:Allgemein2}) becomes
\begin{gather*}
    \sum_{\gamma=\bq{1}{h}+1}^{\bq{1}{j}} C_\gamma^h(\phi^\alpha_i) C_\delta^h(\phi^\gamma_j)
\end{gather*}
where $\bq{1}{h}+1 \le \delta \le m$.  According to Lemma
\ref{lem:Monster}, line (\ref{line:Monster-hh}), this vanishes
for a system with an involutive symbol.

\item[(d2)] At last consider $h<k<i$.  Then $k<j$, and
(\ref{gl:Allgemein2}) becomes
\begin{gather}\label{gl:Spaltung2}
    \sum_{\gamma=\bq{1}{h}+1}^{\bq{1}{j}}
        C_\gamma^h(\phi^\alpha_i)
        C_\delta^k(\phi^\gamma_j)
\end{gather}
where $\bq{1}{k}+1 \le \delta \le m$.  Since we have $h<k<i$, for any
such $\delta$ the cross derivative $u^\delta_{kh}=u^\delta_{hk}$
labels two columns in the complete matrix: the one with label
$\delta$ in $\mbl{[i,h]}{1}{\bq{1}{i}}{}{k}$ and the one with label
$\delta$ in $\mbl{[i,k]}{1}{\bq{1}{i}}{}{h}$ which according to
(\ref{gl:Allgemein2}) has the new entries
\begin{gather}
    -\delta_{ki} C^h_\delta(\phi_j^\alpha)
        + \sum_{\gamma=\bq{1}{k}+1}^{\bq{1}{j}}
            C_\gamma^k(\phi^\alpha_i)
            C_\delta^h(\phi^\gamma_j)
                + \delta_{hj} C^k_\delta(\phi^\alpha_i)
    \label{gl:Allgemein5}
\end{gather}
where $\bq{1}{h}+1 \le \delta \le m$.  In the current subcase,
(\ref{gl:Allgemein5}) becomes
\begin{gather}\label{gl:Contraction2}
    \sum_{\gamma=\bq{1}{k}+1}^{\bq{1}{j}}
        C_\gamma^k(\phi^\alpha_i)
        C_\delta^h(\phi^\gamma_j)
\end{gather}
where $\bq{1}{h}+1 \le \delta \le m$.  Since the columns of both
(\ref{gl:Spaltung2}) and (\ref{gl:Contraction2}) are labelled by the same
second-order derivatives, we have to contract them, which means
adding their new entries.  This yields
\begin{gather*}
    \sum_{\gamma=\bq{1}{k}+1}^{\bq{1}{j}}
        C_\gamma^k(\phi^\alpha_i)
        C_\delta^h(\phi^\gamma_j)
\end{gather*}
for $\bq{1}{h}+1 \le \delta \le \bq{1}{k}$ and vanishes according to Lemma \ref{lem:Monster}, line
(\ref{line:Monster-hk-a}), and
\begin{gather*}
    \sum_{\gamma=\bq{1}{h}+1}^{\bq{1}{j}}
        C_\gamma^h(\phi^\alpha_i)
        C_\delta^k(\phi^\gamma_j)
    +
    \sum_{\gamma=\bq{1}{k}+1}^{\bq{1}{j}}
        C_\gamma^k(\phi^\alpha_i)
        C_\delta^h(\phi^\gamma_j)
\end{gather*}
for $\bq{1}{k}+1 \le \delta \le m$ and vanishes according to Lemma \ref{lem:Monster}, line
(\ref{line:Monster-hk-b}).
\end{enumerate}

What we have shown is that in the contracted matrix all the entries
without some entry $-1$ in a negative unit block to their left may be
eliminated by elementary row transformations if the equation has an
involutive symbol.  Thus under this assumption the rank condition
(\ref{gl:Rangbed}) is satisf\/ied.

Now for the augmented rank condition, equation (\ref{gl:RangbedAugmented}).
To transform the augmented complete matrix into row echelon form, we use for
each $j$ and $i$ the same procedure as for the transformation of the
non-augmented complete matrix, except that now the matrices $(c_{\alpha t})$
and $(d_{st})$ are augmented by one more column each as follows.  Fix $1 < j
\le n$.  Let $1 \le i < j$.  Then, according to the structure of the augmented
complete matrix, for its
transformation into row echelon form we have to eliminate for $1 \le h \le i$
the entries in the blocks $\mbl{[\mathnormal\Xi_i]}{}{}{}{h}$, as we did for
the non-augmented complete matrix given in equation (\ref{def:CompleteMatrix}).  We have to consider the ef\/fect of these
transformations on the additional entries which make up the rightmost column
in the augmented complete matrix.  These are $-\mathnormal\Theta^1_{ij}$,
$-\mathnormal\Theta^2_{ij}$, \ldots, $-\mathnormal\Theta^m_{ij}$, given in
equation (\ref{gl:VektorTheta}).  Of these entries, only
$-\mathnormal\Theta^1_{ij}$, $-\mathnormal\Theta^2_{ij}$, \ldots,
$-\mathnormal\Theta^{\bq{1}{i}}_{ij}$ are af\/fected (since we eliminate the
entries which are in rows $1$ to $\bq{1}{i}$ of the matrices
$\mbl{[\mathnormal\Xi_i]}{}{}{}{h}$).  We add them as the rightmost column in
the augmented matrix~$(c_{\alpha t})$.  Now f\/ix $1 \le h \le i$.  Then augment
the matrix $(d_{st})$, used in the process of eliminating the entries in
$\mbl{[\mathnormal\Xi_i]}{}{}{}{h}$, by adding the entries
$-\mathnormal\Theta^{\bq{1}{h}+1}_{hj}$,
$-\mathnormal\Theta^{\bq{1}{h}+2}_{hj}$, \ldots, $-\mathnormal\Theta^m_{hj}$
as its rightmost column, in accordance with the structure of the augmented
complete matrix.  Then the
substitution~(\ref{lem:substitution}) yields as the transformed entries
$c_{\alpha t}+\sum_{s=1}^S b_{\alpha s}d_{s t}$ the same as for the
transformation of the non-augmented matrix except, of course, for the new last
column.  Now let denote $t$ the index of this last column.  Fix some row index
$1 \le \alpha \le \bq{1}{i}$.  Then the entry $c_{\alpha
  t}=-\mathnormal\Theta^\alpha_{ij}$ transforms as follows: For all $1 \le h
\le i$ we have to add to it
\begin{gather*}
    \sum_{s=1}^S b_{\alpha s}d_{st}
    =
    \sum_{\gamma=\bq{1}{h}+1}^m b_{\alpha \gamma}d_{\gamma t}.
\end{gather*}
Since here the matrix $(b_{\alpha \gamma})=\mbl{[\mathnormal\Xi_i]}{}{}{}{h}$,
this is the product of the row with index $\alpha$ in the matrix
$\mbl{[\mathnormal\Xi_i]}{}{}{}{h}$ and the transpose of the vector
$(-\mathnormal\Theta^{\bq{1}{h}+1}_{hj},
-\mathnormal\Theta^{\bq{1}{h}+2}_{hj},\ldots, -\mathnormal\Theta^m_{hj})$.
According to the structure of~$\mbl{[\mathnormal\Xi_i]}{}{}{}{h}$ as def\/ined
in equation (\ref{def:Xi-oben}), this equals
\begin{gather*}
    \sum_{\gamma=\bq{1}{h}+1}^m C^h_\gamma(\phi^\alpha_i) \mathnormal\Theta^\gamma_{hj}
   .
\end{gather*}
The entries $\mathnormal\Theta^\gamma_{hj}$ can be taken from
equation (\ref{gl:VektorTheta}).  Since $\bq{1}{h}+1 \le \gamma \le
m$, we have $\mathnormal\Theta^\gamma_{hj}=C_h^{(1)}(\phi^\alpha_j)$
if $\bq{1}{h}+1 \le \gamma \le \bq{1}{j}$ and
$\mathnormal\Theta^\gamma_{hj}=0$ if $\bq{1}{j}+1 \le \gamma \le m$.
Therefore we get
\begin{gather*}
    \sum_{\gamma=\bq{1}{h}+1}^{\bq{1}{j}}
        C^h_\gamma(\phi^\alpha_i)
            C^{(1)}_h(\phi^\alpha_j)
\end{gather*}
as the summand for each $h$ to be added to $c_{\alpha
t}=-\mathnormal\Theta^\alpha_{ij}$.  Since $1 \le \alpha \le
\bq{1}{i}$, we have
$-\mathnormal\Theta^\alpha_{ij}=-C_i^{(1)}(\phi^\alpha_j)+C_j^{(1)}(\phi^\alpha_i)$.
Therefore the entry $c_{\alpha t}$ transforms into
\begin{gather*}
    -C_i^{(1)}(\phi^\alpha_j)
    +C_j^{(1)}(\phi^\alpha_i)
    +
    \sum_{h=1}^i
        \sum_{\gamma=\bq{1}{h}+1}^{\bq{1}{j}}
            C^h_\gamma(\phi^\alpha_i)
                C^{(1)}_h(\phi^\gamma_j)
     ,
\end{gather*}
which is the integrability condition in line (\ref{line:Monster-IntegBed}) of
Lemma~\ref{lem:Monster}, except for the sign.  The dif\/ference in sign comes
from the fact that the augmented complete matrix describes the system of
equations~(\ref{gl:Schritti}), where the entries from the row with index
$\alpha$ in the matrix $\mathnormal\Xi_i$ are on the opposite side from the
entries from the row with index $\alpha$ in the matrix $\mathnormal\Xi_j$ and
the inhomogeneous term~$-\mathnormal\Theta^\alpha_{ij}$, while in Lemma
\ref{lem:Monster} for each $i$, $j$ and $\alpha$ the corresponding equation is
set to zero, if the dif\/ferential equation is involutive.

This means that the augmented rank condition holds for all $1 \le i <
j \le n$ if, and only if, the equation has an involutive symbol and
is formally integrable.

\section{Proof of Proposition \ref{prop:DiffBedInvolutiv}}\label{anhang:DiffBedInvolutiv}

We are going to show that the system (\ref{gl:DifferentielleBed}) is
involutive.  Consider the basis $(U_i \colon 1 \le i \le n)$ of the
distribution $\U$ given by $U_i = X_i + \zeta_i^pY_p$ with yet undetermined
coef\/f\/icient functions $\zeta_i^p \in \mathcal{F}(\Rc{1})$.  For all $1 \le i
\le n$ and $1 \le p \le r$, the independent variables of the functions
$\zeta^p_i$ are the coordinates on $\Rc{1}$, which are $\xv,\uv$ and all
$u^\alpha_h$ such that $(\alpha,h) \not\in \mathcal{B}$.  To apply a vector
f\/ield $U_j = \partial_{x^j}+\cdots$ to a function $\zeta^p_i$ includes a
derivation with respect to $x^j$.  We order the independent variables such
that if $j>i$, then $x^j$ is greater than $x^i$, and each $x^i$ is greater
than all the variables $u^\alpha$ and $u^\alpha_h$ where $(\alpha,h) \not\in
\mathcal{B}$.  For any equation $H^p_{ij}$ within the system
(\ref{gl:DifferentielleBed}), the application of the vector f\/ield $U_j =
\partial_{x^j}+\cdots$ to $\zeta_i^p$ yields
$\partial\zeta_i^p/\partial x^j$ as the leader of that
equation; therefore equation $H^p_{ij}$ is of class $j$, and the
equations of maximal class are $H_{in}^p$; the equations of second
highest class in the system are $H_{i n-1}^p$ and so on.  There are
only equations $H_{i j}^p$ of a class indicated by some index $2 \le
j \le n$.

From the Jacobi identity for vector f\/ields
$U_i$, $U_j$ and $U_k$ where $1 \le i<j<k \le n$, we have
\begin{gather*}
    [U_i,[U_j,U_k]]+[U_j,[U_k,U_i]]+[U_k,[U_i,U_j]] = 0 .
\end{gather*}
The structure equations (\ref{gl:UU}) for the vector f\/ields $U_h$ and
the def\/initions of $G_{ij}^c$ and $H_{ij}^p$ in equations
(\ref{gl:AlgebraischeBed}) and (\ref{gl:DifferentielleBed}) imply
that this is
\begin{gather*}
    0=
    [U_i,G_{jk}^cZ_c+H_{jk}^p Y_p]
      +
    [U_j,G_{ki}^cZ_c+H_{ki}^p Y_p]
      +
    [U_k,G_{ij}^cZ_c+H_{ij}^p Y_p]
    \\
\phantom{0}{}     =
    G_{jk}^c[U_i,Z_c]
    +
    U_i(G_{jk}^c)Z_c
    +
    H_{jk}^p[U_i,Y_p]
    +
    U_i(H_{jk}^p)Y_p
    \\
 \phantom{0=}{}   +
    G_{ki}^c[U_j,Z_c]
    +
    U_j(G_{ki}^c)Z_c
    +
    H_{ki}^p[U_j,Y_p]
    +
    U_j(H_{ki}^p)Y_p
    \\
 \phantom{0=}{}    +
    G_{ij}^c[U_k,Z_c]
    +
    U_k(G_{ij}^c)Z_c
    +
    H_{ij}^p[U_k,Y_p]
    +
    U_k(H_{ij}^p)Y_p
 .
\end{gather*}
The combined system (\ref{gl:AlgebraischeBed},
\ref{gl:DifferentielleBed}) means that all $G^c_{ab}=0$ and all
$H^p_{ab}=0$ which implies that $\U$ is involutive, which it is,
being in triangular form, exactly if all $[U_a,U_b]=0$.  This leaves
only
\begin{gather*}
    0  =
	\{U_i(G_{jk}^c)+U_j(G_{ki}^c)+U_k(G_{ij}^c)\} Z_c
     +
    \{U_i(H_{jk}^p)+U_j(H_{ki}^p)+U_k(H_{ij}^p)\} Y_p
     .
\end{gather*}
As part of a basis for $\V^\prime[\Rc{q}]$, the vector f\/ields
$Z_c$ and $Y_p$ are linearly independent, which means their
coef\/f\/icients must vanish individually.  So in particular
\begin{displaymath}
    U_i(H_{jk}^p)+U_j(H_{ki}^p)+U_k(H_{ij}^p) = 0  .
\end{displaymath}
Under the assumption $i<j<k$, the term $U_k(H_{ij}^p)$ contains
derivations with respect to $x_k$ of $U_i(\zeta_j^p)$ and
$U_j(\zeta_i^p)$.  Thus, according to our order, this is a
non-multiplicative prolongation, and the remaining terms are
multiplicative prolongations.  But since any non-multiplicative
prolongation within the system~(\ref{gl:DifferentielleBed}) must be
of such a form, it is a linear combination of multiplicative
prolongations.  Therefore, no integrability conditions arise from
cross-derivatives (and none arise from a
prolongation of lower order equations since all equations of the
system are of f\/irst order).

If we set $\partial\zeta_i^p/\partial
x^j=:(\zeta_i^p)_j$ for the leaders and
\begin{displaymath}
    \tilde{U}_j(\zeta_i^p)
    :=
    U_j(\zeta_i^p)-(\zeta_i^p)_j
\end{displaymath}
and solve each equation of the system (\ref{gl:DifferentielleBed})
for its leader, then it takes the form
\begin{gather}
    (\zeta^p_1)_n  = U_1(\zeta^p_n)-\tilde{U}_n(\zeta_1^p)   , \nonumber\\
    (\zeta^p_2)_n  = U_2(\zeta^p_n)-\tilde{U}_n(\zeta_2^p)   , \nonumber\\
    \cdots \cdots\cdots\cdots\cdots\cdots\cdots\cdots\nonumber\\
    (\zeta^p_{n-1})_n = U_{n-1}(\zeta^p_n)-\tilde{U}_n(\zeta_{n-1}^p)   ,\nonumber \\
        (\zeta^p_1)_{n-1} = U_1(\zeta^p_{n-1})-\tilde{U}_{n-1}(\zeta_1^p)   ,\nonumber \\
    (\zeta^p_2)_{n-1} = U_2(\zeta^p_{n-1})-\tilde{U}_{n-1}(\zeta_2^p)   , 
 \tag{\ref{gl:DifferentielleBed}$^*$}\label{gl:DifferentielleBed*}
    \\
    \cdots \cdots\cdots\cdots\cdots\cdots\cdots\cdots \nonumber\\
    (\zeta^p_{n-2})_{n-1} = U_{n-2}(\zeta^p_{n-1})-\tilde{U}_{n-1}(\zeta_{n-2}^p)   , \nonumber\\
       \cdots \cdots\cdots\cdots\cdots\cdots\cdots\cdots\nonumber\\
        (\zeta^p_2)_3 = U_2(\zeta^p_3)-\tilde{U}_3(\zeta_2^p)  , \nonumber \\
    (\zeta^p_1)_3 = U_1(\zeta^p_3)-\tilde{U}_3(\zeta_1^p)  , \nonumber\\
    (\zeta^p_1)_2 = U_1(\zeta^p_2)-\tilde{U}_2(\zeta_1^p)  ;\nonumber
  \end{gather} 
here for each line $1 \le p \le r$.  Therefore the system
(\ref{gl:DifferentielleBed*}) is in Cartan normal form given in
Def\/inition \ref{def:CNF}.

Now one can prove (see \cite[Lemma 2.4.29]{Fesser:2008}) that such a
dif\/ferential equation $\Rc{1}$ is involutive if, and only if, all
non-multiplicative prolongations of the equations
(\ref{gl:CNF-a})--(\ref{gl:CNF-c}) and all formal derivatives with respect to
all the $x^i$ of the algebraic equations (\ref{gl:CNF-d}) are dependent on the
equations of the system (\ref{gl:CNF}) and its multiplicative prolongations
only.  Therefore, it follows that the system (\ref{gl:DifferentielleBed}) is
involutive.

\section{Proof of Theorem~\ref{thm:CombinedSolvable}}\label{anhang:ExThmFlatVConn}

We consider the combined system of algebraic and dif\/ferential conditions
(\ref{gl:AlgebraischeBed}), (\ref{gl:DifferentielleBed}) and want to show that
it has a solution.  We follow the strategy outlined above and eliminate some
of the unknowns~$\zeta^p_\ell$.  As we consider each of the equations of
(\ref{gl:DifferentielleBed}) as being solved for its derivative
$\partial\zeta^{(\beta,h)}_i/\partial x^j$ of highest class $j$, as given in
equation (\ref{gl:DifferentielleBed*}), we must take a closer look only at
those equations where this leading derivative is of one of the unknowns we
eliminate.  The structure of the vectors $\zeta_i$, given in equations
(\ref{gl:zeta-Vektor-Beziehungen-1}) and (\ref{gl:zeta-rel}), shows which ones
these are.  Let $k$ be such that $2 \le k \le n$.  Then for the subsystem of
the equations of class $k$ in the system (\ref{gl:DifferentielleBed}), the
equations which hold the following terms are concerned:
\begin{gather*}
     U_k\big(\zeta_2^{(\beta,1)}\big), \\
     U_k\big(\zeta_3^{(\beta,1)}\big), \ U_k\big(\zeta_3^{(\beta,2)}\big), \\
   \cdots \cdots  \cdots \cdots\cdots \cdots \cdots \cdots\cdots \cdots\\
     U_k\big(\zeta_{k-1}^{(\beta,1)}\big), \ U_k\big(\zeta_{k-1}^{(\beta,2)}\big), \ \ldots,
      \ U_k\big(\zeta_{k-1}^{(\beta,k-2)}\big)  ;
\end{gather*}
here, for any $U_k(\zeta_i^{(\beta,h)})$, we have $\bq{1}{h}+1 \le
\beta \le m$.  We now show that these equations vanish.  The proof is
by straightforward calculation, though tedious and requiring a case
distinction.  Let $1<i<k$.  Fix some $\bq{1}{h}+1 \le \beta \le m$.
Consider the equation
\begin{equation}\label{gl:Uikh=Ukih}
    U_i\big(\zeta_k^{(\beta,h)}\big) = U_k\big(\zeta_i^{(\beta,h)}\big) .
\end{equation}
Then $h<i<k$.  According to the structure of the vector $\zeta_i$, the
entries of which in its $h$th block are of two kinds, there are two
cases.
\begin{enumerate}\itemsep=0pt
\item The interrelation for $\zeta_i^{(\beta,h)}$ is an equality:
$\zeta_i^{(\beta,h)}=\zeta_h^{(\beta,i)}$.  This is so if, and only
if, $\bq{1}{i}+1 \le \beta \le m$ according to the structure of
$\zeta_i$.  Now there arise two subcases.
\begin{enumerate}\itemsep=0pt
\item The other interrelation is an equality, too:
$\zeta_k^{(\beta,h)}=\zeta_h^{(\beta,k)}$.  This is so if, and only
if, $\bq{1}{k}+1 \le \beta \le m$ according to the structure of
$\zeta_k$.  In this subcase, equation (\ref{gl:Uikh=Ukih}) becomes
\begin{equation}\label{gl:Uihk=Ukhi}
    U_i\big(\zeta_h^{(\beta,k)}\big) = U_k\big(\zeta_h^{(\beta,i)}\big) .
\end{equation}
Since the system (\ref{gl:DifferentielleBed}) contains the equalities
$U_i\big(\zeta_h^{(\beta,k)}\big)=U_h\big(\zeta_i^{(\beta,k)}\big)$ and
$U_k\big(\zeta_h^{(\beta,i)}\big)=U_h\big(\zeta_k^{(\beta,i)}\big)$, equation
(\ref{gl:Uihk=Ukhi}) becomes
\begin{equation*}
    U_h\big(\zeta_i^{(\beta,k)}\big) = U_h\big(\zeta_k^{(\beta,i)}\big)  .
\end{equation*}
Since $i<k$ and $\bq{1}{k}+1 \le \beta \le m$, from the structure of
$\zeta_k$ follows $\zeta_i^{(\beta,k)}=\zeta_k^{(\beta,i)}$.  Thus,
equation (\ref{gl:Uikh=Ukih}) vanishes.

\item The other interrelation is an af\/f\/ine-linear combination:
\begin{displaymath}
    \zeta_k^{(\beta,h)}
    =
    \sum_{a=1}^k
        \sum_{\gamma=\bq{1}{a}+1}^m
            C^a_\gamma\big(\phi^\beta_k\big)\zeta_h^{(\gamma,a)}
    +
    C^{(1)}_h\big(\phi^\beta_k\big)  .
\end{displaymath}
This is so if, and only if, $\bq{1}{i}+1 \le \beta \le \bq{1}{k}$
according to the structure of $\zeta_k$.  In this subcase, the term
$U_i(\zeta_k^{(\beta,h)})$ in equation (\ref{gl:Uikh=Ukih}) becomes
\begin{subequations}\label{gl:Uikh=Uihk}
\begin{gather}
    U_i\big(\zeta_k^{(\beta,h)}\big)
    =
    \sum_{a=1}^k
        \sum_{\gamma=\bq{1}{a}+1}^m
            C_\gamma^a(\phi_k^\beta)
                U_i\big(\zeta_h^{(\gamma,a)}\big)
    \\
    \phantom{U_i\big(\zeta_k^{(\beta,h)}\big)    =}{} +
    \sum_{a=1}^k
        \sum_{\gamma=\bq{1}{a}+1}^m
            U_i\big(C_\gamma^a\big(\phi_k^\beta\big)\big)
                \zeta_h^{(\gamma,a)}
    +
    U_i\big(C^{(1)}_h\big(\phi_k^\beta\big)\big)  .
\end{gather}
\end{subequations}
The term $U_k(\zeta_i^{(\beta,h)})$ in equation (\ref{gl:Uikh=Ukih})
becomes
\begin{subequations}\label{gl:Ukih=Ukhi}
\begin{gather}
    U_k\big(\zeta_i^{(\beta,h)}\big)
     =
    U_k\big(\zeta_h^{(\beta,i)}\big)
   =
    U_h\big(\zeta_k^{(\beta,i)}\big)
    \notag
    \\
  \phantom{U_k\big(\zeta_i^{(\beta,h)}\big)}{}  =
    U_h
    \Bigg(\sum_{a=1}^k
        \sum_{\gamma=\bq{1}{a}+1}^m
            C_\gamma^a\big(\phi_k^\beta\big)
                \zeta_i^{(\gamma,a)}
    +C^{(1)}_i\big(\phi_k^\beta\big)
    \Bigg)
    \notag
    \\
\phantom{U_k\big(\zeta_i^{(\beta,h)}\big)}{}    =
    \sum_{a=1}^k
        \sum_{\gamma=\bq{1}{a}+1}^m
            C_\gamma^a\big(\phi_k^\beta\big)
                U_h\big(\zeta_i^{(\gamma,a)}\big)
    \\
\phantom{U_k\big(\zeta_i^{(\beta,h)}\big)=}{}    +
    \sum_{a=1}^k
        \sum_{\gamma=\bq{1}{a}+1}^m
            U_h\big(C_\gamma^a\big(\phi_k^\beta\big)\big)
                \zeta_i^{(\gamma,a)}
    +
    U_h\big(C^{(1)}_i\big(\phi_k^\beta\big)\big)  ;
\end{gather}
\end{subequations}
here we have the f\/irst equality because we are considering the f\/irst
main case, the second equality because of the structure of the system
(\ref{gl:DifferentielleBed}) and the third equality according to the
structure of $\zeta_k$, since $i<k$ and because $\bq{1}{i}+1 \le
\beta \le \bq{1}{k}$.  Substituting (\ref{gl:Uikh=Uihk}) and
(\ref{gl:Ukih=Ukhi}) in equation (\ref{gl:Uikh=Ukih}) and factoring
out, we get
\begin{subequations}
\begin{gather}
    0  =
        \sum_{a=1}^k
            \sum_{\gamma=\bq{1}{a}+1}^m
                C_\gamma^a(\phi_k^\beta)
                    \big\{
                        U_i\big(\zeta_h^{(\gamma,a)}\big)
                        -
                        U_h\big(\zeta_i^{(\gamma,a)}\big)
                    \big\}
    \label{gl:Uikh=Ukih-1a}
    \\
    \phantom{0  =}{} +
        \sum_{a=1}^k
            \sum_{\gamma=\bq{1}{a}+1}^m
                \big\{%
                    U_i\big(C_\gamma^a\big(\phi_k^\beta\big)\big)
                        \zeta_h^{(\gamma,a)}
                    -
                    U_h\big(C_\gamma^a\big(\phi_k^\beta\big)\big)
                        \zeta_i^{(\gamma,a)}%
                \big\}
    \label{gl:Uikh=Ukih-1b}
    \\
   \phantom{0  =}{}+
    U_i\big(C^{(1)}_h\big(\phi_k^\beta\big)\big)
    -
    U_h\big(C^{(1)}_i\big(\phi_k^\beta\big)\big) .
    \label{gl:Uikh=Ukih-1c}
\end{gather}
\end{subequations}
Line (\ref{gl:Uikh=Ukih-1c}) contains the Lie bracket
$[U_i,C^{(1)}_h](\phi_k^\beta)$.  According to the structure of the system
(\ref{gl:DifferentielleBed}), the term (\ref{gl:Uikh=Ukih-1a}) vanishes.  If
the terms (\ref{gl:Uikh=Ukih-1b}) and (\ref{gl:Uikh=Ukih-1c}) vanish, too,
then so does equation (\ref{gl:DifferentielleBed}).  Otherwise they form a new
algebraic condition for (\ref{gl:DifferentielleBed}), which can be solved for
some function $\zeta_h^{(\beta,a)}$.  Substituting this function in
(\ref{gl:DifferentielleBed}) does not change the classes or the numbers of the
single equations therein.  Thus, equation~(\ref{gl:Uikh=Ukih}) vanishes.
\end{enumerate}

\item The interrelation for $\zeta_i^{(\beta,h)}$ is an af\/f\/ine-linear
combination:
\begin{displaymath}
    \zeta_i^{(\beta,h)}
    =
    \sum_{a=1}^i
        \sum_{\gamma=\bq{1}{a}+1}^m
            C^a_\gamma\big(\phi^\beta_i\big)\zeta_h^{(\gamma,a)}
            +
            C^{(1)}_h\big(\phi^\beta_i\big).
\end{displaymath}
This is so if, and only if, $\bq{1}{h}+1 \le \beta \le \bq{1}{i}$
according to the structure of $\zeta_i$.  Since we have $h<i<k$ and
$\bq{1}{i} \le \bq{1}{k}$, according to the structure of $\zeta_k$
the other interrelation is an af\/f\/ine-linear combination, too:
\begin{displaymath}
    \zeta_k^{(\beta,h)}
    =
    \sum_{b=1}^i
        \sum_{\delta=\bq{1}{b}+1}^m
            C^b_\delta\big(\phi^\beta_k\big)\zeta_h^{(\delta,b)}
            +
            C^{(1)}_h\big(\phi^\beta_k\big) .
\end{displaymath}
Thus, equation (\ref{gl:Uikh=Ukih}) becomes
\begin{subequations}\label{gl:Uikh=Ukih-2}
\begin{gather}
    0
     =
    \sum_{a=1}^i
        \sum_{\gamma=\bq{1}{a}+1}^m
            C_\gamma^a\big(\phi_i^\beta\big)
                U_k\big(\zeta_h^{(\gamma,a)}\big)
    -
    \sum_{b=1}^k
        \sum_{\delta=\bq{1}{b}+1}^m
            C_\delta^b\big(\phi_k^\beta\big)
                U_i\big(\zeta_h^{(\delta,b)}\big)
    \label{gl:Uikh=Ukih-2a}
    \\
   \phantom{0=}{} +
    \sum_{a=1}^i
        \sum_{\gamma=\bq{1}{a}+1}^m
            U_k\big(C_\gamma^a\big(\phi_i^\beta\big)\big)
                \zeta_h^{(\gamma,a)}
    -
    \sum_{b=1}^k
        \sum_{\delta=\bq{1}{b}+1}^m
            U_i\big(C_\delta^b\big(\phi_k^\beta\big)\big)
                \zeta_h^{(\delta,b)}
    \label{gl:Uikh=Ukih-2b}
    \\
 \phantom{0=}{}+
    U_k\big(C^{(1)}_h\big(\phi_i^\beta\big)\big)
    -
    U_i\big(C^{(1)}_h\big(\phi_k^\beta\big)\big)  .
    \label{gl:Uikh=Ukih-2c}
\end{gather}
\end{subequations}
In part (\ref{gl:Uikh=Ukih-2a}), the terms
$U_k(\zeta_h^{(\gamma,a)})$ and $U_i(\zeta_h^{(\delta,b)})$ are equal
to $U_h(\zeta_k^{(\gamma,a)})$ and $U_h(\zeta_i^{(\delta,b)})$
according to the structure of the system
(\ref{gl:DifferentielleBed}).  Thus, equation (\ref{gl:Uikh=Ukih-2})
becomes
\begin{gather}
    0
    =
        \sum_{a=1}^i
            \sum_{\gamma=\bq{1}{a}+1}^m
                C_\gamma^a\big(\phi_i^\beta\big)
                    U_h\big(\zeta_k^{(\gamma,a)}\big)
    -
        \sum_{b=1}^k
            \sum_{\delta=\bq{1}{b}+1}^m
                C_\delta^b\big(\phi_k^\beta\big)
                    U_h\big(\zeta_i^{(\delta,b)}\big)
    \tag{\ref{gl:Uikh=Ukih-2a}$^\prime$}\label{gl:Uikh=Ukih-3}
    \\
 \phantom{0=}{}+ \text{(\ref{gl:Uikh=Ukih-2b})} + \text{(\ref{gl:Uikh=Ukih-2c})} .
    \notag
\end{gather}
Factoring out the vector f\/ield $U_h$ in part (\ref{gl:Uikh=Ukih-3}),
this equals
\begin{subequations}\label{gl:Uikh=Ukih-4}
\begin{gather}
    0
     =
        U_h
        \Bigg(
        \sum_{a=1}^i
            \sum_{\gamma=\bq{1}{a}+1}^m
                C_\gamma^a\big(\phi_i^\beta\big)
                    \zeta_k^{(\gamma,a)}
    -
        \sum_{b=1}^k
            \sum_{\delta=\bq{1}{b}+1}^m
                C_\delta^b\big(\phi_k^\beta\big)
                    \zeta_i^{(\delta,b)}
        \Bigg)
    \label{gl:Uikh=Ukih-4a}
    \\
 \phantom{0=}{} -
    \Bigg(
    \sum_{a=1}^i
            \sum_{\gamma=\bq{1}{a}+1}^m
                U_h\big(C_\gamma^a\big(\phi_i^\beta\big)\big)
                    \zeta_k^{(\gamma,a)}
    -
        \sum_{b=1}^k
            \sum_{\delta=\bq{1}{b}+1}^m
                U_h\big(C_\delta^b\big(\phi_k^\beta\big)\big)
                    \zeta_i^{(\delta,b)}
    \Bigg)
    \label{gl:Uikh=Ukih-4b}
    \\
 \phantom{0=}{}    + (\ref{gl:Uikh=Ukih-2}\text{b}) + (\ref{gl:Uikh=Ukih-2}\text{c})  .
    \label{gl:Uikh=Ukih-4c}
\end{gather}
\end{subequations}
Now one can show (see \cite[Corollary 3.3.25 (for $j=k$)]{Fesser:2008}) that the
term (\ref{gl:Uikh=Ukih-4a}) equals
\begin{displaymath}
    U_h\big(C^{(1)}_i\big(\phi^\beta_k\big)-C^{(1)}_k\big(\phi^\beta_i\big)\big) ,
\end{displaymath}
which does not contain any $\zeta_k^{(\gamma,a)}$ or
$\zeta_i^{(\delta,b)}$ any more; it is an algebraic expression
instead of the dif\/ferential expression that it seems to be when
written in the form (\ref{gl:Uikh=Ukih-4a}).  The other terms,
(\ref{gl:Uikh=Ukih-4b}) and (\ref{gl:Uikh=Ukih-4c}), are algebraic,
too.  So all of equation (\ref{gl:Uikh=Ukih}) has shown to be an
algebraic condition when the interrelations between the entries of
the vectors~$\zeta_h$,~$\zeta_i$ and~$\zeta_k$, as noted in equations
(\ref{gl:zeta-Vektor-Beziehungen-1}) and (\ref{gl:zeta-rel}), are taken into account.

If this new algebraic condition for the system
(\ref{gl:DifferentielleBed}) vanishes, equation (\ref{gl:Uikh=Ukih})
vanishes.  Otherwise, this new algebraic condition given in equation
(\ref{gl:Uikh=Ukih-4}) now appears as
\begin{subequations}\label{gl:Uikh=Ukih-5}
\begin{gather}
    0
     =
    U_h%
    \big(
        C^{(1)}_i\big(\phi^\beta_k\big)
        -
        C^{(1)}_k\big(\phi^\beta_i\big)
    \big)
    \tag{\ref{gl:Uikh=Ukih-4a}$^\prime$}
    \\
 \phantom{0=}{} -
    \Bigg(
        \sum_{a=1}^i
            \sum_{\gamma=\bq{1}{a}+1}^m
                U_h\big(C_\gamma^a\big(\phi_i^\beta\big)\big)
                    \zeta_k^{(\gamma,a)}
        -
        \sum_{b=1}^k
            \sum_{\delta=\bq{1}{b}+1}^m
                U_h\big(C_\delta^b\big(\phi_k^\beta\big)\big)
                    \zeta_i^{(\delta,b)}
    \Bigg)
    \tag{\ref{gl:Uikh=Ukih-4b}}
    \\
 \phantom{0=}{}+
        \sum_{a=1}^i
            \sum_{\gamma=\bq{1}{a}+1}^m
                U_k\big(C_\gamma^a\big(\phi_i^\beta\big)\big)
                    \zeta_h^{(\gamma,a)}
        -
        \sum_{b=1}^k
            \sum_{\delta=\bq{1}{b}+1}^m
                U_i\big(C_\delta^b\big(\phi_k^\beta\big)\big)
                    \zeta_h^{(\delta,b)}
    \tag{\ref{gl:Uikh=Ukih-2b}}
    \\
 \phantom{0=}{}+
    U_k\big(C^{(1)}_h\big(\phi_i^\beta\big)\big)
    -
    U_i\big(C^{(1)}_h\big(\phi_k^\beta\big)\big)
  .
    \tag{\ref{gl:Uikh=Ukih-2c}}
\end{gather}
\end{subequations}
Collecting terms in lines (\ref{gl:Uikh=Ukih-4a}$^\prime$) and
(\ref{gl:Uikh=Ukih-2c}), this yields
\begin{gather*}
    0
    =
    \text{(\ref{gl:Uikh=Ukih-2b})}
    +
    \text{(\ref{gl:Uikh=Ukih-4b})}
    \\
\phantom{0=}{} +
    U_h\big(C^{(1)}_i\big(\phi^\beta_k\big)\big)
    -
    U_i\big(C^{(1)}_h\big(\phi_k^\beta\big)\big)
    +
    U_k\big(C^{(1)}_h\big(\phi_i^\beta\big)\big)
    -
    U_h\big(C^{(1)}_k\big(\phi^\beta_i\big)\big)
    .
\end{gather*}
The lower line contains the Lie brackets
$[U_h,C^{(1)}_i](\phi_k^\beta)$ and $[U_k,C^{(1)}_h](\phi_i^\beta)$.
There must be some non-vanishing summand containing a factor
$\zeta_k^{(\gamma,a)}$, $\zeta_i^{(\delta,b)}$,
$\zeta_h^{(\gamma,a)}$ or $\zeta_h^{(\delta,a)}$.  As we did in case
1.\,(b), we solve (\ref{gl:Uikh=Ukih-2}) for this non-vanishing
factor and substitute it into the system~(\ref{gl:DifferentielleBed}), which does not change the class of any
equation therein.  Therefore equation~(\ref{gl:Uikh=Ukih}) drops out
from the system (\ref{gl:DifferentielleBed}).
\end{enumerate}

Now we have shown that all those equations vanish where the leading derivative
is subject to being substituted through the interrelations concerning the
coef\/f\/icient function $\zeta_k^{(\beta,i)}$.  In the system
(\ref{gl:DifferentielleBed*}), these are the equations with the leaders
\begin{gather*}
     \big(\zeta_2^{(\beta,1)}\big)_k  , \\
     \big(\zeta_3^{(\beta,1)}\big)_k, \ \big(\zeta_3^{(\beta,2)}\big)_k  , \\
     \cdots  \cdots \cdots \cdots \cdots \cdots \cdots \cdots \cdots\\
     \big(\zeta_{k-1}^{(\beta,1)}\big)_k, \ \big(\zeta_{k-1}^{(\beta,2)}\big)_k, \ \ldots, \ \big(\zeta_{k-1}^{(\beta,k-2)}\big)_k  ;
\end{gather*}
here $2 \le k \le n$ and $\bq{1}{h}+1 \le \beta \le m$.  The remaining
equations still form an involutive system (we may enumerate the remaining
$\zeta^p_i$ in such a way that no gaps appear) as the considerations for the
system (\ref{gl:DifferentielleBed}) in Proposition~\ref{prop:DiffBedInvolutiv}
apply likewise.  Thus we eventually arrive at an analytic involutive
dif\/ferential equation for the coef\/f\/icient functions $\zeta^k_i$ which is
solvable according to the Cartan--K{\"a}hler theorem (Theorem~\ref{thm:Cartan-Kaehler}).

\subsection*{Acknowledgements}

This work received partial f\/inancial support by the European NEST-Adventure
grant 5006, \emph{Global Integrability of Field Theories}.

\addcontentsline{toc}{section}{References}
\LastPageEnding

\end{document}